\documentclass{article}

\usepackage{proof}
\usepackage{latexsym}

\newcommand{\blem}{\begin{lemma}}
\newcommand{\elem}{\end{lemma}}
\newcommand{\bth}{\begin{theorem}}
\newcommand{\eth}{\end{theorem}}
\newcommand{\benu}{\begin{enumerate}}
\newcommand{\eenu}{\end{enumerate}}
\newcommand{\bdes}{\begin{description}}
\newcommand{\edes}{\end{description}}
\newcommand{\bdf}{\begin{definition}}
\newcommand{\edf}{\end{definition}}
\newcommand{\bcor}{\begin{cor}}
\newcommand{\ecor}{\end{cor}}
\newcommand{\bprp}{\begin{proposition}}
\newcommand{\eprp}{\end{proposition}}
\newcommand{\bmlem}{\begin{mlemma}}
\newcommand{\emlem}{\end{mlemma}}
\newcommand{\bclm}{\begin{claim}}
\newcommand{\eclm}{\end{claim}}
\newcommand{\bprf}{{\bf Proof}.\hspace{2mm}}

\newcommand{\eprf}{\hspace*{\fill} $\Box$}

\newcommand{\ovl}{\overline}
\newcommand{\beqn}{\begin{equation}}
\newcommand{\eeqn}{\end{equation}}
\newcommand{\beqnarr}{\begin{eqnarray}}
\newcommand{\eeqnarr}{\end{eqnarray}}
\newcommand{\beqnarrs}{\begin{eqnarray*}}
\newcommand{\eeqnarrs}{\end{eqnarray*}}

\newcommand{\smlskp}{\\ \smallskip \\}
\newcommand{\spand}{\,\&\,}

\newtheorem{theorem}{Theorem}[section]
\newtheorem{definition}{Definition}[section]
\newtheorem{proposition}{Proposition}[section]
\newtheorem{lemma}{Lemma}[section]
\newtheorem{cor}{Corollary}[section]
\newtheorem{mlemma}{Main Lemma}[section]
\newtheorem{claim}{Claim}[section]

\newcommand{\alp}{\alpha}
\newcommand{\veps}{\varepsilon}
\newcommand{\del}{\delta}
\newcommand{\Del}{\Delta}
\newcommand{\ome}{\omega}
\newcommand{\Ome}{\Omega}
\newcommand{\bet}{\beta}
\newcommand{\gam}{\gamma}
\newcommand{\Gam}{\Gamma}
\newcommand{\kap}{\kappa}
\newcommand{\sig}{\sigma}
\newcommand{\Sig}{\Sigma}

\newcommand{\lam}{\lambda}
\newcommand{\Lam}{\Lambda}
\newcommand{\vphi}{\varphi}

\newcommand{\fal}{\forall}
\newcommand{\exi}{\exists}
\newcommand{\Rarw }{\Rightarrow}

\newcommand{\lrarw}{\leftrightarrow}
\newcommand{\Lrarw}{\Leftrightarrow}
\newcommand{\uarw}{\uparrow}
\newcommand{\darw}{\downarrow}

\newcommand{\cala}{{\cal A}}
\newcommand{\calb}{{\cal B}}
\newcommand{\calc}{{\cal C}}
\newcommand{\cald}{{\cal D}}

\newcommand{\calr}{{\cal R}}

\newcommand{\calt}{{\cal T}}
\newcommand{\calk}{{\cal K}}
\newcommand{\calL}{{\cal L}}

\newcommand{\msten}{\mbox{\hspace{10mm}}}
\newcommand{\msfiv}{\mbox{\hspace{5mm}}}
\newcommand{\setm}{\setminus}
\newcommand{\incl}{\subseteq}

\title{Proof Theory for Theories of Ordinals III:$\Pi_{N}$-Reflection}
\author{Toshiyasu Arai
\\
Graduate School of Science,
Chiba University
\\
1-33, Yayoi-cho, Inage-ku,
Chiba, 263-8522, JAPAN
\\
tosarai@faculty.chiba-u.jp
}
\date{July 5, 2010}
\begin{document}
\maketitle
\begin{abstract}
This paper deals with a proof theory for a theory $\mbox{{\rm T}}_{N}$ of $\Pi_{N}$-reflecting ordinals using a system $Od(\Pi_{N})$ of ordinal diagrams. 
This is a sequel to the previous one \cite{ptpi3} in which a theory for $\Pi_{3}$-reflecting ordinals is analysed proof-theoretically.
\end{abstract}

\section{Prelude}\label{sec:prlpt5a}
This is a sequel to the previous ones \cite{ptMahlo} and \cite{ptpi3}. 
Namely our aim here is to give finitary analyses of finite proof figures in a theory for $\Pi_{N}$-reflecting ordinals, 
 \cite{Richter-Aczel74}
via cut-eliminations
as in Gentzen-Takeuti's consistency proofs, \cite{Gentzen38} and \cite{Takeuti87}.
Throughout this paper $N$ denotes a positive integer such that $N\geq 4$.

 



Let $T$ be a theory of ordinals. Let $\Omega$ denote the (individual constant corresponding to the) ordinal $\omega^{CK}_1$.
We say that $T$ is a $\Pi^{\Ome}_{2}${\it -sound theory} if 
\[
\forall\Pi_2\, A(T\vdash A^{\Ome}  \Rightarrow  A^{\Ome})
.\]
\begin{definition}{\rm (}$\Pi_2^\Ome$-ordinal of a theory{\rm )}  
{\rm Let} $T$ {\rm be a} $\Pi^{\Ome}_{2}${\rm -sound and recursive theory of ordinals.} {\rm For a sentence} $A$ {\rm let}  $A^{\alpha}$ {\rm denote the result of replacing unbounded quantifiers} $Qx \, (Q\in\{ \forall, \exists\})$ {\rm in} $A$ {\rm by} $Qx<\alpha$. {\rm Define the}
$\Pi_2^\Ome$-ordinal $\mid \! T \! \mid_{\Pi_{2}^{\Ome}}$ of  $T$ {\rm by}
\[\mid \! T \! \mid_{\Pi_{2}^{\Ome}}:=\inf\{\alpha\leq\omega^{CK}_1 :\forall\Pi_2\mbox{ {\rm sentence} } A(T\vdash A^{\Ome}  \Rightarrow  A^{\alp})\}<\omega^{CK}_1\]
\edf
Roughly speaking, the aim of proof theory for theories $T$ of ordinals is to describe the ordinal $\mid \! T \! \mid_{\Pi_{2}^{\Ome}}$. This gives $\Pi_2^\Ome$-ordinal of an equivalent theory of sets, cf. \cite{ptMahlo}.

Let $\mbox{{\rm KP}}\Pi_{N}$ denote the set theory for $\Pi_{N}$-reflecting universes. $\mbox{{\rm KP}}\Pi_{N}$ is obtained from the Kripke-Platek set theory with the Axiom of Infinity by adding the axiom: for any $\Pi_{N}$ formula $A(u)$
\[
A(u)\to \exi z(u\in z\spand A^{z}(u))
.\]

In \cite{WienpiN} we introduced a recursive notation system $Od(\Pi_{N})$ of ordinals,
which we studied first in \cite{hndWienpiN}. 
An element of the notation system is called an {\it ordinal diagram\/} (henceforth abbreviated by {\it o.d.\/}). 
The system is designed for proof theoretic study of theories of $\Pi_{N}$-reflection.  
We \cite{WienpiN} showed that for each $\alp<\Ome$ in $Od(\Pi_{N})$ 
$\mbox{{\rm KP}}\Pi_{N}$ proves that the initial segment of $Od(\Pi_{N})$ determined by $\alp$ is a well ordering.

Let $\mbox{{\rm T}}_{N}$ denote a theory of $\Pi_{N}$-reflecting ordinals. 
The aim of this paper is to show an upper bound theorem for the ordinal 
$\mid \! \mbox{{\rm T}}_{N} \! \mid_{\Pi_{2}^{\Ome}}$:

\bth\label{th:5}
$\fal\Pi_2 \, A(\mbox{{\rm T}}_{N}\vdash A^{\Ome} \, \Rarw \, \exi\alp\in Od(\Pi_{N})\mid \Ome \: A^\alp)$
\eth

Combining Theorem\ref{th:5} with the result in \cite{WienpiN} mentioned above yields the:

\bth\label{th:ptpiN}
$\mid \! \mbox{{\rm KP}}\Pi_{N} \! \mid^{\Ome}_{\Pi_{2}}=\mid \! \mbox{{\rm T}}_{N} \! \mid^{\Ome}_{\Pi_{2}}=\mbox{{\rm the order type of }} Od(\Pi_{N})\mid\Ome$
\eth

Proof theoretic study for $\Pi_{N}$-reflecting ordinals via ordinal diagrams were first obtained in a handwritten note \cite{hndptpiN}.

For an alternative approach to ordinal analyses of set theories, see
M. Rathjen's papers \cite{Rathjen94}, \cite{RathjenAFML1} and \cite{RathjenAFML2}.

Let us mention the contents of this paper.

In Section \ref{sec:5apreview} a preview of our proof-theoretic analysis for $\Pi_{N}$-reflection is given. 
As in \cite{ptpi3} inference rules $(c)^{\sig}_{\alp_{1}}$ are added to analyse an inference rule $(\Pi_{N}\mbox{{\rm -rfl}})$ saying the universe of the theory $\mbox{{\rm T}}_{N}$ is $\Pi_{N}$-reflecting. 
A {\it chain} is defined to be a consecutive sequence of rules $(c)$.

In Subsection \ref{subsec:5amerge} we expound that chains have to merge each other for a proof theoretic analysis of $\mbox{{\rm T}}_{N}$ for $N\geq 4$. 
An ordinal diagram in the system $Od(\Pi_{N})$ defined in \cite{WienpiN} may have its $Q$ part, 
which has to obey complicated requirements. 
In Subsection \ref{subsec:5aQpart} we explain what parts correspond to the $Q$ part in proof figures. 

In Section \ref{sec:5} the theory $\mbox{{\rm T}}_{N}$ for $\Pi_{N}$-reflecting ordinals is defined.
In Section \ref{sec:WienpiN} let us recall briefly the system $Od(\Pi_{N})$ of ordinal diagrams (abbreviated by o.d.'s) in \cite{WienpiN}. 

In Section \ref{sec:TNc} we extend $\mbox{{\rm T}}_{N}$ to a formal system $\mbox{{\rm T}}_{Nc}$. The language is expanded so that individual constants $c_{\alp}$ for o.d.'s $\alp\in Od(\Pi_{N})\mid\pi$ are included. Inference rules $(c)^{\sig}_{\alp_{1}}$ are added. {\it Proofs} in $\mbox{{\rm T}}_{Nc}$ defined in Definition \ref{df:5prf} are proof figures enjoying some provisos and obtained from given proofs in $\mbox{{\rm T}}_{N}$ by operating rewriting steps. Some lemmata for proofs are established. These are needed to verify that rewrited proof figures enjoy these provisos. To each proof $P$ in $\mbox{{\rm T}}_{Nc}$ an o.d. $o(P)\in Od(\Pi_{N})\mid\Ome$ is attached. Then the Main Lemma \ref{mlem:5} is stated as follows: If $P$ is a proof in $\mbox{{\rm T}}_{Nc}$, then the endsequent of $P$ is true.

In Section \ref{sec:ml5} the Main Lemma \ref{mlem:5} is shown by a transfinite induction on $o(P)\in Od(\Pi_{N})\mid\Ome$. 

This paper relies heavily on the previous ones \cite{ptMahlo} and \cite{ptpi3}.
\smlskp
{\bf General Coventions}. Let $(X,<)$ be a quasiordering. Let $F$ be a function $F:X\ni\alp\mapsto F(\alp)\incl X$. For subsets $Y,Z\subset X$ of $X$ and elements $\alp,\bet\in X$, put
\benu
\item $\alp\leq\bet\Lrarw \alp<\bet \mbox{ or } \alp=\bet$
\item $Y\mid \alp=\{\bet\in Y:\bet<\alp\}$
\item $Y<Z  :\Lrarw  \exi\bet \in Z\fal \alp \in Y(\alp<\bet)$
\item $Y<\bet  :\Lrarw  Y<\{\bet\} \Lrarw \fal \alp \in Y(\alp<\bet)$; $\alp< Z  :\Lrarw \{\alp\}< Z$
\item $Z\leq Y  :\Lrarw  \fal \bet \in Z\exi \alp \in Y(\bet\leq\alp)$
\item $\bet\leq Y  :\Lrarw  \{\bet\}\leq Y \Lrarw \exi\alp\in Y(\bet\leq\alp)$; $Z\leq \alp :\Lrarw Z\leq\{\alp\}$
\item $F(Y)=\bigcup\{F(\alp):\alp\in Y\}$
\eenu

\section{A preview of proof-theoretic analysis}\label{sec:5apreview}
In this section a preview of our proof-theoretic analysis for $\Pi_{N}$-reflection is given.

Let us recall briefly the system $Od(\Pi_{N})$ of o.d.'s in \cite{WienpiN}. 
%
 The main constructor in $Od(\Pi_{N})$ is to form an o.d. $d_{\sig}^{q}\alp$ 
 from a symbol $d$ and o.d.'s in $\{\sig,\alp\}\cup q$, 
 where $\sig$ denotes a recursively Mahlo ordinal and 
 $q=Q(d_{\sig}^{q}\alp)$ a finite sequence of quadruples of o.d.'s called 
 $Q$ {\it part\/} of $d_{\sig}^{q}\alp$. 
 By definition we set $d_{\sig}^{q}\alp<\sig$. 
 Let $\gam\prec_{2}\del$ denote the transitive closure of the relation 
 $\{(\gam,\del): \exi q,\alp(\gam=d_{\del}^{q}\alp)\}$, 
 and $\preceq_{2}$ its reflexive closure. 
 Then the set $\{\tau:\sig\prec_{2}\tau\}$ is finite and linearly ordered by $\prec_{2}$ 
 for each $\sig$.

An o.d. of the form $\rho=d_{\sig}^{q}\alp$ is introduced in proof figures 
only when an inference rule $(\Pi_{N}\mbox{{\rm -rfl}})$ for $\Pi_{N}$-reflection is 
resolved by using an inference rule $(c)_{\rho}$.


$q$ in $\rho=d_{\sig}^{q}\alp$ includes some data $st_{i}(\rho), rg_{i}(\rho)$ for $2\leq i<N$. $st_{N-1}(\rho)$ is an o.d. less than $\veps_{\pi+1}$ and $rg_{N-1}(\rho)=\pi$, 
while $st_{i}(\rho), rg_{i}(\rho)$ for $i<N-1$ may be undefined. 
If these are defined, then we write $rg_{i}(\rho)\darw$, etc. 
and $\kap=rg_{i}(\rho)$ is an o.d. such that $\rho\prec_{i}\kap$, 
where $\gam\prec_{i}\del$ is a transitive closure of the relation 
$pd_{i}(\gam)=\del$ on o.d.'s such that $\prec_{i+1}\incl\prec_{i}$ 
and $\prec_{2}$ is the same as one mentioned above. $q$ also includes data $pd_{i}(\rho)$. $st_{N-1}(\rho)$ is defined so that 
\beqn\label{eqn:1}
\gam\prec_{N-1}\rho\Rarw st_{N-1}(\gam)<st_{N-1}(\rho)
\eeqn

In Subsection \ref{subsec:5aQpart} we explain what parts correspond to the $Q$ part in proof figures. 
\smallskip

A theory $\mbox{{\rm T}}_{N}$ for $\Pi_{N}$-reflection is formulated in Tait's logic calculus, i.e., one-sided sequent calculus and $\Gamma,\Delta\ldots$ denote a {\it sequent\/}, i.e., 
a finite set of formulae. $\mbox{{\rm T}}_{N}$ has the inference rule $(\Pi_{N}\mbox{{\rm -rfl}})$:
 \[\infer[(\Pi_{N}\mbox{{\rm -rfl}})] {\Gamma}{\Gamma,A & \lnot\exi z A^{z},\Gamma}\]                  
 where $A\equiv\forall x_{N}\exists x_{N-1}\cdots Q x_{1}B$ with a bounded formula $B$.

So $(\Pi_{N}\mbox{{\rm -rfl}})$ says   $A\to \exists z A^{z}$.
\footnote{For simplicity we suppress the parameter. 
Correctly $\forall u(A(u) \, \to \, \exists z(u<z \, \& \,A^z(u)))$.}

To deal with the inference rule $(\Pi_{N}\mbox{{\rm -rfl}})$ we introduce 
new inference rules $(c)^{\sig}_{\rho}$ and $(\Sig_{i})^{\sig}\, (1\leq i\leq N)$ as in \cite{ptpi3}:
\[
\infer[(c)^{\sig}_{\rho}]{\Gam,\Lam^{\rho}}{\Gam,\Lam^{\sig}}
\]
where $\Lam$ is a set of $\Pi_{N}$-sentences as above, $\Lam^{\sig}=\{A^{\sig}:A\in\Lam\}$,
the side formulae $\Gam$ consists solely of $\Sig_{1}^{\sig}$-sentences 
and 
$\rho$ is of the form $d_{\sig}^{q}\alp$.

\[
\infer[(\Sig_{i})^{\sig}]{\Gam,\Lam}{\Gam,\lnot A^{\sig} & A^{\sig},\Lam}
\]
where $A$ is a $\Sig_{i}$ sentence. Although this rule $(\Sig_{i})^{\sig}$ is essentially 
a $(\mbox{{\rm cut}})$ inference, we need to distinguish between this and $(\mbox{{\rm cut}})$ to remember that a 
$(\Pi_{N}\mbox{{\rm -rfl}})$ was resolved.

When we apply the rule $(c)^{\sig}_{\rho}$ it must be the case:
\beqnarr
&& \mbox{ any instance term } \bet<\sig \mbox{ for the existential quantifiers } \exi x_{N-i}<\sig (i\mbox{:odd}) \nonumber \\
&& \mbox{ in } A^{\sig}\equiv\forall x_{N}<\sig\exists x_{N-1}<\sig\cdots Q x_{1}<\sig B \mbox{ is less than } \rho \label{eq:(c)}
\eeqnarr
As in \cite{ptpi3} an inference rule $(\Pi_{N}\mbox{{\rm -rfl}})$ is resolved 
by forming a succession of rules $(c)$'s, called a {\it chain\/}, 
which grows downwards in proof figures. 
We have to pinpoint, for each $(c)$, the unique chain, 
which describes how to introduce the $(c)$. 
To retain the uniqueness of the chain, i.e., not to branch or split a chain, 
we have to be careful in resolving rules with two uppersequents.  Our guiding principles are:
\bdes
\item[(ch1)] For any $\infer[(c)^\sig_\tau]{A^\tau}{A^\sig}$ with $\tau=d_{\sig}^{q}\alp$, if an o.d. $\bet$ is substituted for an existential quantifier $\exi y<\sig$ in $A^\sig$, i.e., $\bet$ is a realization for $\exi y<\sig$, then $\bet<\tau$, cf. (\ref{eq:(c)}),\\
 and
\item [(ch2)] Resolving rules having several upperseuquents must not branch a chain.
\edes

\subsection{Merging chains}\label{subsec:5amerge}
As contrasted with \cite{ptpi3} for $\Pi_{3}$-reflection we have to merge chains here. 
Let us explain this phenomenon.

We omit side formulae in this subsection.
\\
1)  First resolve a $(\Pi_{N}\mbox{{\rm -rfl}})$ in the left figure, and resolve the $(\Sig_{N})^{\sig} \, J_0$ to the right figure
with a $\Sig_{N-1} \: A_1$:
\[\infer[(\Sig_{N})^{\sig}\, J_0]{}
 {
  \infer[(c)^{\pi}_{\sig}]{A^\sig}{A} 
 &
 \lnot A^\sig
 }
\hskip1.0cm
\infer[J_1\, (\Sig_{N-1})^{\sig}]{}
{
 \infer{\lnot A_1^\sig}
 {
  \infer{A^\sig}{A}
 &
  \lnot A^\sig,\lnot A_1^\sig
 }
&
 \infer[(c)^{\pi}_{\sig}\, I_0]{A_1^\sig}{A_1}
}
\]
with $A\equiv \fal x_{N}\exi x_{N-1}\fal x_{N-2}A_{3}, \, \sig=d_{\pi}^{\alp}\alp$, where $A_{3}\equiv\exi x_{N-3}A_{4}$ is a $\Sig_{N-3}$-formula and $\alp$ denotes the o.d. attached to the uppersequent $A$ of $(c)^{\pi}_{\sig}$.
\\

\noindent
2)  Second resolve a $(\Pi_{N}\mbox{{\rm -rfl}})$ above the $(c)^{\pi}_{\sig}\,  I_0$ and 
a $(\Sig_{N})$ as in 1):
\[
\infer{}
 {
  \infer{\lnot B_1^\tau}
   {
     \infer{B^\tau}
      {
        \infer[(\Sig_{N-1})^{\sig}]{B^\sig}
         { 
           \lnot A_1^\sig        
          &
            \infer[(c)^\pi_\sig\, \tilde{I}_0]{A_1^\sig,B^\sig}
                                   {\deduce[P_{1}]{A_1,B}{}}
          }
       }
    &
      \infer{\lnot B^\tau,\lnot B_1^\tau}
        {
           \lnot A_1^\sig
         &
           \infer{A_1^\sig,\lnot B^\tau,\lnot B_1^\tau}
                 {A_1,\lnot B^\tau,\lnot B_1^\tau}
         }
    }
 &
  \infer[(c)^\sig_\tau\, I_{1}]{B_1^\tau}
    {
      \infer[J_1]{B_1^\sig}
        {
         \lnot A_1^\sig
         &
          \infer[(c)^\pi_\sig]{A_1^\sig,B_1^\sig}{A_1,B_1}
         }
     }
 }
Fig.1
\]
\noindent
with a $\tau=d_{\sig}^{\nu}\bet$ and a $\Sig_{N-1} \: B_{1}\equiv\exi y_{N-1}\fal y_{N-2}B_3$, where $\nu$ denotes the o.d. attached to the subproof $P_{1}$ ending with the uppersequent $A_{1},B$ of $(c)^\pi_\sig\, \tilde{I}_0$.
 \\          
After that resolve the $(\Sig_{N-1})^{\sig}\, J_1$:
\[
\infer{}
 {
  \infer{\lnot B_1^\tau}
   {
     \infer{B^\tau}
      {
        \infer{B^\sig}
         { 
           \lnot A_1^\sig        
          &
            \infer[\tilde{I}_0]{A_1^\sig,B^\sig}
                                   {\deduce[P_{2}]{A_1,B}{}}
          }
        }
    &
      \lnot B^\tau,\lnot B_1^\tau
    }
 &
  \infer{B_1^\tau}
    {
      \infer{B_1^\sig}
        {
           \infer{B_1^\sig,A_2^\sig}
            {
               \lnot A_1^\sig
             &
               \infer{A_1^\sig,B_1^\sig,A_2^\sig}{A_1,B_1,A_2}
             }
         &
            \infer[J_0']{\lnot A_2^\sig}
             {
                \infer{A^\sig}{A}
              &
                \lnot A^\sig,\lnot A_2^\sig
              }
         }
     }
 }
\]
Then resolve the $(\Sig_{N})^{\sig} \, J_0'$ :
\[
\infer[K]{}
 {
  \infer{\lnot B_1^\tau}
   {
     \infer{B^\tau}
      {
        \infer{B^\sig}
         { 
           \lnot A_1^\sig        
          &
            \infer[\tilde{I}_0]{A_1^\sig,B^\sig}
                                   {\deduce[P_{3}]{A_1,B}{}}
          }
       }
    &
      \lnot B^\tau,\lnot B_1^\tau
    }
 &
  \infer[(c)^\sig_\tau\,I_1]{B_1^\tau}
    {
      \infer[(\Sig_{N-2})^{\sig}\, J_2]{B_1^\sig}
        {
           B_1^\sig,A_2^\sig
         &
            \infer{\lnot A_2^\sig}
             {
                \lnot A_2^\sig,\lnot \tilde{A}_1^\sig
              &
                \infer[(c)^\pi_\sig\, I_0']{\tilde{A}_1^\sig}
                  {
                    \infer*{\tilde{A}_1}{(\Pi_{N}\mbox{{\rm -rfl}}) \, H}
                  }
              }
         }
     }
 }
\]
3) Thirdly resolve a $(\Pi_{N}\mbox{{\rm -rfl}}) \, H$ above the $(c)^\pi_\sig \, I_0'$. 
One cannot resolve the $(\Pi_{N}\mbox{{\rm -rfl}}) \, H$ by introducing a 
$(c)^{\sig}_{\rho}$ with $\rho<\tau$. 
Let me explain the reason.

Suppose that we introduce a new $(c)^\sig_\rho \, I_1'$ with $\rho=d_{\sig}^{\eta}\gam$ immediately above the $(\Sig_{N-2})^{\sig}\,J_2$ as in \cite{ptpi3}. 
Then the new $(c)^\sig_\rho \, I_1'$ is introduced after the $(c)^\sig_\tau\,  I_1$  
and so $\rho=d_{\sig}^{\eta}\gam<\tau$. 
Hence a new $(\Sig_{N})^{\rho} \, K'$ is introduced below the $(\Sig_{N-1})^{\tau} \, K$:
\[
\infer[K']{}
{
  \infer[K]{D^\rho}
   {
     \lnot B_1^\tau
   &
     \infer[(c)^\sig_\tau\, I_{1}]{B_1^\tau,D^\rho}
    {
      \infer[J_2]{B_1^\sig,D^\rho}
        {
           B_1^\sig,A_2^\sig
         &
            \infer[(c)^\sig_\rho\, I_{1}']{\lnot A_2^\sig,D^\rho}
             {
               \infer{\lnot A_2^\sig,D^\sig}
                 {
                   \lnot A_2^\sig,\lnot \tilde{A}_1^\sig
                 &
                   \infer[(c)^\pi_\sig]{\tilde{A}_1^\sig,D^\sig}
                         {\tilde{A}_1,D}
                  }
               }
           }
       }
     }
 &
 \lnot D^\rho
 }
\msfiv Fig. 2\]
\\
with $D\equiv\fal z_{N}\exi z_{N-1}\fal z_{N-2}D_3$. 
Nevertheless this does not work, 
because $\lnot A_{2}\equiv\exi x_{N-3}\fal x_{N-4}\lnot A_{4}$ is 
a $\Sig_{N-2}$ sentence with $N-2\geq 2$. 
Namely the principle {\bf (ch1)} may break down for the $(c)^\sig_\rho\, I_1'$ 
since any o.d. $\del<\sig$, i.e., possibly $\del\geq\rho$ may be 
an instance term for the existential quantifier $\exi x_{N-3}$ in 
$A_{2}\equiv\fal x_{N-2}\exi x_{N-3} A_{4}$ and may be substituted for the variable 
$x_{N-3}$ in $\lnot A_{2}^{\sig}$. 
Only we knows that such a $\del$ is less than $\sig$ 
and comes from the left upper part of $J_{2}$.
\\

\noindent
4)  Therefore the chain for $H$ has to connect or merge with the chain $I_0-I_1$ for $B$: 
 \[\infer[(c)^\tau_\rho\, I_2]{D^\rho}
   {
    \infer{D^\tau}
     {
      \infer{\lnot B_1^\tau}
      {
        \infer[I_1']{B^\tau}
         {
          \infer{B^\sig}
           { 
            \lnot A_1^\sig        
           &
            \infer[\tilde{I}_0]{A_1^\sig,B^\sig}
                                   {\deduce[P_{4}]{A_1,B}{}}
            }
          }
       &
         \lnot B^\tau,\lnot B_1^\tau
       }
      &
       \infer[I_1]{B_1^\tau,D^\tau}
       {
        \infer[J_{2}]{B_1^\sig,D^\sig}
        {
           B_1^\sig,A_2^\sig
         &
           \hskip-1.0cm
           \infer{\lnot A_2^\sig,D^\sig}
             {
                \infer[J_{0}']{\lnot A_2^\sig,\lnot \tilde{A}_1^\sig}
                 {
                   \infer[I_{0}'']{A^\sig}{A}
                  &
                   \lnot A^\sig, \lnot A_2^\sig,\lnot \tilde{A}_1^\sig
                  }
              &
                \infer[I_0']{\tilde{A}_1^\sig,D^\sig}
                                         {\tilde{A}_1,D}
              }
          }
        }
      }
    }
\]
\hspace*{\fill} $Fig. 3$
\\
with $\rho=d_{\tau}^{\eta}\gam$ and a $(\Sig_{N})^{\rho}$ 
with the cut formula $D^\rho$ follows this figure as in $Fig.2$, 
where $\eta$ denotes the o.d. attached to the uppersequent 
$\tilde{A}_1,D$ of $(c)^{\pi}_{\sig}\, I_{0}'$. $(\Sig_{N-2})^{\sig}\, J_{2}$ 
is a merging point for chains $I_{0}-I_{1}$ and $I_{0}'-I_{1}-I_{2}$.
\smallskip

The principle {\bf (ch1)} for the new $(c)^\tau_\rho \, I_2$ will be retained 
for the simplest case $N=4$ as in \cite{ptpi3}. 
The problem is that the proviso (\ref{eqn:1}) may break down:
it may be the case $\nu=st_{N-1}(\tau)\leq st_{N-1}(\rho)=\eta$ 
since we cannot expect the upper part of $(c)^\pi_\sig\, I_0'$ is simpler than the one of $(c)^\pi_\sig\, \tilde{I}_0$. 

In other words a new succession $I_0'-I_1-I_2$ of collapsings starts. 
This is required to resolve $\Sig_{N-2}^{\sig}$ sentence $\lnot A_{2}^{\sig}\, (N-2\geq 2)$ 
and hence $\sig$ has to be $\Pi_{N-1}$-reflecting.

If this chain $I_0'-I_1-I_2$ would grow downwards as in $\Pi_{3}$-reflection, i.e., 
in a chain 
$I_0'-I_1-I_2-\cdots-I_n$, $I_n$ would come only from the upper part of $I_0'$, 
then the proviso (\ref{eqn:1}) would suffice to kill this process. 
But the whole process may be iterated : in $Fig.3$ another succession 
$I_{0}''-I_1-I_2-I_3$ may arise by resolving the $(\Sig_{N})^{\sig} \, J_0'$.

    Nevertheless still we can find a reducing part, that is, 
the upper part of the $(c)^{\tau}_{\rho} \, I_{2}$: 
the upper part of the $(c)^{\tau}_{\rho} \, I_{2}$ becomes simpler in 
the step $I_2-I_3$. 
Furthermore in the general case $N>4$ merging processes could be iterated, vz.
the merging point $(\Sig_{N-2})^{\sig}\, J_{2}$ may be resolved into a 
$(\Sig_{N-3})^{\rho_{1}}$, 
which becomes a new merging point to analyse a $\Sig_{N-3}$ sentence $A_{3}^{\rho_{1}}$ where $\rho_{1}\preceq\rho$ is a $\Pi_{N-2}$-reflecting and so on. 
Therefore in $Od(\Pi_{N})$ the $Q$ part of an o.d. may consist of several factors: 
\[
(\tau,\alp,q=\{\nu_{i},\kap_{i},\tau_{i}:i\in In(\rho)\}) \mapsto d_{\tau}^{q}\alp=\rho
\]
with $\kap_{N-1}=rg_{N-1}(\rho)=\pi$. 
$In(\rho)$ denotes a set such that 
\[
N-1\in In(\rho)\incl\{i: 2\leq i\leq N-1\}.
\]
We set for $i\in In(\rho)$:
\[
st_{i}(\rho)=\nu_{i}, rg_{i}(\rho)=\kap_{i}, pd_{i}(\rho)=\tau_{i}.
\]
If $i\not\in In(\rho)$, 
set
\[
pd_{i}(\rho)=pd_{i+1}(\rho), st_{i}(\rho)\simeq st_{i}(pd_{i}(\rho)), rg_{i}(\rho)\simeq rg_{i}(pd_{i}(\rho)).
\]
Also these are defined so that 
$pd_{2}(\rho)=\tau$ for $\rho=d_{\tau}^{q}\alp$.

For the o.d. $\rho=d_{\tau}^{q}\gam$ in the $Fig. 3$, $In(\rho)=\{N-2,N-1\}, st_{N-1}(\rho)=\eta, pd_{N-1}(\rho)=\sig, rg_{N-2}(\rho)=\tau=pd_{N-2}(\rho), st_{N-2}(\rho)=\gam=st_{2}(\gam)$.

Thus $\nu_{i}=st_{i}(\rho)$ corresponds to the upper part of a $(c)^{rg_{i}(\rho)}$ 
while $\tau_{N-1}=pd_{N-1}(\rho)$ indicates that the first, i.e., 
uppermost merging point for a chain ending with a $(c)_{\rho}$ is a rule 
$(\Sig_{N-2})^{\tau_{N-1}}$, e.g., the rule $J_{2}$ in $Fig. 3$. 
Note that $st_{N-1}(\rho)=\eta<st_{N-1}(pd_{N-1}(\rho))$, cf. (\ref{eqn:1}). 
$\kap_{i}=rg_{i}(\rho)$ is an o.d. such that there exists a $(c)^{\kap_{i}}$ 
in the chain for $(c)_{\rho}$. 
We will explain how to determine the rule $(c)^{rg_{i}(\rho)}$, i.e., 
the point to which we direct our attention in Subsection \ref{subsec:5aQpart}.

The case $In(\rho)=\{N-1\}$ corresponds to the case when a $(c)^{pd_{N-1}(\rho)}_{\rho}$ 
is introduced without merging points, i.e., as a resolvent of a $(\Pi_{N}\mbox{{\rm -rfl}})$ 
above the top of the chain whose bottom is a $(c)_{pd_{N-1}(\rho)}$.
The case $In(\rho)=\{N-2,N-1\}$ corresponds to the case when a 
$(c)^{pd_{2}(\rho)}_{\rho}\, (pd_{2}(\rho)=pd_{N-2}(\rho))$ is introduced 
with a merging point $(c)^{pd_{N-1}(\rho)}$.

In $Fig.3$ a new succession with a merging point $(c)^\tau_\rho\, I_2$ 
arises by resolving a $(\Sig_{N})^{\tau}$ below the $(c)^\sig_\tau \, I_1'$, i.e., 
$\tilde{I}_0-I_1'-I_2-I_3 \, (c)^\rho_\kap$ for a $\kap$ with a $\lam=st_{N-1}(\kap)$. 
But in this case we have
  \[
  \lam=st_{N-1}(\kap)<st_{N-1}(\tau)=\nu
  .\]
$st_{N-1}(\kap)$ corresponds to the upper part $P_{1}$ of a 
$(c)^\pi_\sig \, \tilde{I}_0$ in $Fig.1$, 
when the $(c)^\sig_\tau$ was originally introduced. 
This part $P_{1}$ is unchanged up to $Fig.3$:\\
$P_{1}=P_{2}=P_{3}=P_{4}$. 
Roughly speaking, $\tilde{I}_0-I_1'-I_3$ can be regarded as a $\Pi_{N-1}$ resolving series 
$I_0-I_1-I_3$. 
This prevents the new merging points from going downwards unlimitedly.

\subsection{The $Q$ part of an ordinal diagram}\label{subsec:5aQpart}
In this subsetion we explain how to determine the $Q$ part $q$ of 
$\rho=d_{\sig}^{q}\alp$ from a proof figure when an inference rule $(c)^{\sig}_{\rho}$ 
is introduced.

In general such a $(c)^{\sig}_{\rho}$ is formed when we resolve an inference rule 
$(\Pi_{N}\mbox{{\rm -rfl}})\, H$:

\[
\deduce
{\hskip2.6cm
 \infer[(c)^{\sig_{n-1}}_{\sig}\, J_{n-1}]{\Gam_{n-1}'}
 {
  \infer*{\Gam_{n-1}}
  {
   \infer[(c)^{\sig_{n_{m+1}}}\, J_{n_{m+1}}]{\Gam_{n_{m+1}}'}
   {
    \infer*[\calc_{n_{m+1}}]{\Gam_{n_{m+1}}}
    {}
   }
  }
 }
}
{\hskip0.0cm
  \infer[(\Sig_{i_{m}})^{\sig_{n_{m}+1}}\, K_{m}]{\Phi_{m},\Psi_{m}}
  {
   \deduce
   {\hskip0.5cm
    \infer*[\calc_{n_{m+1}}]{\Phi_{m},\lnot A_{m}}
    {
     \infer[(c)^{\sig_{n_{m}}}_{\sig_{n_{m}+1}}\, J^{m+1}_{n_{m}}]{(\Gam_{n_{m}}')^{m+1}}
     {}
     }
    }
    {\hskip-2.1cm
     \infer*{\Gam^{m+1}_{n_{m}}}
     {
      \infer[(c)^{\sig_{p}}_{\sig_{p+1}}\, J^{m+1}_{p}]{(\Gam_{p}')^{m+1}}
      {
       \infer*[\calc_{n_{m+1}}]{\Gam^{m+1}_{p}}{}
      }
     }
    }
  &
   \deduce
   {\hskip0.0cm
    \infer*{A_{m},\Psi_{m}}
    {
     \infer[(c)^{\sig_{n_{m}}}_{\sig_{n_{m}+1}}\, J_{n_{m}}]{\Gam_{n_{m}}'}
     {
      \infer*{\Gam_{n_{m}}}
      {}
      }
     }
    }
    {\deduce
     {\hskip-0.3cm
      \infer[(c)^{\sig_{p}}_{\sig_{p+1}}\, J_{p}]{\Gam_{p}'}
      {
       \infer*{\Gam_{p}}
       {}
       }
      }
      {\hskip-2.0cm
       \infer[(c)^{\pi}_{\sig_{1}}\, J_{0}]{\Gam_{0}'}
       {
        \infer*{\Gam_{0}}{(\Pi_{N}\mbox{{\rm -rfl}})\, H}
        }
       }
     }
   }
}
Fig.4\]
where $\calr=J_0,\ldots,J_{n-1}$ denotes a series of rules $(c)^{\sig_p}_{\sig_{p+1}}\, J_{p}$ with $\pi=\sig_{0}, \sig=\sig_{n}$. $(\Pi_{N}\mbox{{\rm -rfl}})\, H$ is resolved into a $(c)^{\sig}_{\rho}\, J_{n}$ and a $(\Sig_{N})^{\rho}$ below $J_{n-1}$.

This series $\calr$ is devided into intervals $\{\calr_{m}=J_{n_{m-1}+1},\ldots,J_{n_{m}}:m\leq l\}$ with an increasing sequence $n_{-1}+1=0\leq n_0<n_1<\cdots<n_l=n-1\, (l\geq 0)$ of numbers so that 
\benu
\item $\calr_{0}=J_{0},\ldots,J_{n_{0}}$ is a chain $\calc_{n_{0}}$ leading to $J_{n_{0}}$.
\item For $m<l$ $\calr_{m+1}=J_{n_{m}+1},\ldots,J_{n_{m+1}}$ is a tail of a chain $\calc_{n_{m+1}}=J^{m+1}_{0},\ldots,J^{m+1}_{n_{m}},J_{n_{m}+1},\ldots,J_{n_{m+1}}$ leading to $J_{n_{m+1}}$ such that the chain $\calc_{n_{m+1}}$ passes through the left side of an inference rule $(\Sig_{i_m})^{\sig_{n_{m}+1}}\, K_{m}$ with $2\leq i_{m}<N-1$, $J_{n_m}$ is above the right uppersequent $A_{m},\Psi_{m}$ and $J^{m+1}_{n_{m}}$ is above the left uppersequent $\Phi_{m},\lnot A_{m}$ of $K_{m}$, resp. $A_{m}$ is a $\Sig_{i_{m}}$ sentence. Each rule $J^{m+1}_{p}$ for $p\leq n_{m}$ is again an inference rule $(c)^{\sig_{p}}_{\sig_{p+1}}$. $K_{m}$ will be a merging point of chains $\calc_{n_{m+1}}$ and a new chain $\calc_{\rho}=J_{0},\ldots,J_{n-1},J_{n}$ leading to $(c)_{\rho}\, J_{n}$.
\item There is no such a merging point below $J_{n-1}$, vz. there is no $(\Sig_{k})^{\sig}$ with $1<k<N-1$ such that $J_{n-1}$ is in the right upper part of the inference rule and there exists a chain passing through its left side.
\eenu

Set $N-1\in In(\rho), rg_{N-1}(\rho)=\pi$ and $st_{N-1}(\rho)$ is the o.d. attached to the upper part of $(c)^{\pi}\, J_{0}$, where by the upper part we mean the part after resolving $(\Pi_{N}\mbox{{\rm -rfl}})\, H$. 

First consider the case $l=0$, i.e., there is no merging point for the new chain $\calc_{\rho}$ leading to the new $J_{n}$. Then set $In(\rho)=\{N-1\}$ and $pd_{N-1}(\rho)=\sig$.

Suppose $l>0$ in what follows.
Then set $pd_{N-1}(\rho)=\sig_{n_{0}+1}$, i.e., $pd_{N-1}(\rho)$ is the superscript of the first uppermost merging point $(\Sig_{i_{0}})^{\sig_{n_{0}+1}}\, K_{0}$.

In any cases we have $st_{N-1}(\rho)<st_{N-1}(pd_{N-1}(\rho))$, cf. (\ref{eqn:1}).
$st_{i}(\rho)$ always corresponds to the upper part of a $(c)^{rg_{i}(\rho)}$ in the chain $\calc_{\rho}$ for $i\in In(\rho)$

\subsubsection{The simplest case $N=4$}\label{subsubsec:N4}
Here suppose $N=4$ and we determine the $Q$ part of $\rho$. First set $2\in In(\rho)$, vz. $In(\rho)=\{2,3\}$ and $pd_{2}(\rho)=\sig$. It remains to determine the o.d. $rg_{2}(\rho)$. In other words to specify a rule $(c)^{\sig_{q}}\, J_{q}$ with $rg_{i}(\rho)=\sig_{q}$.

Note that $i_{m}=2$ for any $m$ with $0<m\leq l$ since $2\leq i_{m}<N-1=3$ in this case. There are two cases to consider. First suppose there is a $p<n$ such that
\benu
\item $p>n_{0}$, i.e., $\sig_{p+1}\prec_{2}\sig_{n_{0}+1}=pd_{3}(\rho)$ and
\item $2\in In(\sig_{p+1})$, i.e., there was a merging point of the chain leading to $(c)_{\sig_{p+1}}\, J_{p}$.
\eenu
Then pick the minimal $q$ satisfying these two conditions, vz. the uppermost rule  $(c)^{\sig_{q}}_{\sig_{q+1}}\, J_{q}$ below the first uppermost merging point $(\Sig_{i_{0}})^{\sig_{n_{0}+1}}\, K_{0}$ with $2\in In(\sig_{q+1})$.
Then set 
\bdes
\item[Case 1] $rg_{2}(\rho)=rg_{2}(\sig_{q+1})$.
\edes
Otherwise set
\bdes
\item[Case 2] $rg_{2}(\rho)=\sig=pd_{2}(\rho)$.
\edes

Consider the first case {\bf Case 1} $rg_{2}(\rho)=rg_{2}(\sig_{q+1})\neq pd_{2}(\rho)$. From the definition we see $rg_{2}(\rho)=rg_{2}(\sig_{q+1})=pd_{2}(\sig_{q+1})=\sig_{q}$.
We have $\sig_{q}=rg_{2}(\rho)\preceq_{3}pd_{3}(\rho)=\sig_{n_{0}+1}$. This follows from the minimality of $q$, i.e., $\fal t[n_{0}<t<q \to 2\not\in In(\sig_{t+1})]$ and hence $\fal t[n_{0}<t<q \to \sig_{t}=pd_{2}(\sig_{t+1})=pd_{3}(\sig_{t+1})]$.

Furthermore $q$ is minimal, i.e, $\sig_{q}$ is maximal in the following sense:
\beqnarr
&& \fal t[n_{0}<t<n(\lrarw pd_{2}(\rho)=\sig\preceq_{2}\sig_{t+1}\prec_{2}pd_{3}(\rho))\spand rg_{2}(\sig_{t+1})\darw 
\nonumber \\
&& \to rg_{2}(\sig_{t+1})\preceq_{2}\sig_{q}] \label{eqn:N4.1}
\eeqnarr
In general we have the following fact. 
\bprp\label{prp:N4.1}
Let $\calc=J_{0},\ldots,J_{n-1}$ be a chain leading to a $(c)^{\sig_{n-1}}_{\sig_{n}}\, J_{n-1}$. Each $J_{p}$ is a rule $(c)^{\sig_{p}}_{\sig_{p+1}}$ with $\sig_{0}=\pi$. Suppose that $2\in In(\sig_{n})$ and the chain passes through the left side of a $(\Sig_{2})^{\sig_{p}}\, K$ for a $p$ with $0<p<n$ so that $J_{p-1}$ is in the left upper part of $K$ and $J_{p}$ is below $K$. Then $\sig_{q}=rg_{2}(\sig_{n})\preceq_{2}\sig_{p}$, i.e., $q\geq p$.
\eprp
\[
\deduce
{\hskip0.0cm
 \infer[(c)^{\sig_{n-1}}_{\sig}\, J_{n-1}]{\Gam_{n-1}'}
 {
  \infer*{\Gam_{n-1}}
  {
   \infer[(c)^{\sig_{p}}\, J_{p}]{\Gam_{p}'}
   {
    \infer*[\calc]{\Gam_{p}}
    {}
   }
  }
 }
}
{\hskip-2.1cm
  \infer[(\Sig_{2})^{\sig_{p}}\, K]{\Phi,\Psi}
  {
   \infer*{\Phi,\lnot A^{\sig_{p}}}
   {
    \infer[(c)^{\sig_{p-1}}_{\sig_{p}}\, J_{p-1}]{\Gam_{p-1}'}
    {
     \infer*[\calc]{\Gam_{p-1}}{}
     }
    }
  &
   \infer*{A^{\sig_{p}},\Psi}{}
 }
}
\]

This means that when, in $Fig.4$ a $(\Sig_{3})^{\sig_{t}}\, K^{3}\, (0<t\leq n)$ 
in the new chain $\calc_{\rho}=J_{0},\ldots,J_{n-1},J_{n}$ leading to 
$(c)^{\sig}_{\rho}\, J_{n}$ is to be resolved into a $(\Sig_{2})^{\sig_{t}}\, K^{2}$, 
then $t\leq q$, i.e., $rg_{2}(\rho)=\sig_{q}\preceq_{2}\sig_{t}$. 
In other words any $(\Sig_{3})^{\sig_{t}}$ with $q<t\leq n$, 
equivalently $(\Sig_{3})^{\sig_{t}}$ which is below $(c)^{\sig_{q}}\, J_{q}$ has to wait 
to be resolved, until the chain $\calc_{\rho}$ will disapper by inversion.

For example consider, in $Fig.4$, an inference rule $(\Sig_{3})^{\sig_{t}}\, K^{3}$ for 
$t=n_{m+1}+1$. 
Its right cut formula is a $\Sig_{3}^{\sig_{t}}$ sentence $C^{\sig_{t}}$ and 
a descendent of a $\Sig_{3}$ sentence $C$: 
a series of sentences from $C$ to $C^{\sig_{t}}$ are in the chain 
$\calc_{n_{m+1}}=J^{m+1}_{0},\ldots,J^{m+1}_{n_{m}},J_{n_{m}+1},\ldots,J_{n_{m+1}}$ leading to $J_{n_{m+1}}$. 
Then the chain $\calc_{n_{m+1}}$ passes through the left side of the inference rule 
$(\Sig_{i_m})^{\sig_{n_{m}+1}}\, K_{m}$ and hence $K^{3}$ will not be resolved until 
$K_{m}$ will be resolved and its right upper part will disapper 
since we always perform rewritings of proof figure on the rightmosot branch. 
But then the chain $\calc_{\rho}$ will disapper by inversion 
since it passes through the right side of $K_{m}$.
In this way we see Proposition \ref{prp:N4.1}, cf. Lemma \ref{lem:5.3.22.A} 
in Section \ref{sec:TNc} for a full statemnt and a detailed proof.

(\ref{eqn:N4.1}) is seen from Proposition \ref{prp:N4.1} and the minimality of $q$.
Thus we have shown, cf. the conditions $(\cald^{Q}.1)$ for $Od(\Pi_{4})$ in \cite{WienpiN} or Section \ref{sec:WienpiN},
\[
rg_{2}(\sig)=rg_{2}(pd_{2}(\rho))\preceq_{2}rg_{2}(\rho)\preceq_{3}pd_{3}(\rho)
\]
and
\[
\fal t[rg_{2}(pd_{2}(\rho))\preceq_{2}\sig_{t}\prec_{2}\sig_{q} \Rarw rg_{2}(\sig_{t})\preceq_{2}\sig_{q}]
.\]
Furthermore we have
\beqn\label{eqn:N4.4}
st_{2}(\rho)<st_{2}(\sig_{p+1})<\sig_{q}^{+}
\eeqn
for the maximal $p$, vz. for the latest $(c)_{\sig_{p+1}}\, J_{p}$ with 
$rg_{2}(\sig_{p+1})=\sig_{q}\spand 2\in In(\sig_{p+1})$. 

Let $m<l$ denote the number such that $n_{m}<q\leq n_{m+1}$, i.e., 
$J_{q}$ is a member of the tail $\calr_{m+1}=J_{n_{m}+1},\ldots,J_{n_{m+1}}$ of the chain 
$\calc_{n_{m+1}}$. 
Then from Proposition \ref{prp:N4.1} we see that $J_{p}$ is also a member of 
$\calr_{m+1}$ and further that $J_{q}$ is a member of a chain 
$\calc_{p}$ leading to $J_{p}$. Thus the upper part of $(c)^{\sig_{q}}\, J_{q}$ 
corresponding to $st_{2}(\rho)$ is a result of perfoming 
several non-void rewritings to the upper part of a $(c)^{\sig_{q}}$ 
which determined $st_{2}(\sig_{p+1})$ when $(c)_{\sig_{p+1}}\, J_{p}$ was introduced originally. 
This yields (\ref{eqn:N4.4}).

Thus we have established the conditions $(\cald^{Q}.1)$ in \cite{WienpiN} or Section \ref{sec:WienpiN} for the newly introduced $\rho$.

Why we choose such a $\sig_{q}$ as $rg_{2}(\rho)$? 
First introducing $\sig_{q}=rg_{2}(\rho)$ is meant to express the fact that $\sig_{q}$ is 
(iterated) $\Pi_{3}$-reflecting and it is responsible to $\Sig_{2}^{\sig_{q}}$ sentences 
occurring above a $(c)^{\sig_{q}}$. 
Therefore even if there exists a $\sig_{p+1}$ above $pd_{3}(\rho)$, i.e., 
$p\leq n_{0}$ such that $2\in In(\sig_{p+1})$, 
we ignore these in determining $rg_{2}(\rho)$. 
Second in the {\bf Case 1} the reason why we chose $\sig_{q}$ as the uppermost one is 
explained by Proposition \ref{prp:N4.1}: 
any $(\Sig_{3})^{\sig_{t}}$ in the new chain $\calc_{\rho}$ will not be resolved for 
$q<t\leq n$ until the chain $\calc_{\rho}$ will disapper by inversion. 
Hence any $\sig_{q_{1}}$ with $rg_{2}(\sig_{p_{1}+1})=\sig_{q_{1}}\prec_{2}\sig_{q}$ 
for some $p_{1}\leq n$ will not be $rg_{2}(\kap)$ for $\kap\prec_{2}\rho$ in the future. 
This means that a collapsing series $\{(c)_{\kap}: rg_{2}(\kap)=\sig_{q_{1}}\}$ 
expressing the fact that $\sig_{q_{1}}$ is $\Pi_{3}$-reflecting is killed by introducing 
$\rho$ such that $\rho\prec_{2}\sig_{q_{1}}\prec_{2}\sig_{q}=rg_{2}(\rho)$. 
Therefore once we introduce such a $\rho$, then we can ignore 
$rg_{2}(\sig_{p_{1}+1})=\sig_{q_{1}}$ between $rg_{2}(\rho)$ and $\rho$.

\subsubsection{The general case $N>4$}\label{subsubsec:N5}
Here suppose $N>4$ and we determine the $Q$ part of $\rho$, i.e., 
determine the set $In(\rho)$ and o.d.'s $pd_{i}(\rho), rg_{i}(\rho)$ for 
$i\in In(\rho)$ by referring $Fig.4$.

First set $i_{0}\in In(\rho)$ where $i_{0}$ denotes the number such that the 
first merging point is a $(\Sig_{i_{0}})^{\sig_{n_{0}+1}}\, K_{0}$. 
Now let us assume inductively that for $k_{0}\geq 0$ we have specified merging points 
$\{K_{m_{k}}:k\leq k_{0}\}$ so that $0=m_{0}<\cdots<m_{k_{0}}$, 
$N-1>i_{m_{0}}>\cdots>i_{m_{k_{0}}}\geq 2$ and 
$\fal m\fal k<k_{0}[m_{k}<m<m_{k+1} \to i_{m}\geq i_{m_{k}}]$, and have setted 
$\{i_{m_{k}}:k\leq k_{0}\}\incl In(\rho)$. 
Namely $K_{m_{0}},\ldots, K_{m_{k_{0}}}$ is a series of merging points going downwards 
with decreasing indices $i_{m_{k}}$ and $K_{m_{k}}$ is the uppermost merging point 
with $i_{m_{k}}<i_{m_{k-1}}\, (i_{m_{-1}}:=N-1)$. 

If there exists an $m<l$ such that $m_{k_{0}}<m\spand i_{m_{k_{0}}}>i_{m}\geq 2$, 
then let $m$ denote the minimal one, vz. the uppermost merging point $K_{m}$ 
below the latest one $K_{m_{k_{0}}}$ with $i_{m_{k_{0}}}>i_{m}$, 
and set $i_{m}\in In(\rho)$. Otherwise set 
\[
In(\rho)=\{i_{m_{k}}:k\leq k_{0}\}\cup\{N-1\}
.\]
 This completes a description of the set 
\beqnarrs
In(\rho) & = & \{N-1=i_{m_{-1}}\}\cup\{i_{m_{k}}: 0\leq k\leq k_{1}\} \\
& = & \{N-1=i_{m_{-1}}>i_{m_{0}}>\cdots>i_{m_{k_{1}}}\}.
\eeqnarrs
Observe that for $i<N-1$
\[
i\in In(\rho) \Lrarw \exi m<l[i_{m}=i\spand \fal p<m(i_{p}\geq i)]
.\]

Now set $pd_{i_{m_{k}}}(\rho)=\sig_{n_{m_{k+1}}+1}$ for $-1\leq k\leq k_{1}$ with
 $m_{k_{1}+1}:=l$, vz. the merging point $K_{m_{k}}$ chosen for $i_{m_{k}}\in In(\rho)$ is a 
 $(\Sig_{i_{m_{k}}})^{pd_{i_{m_{k-1}}}(\rho)}$ for $0\leq k\leq k_{1}$ and 
 $pd_{2}(\rho)=pd_{i_{m_{k_{1}}}}(\rho)=\sig_{n_{l}+1}=\sig_{n}=\sig$. 
 Observe that for any $i$ with $2\leq i\leq N-1$ there exists an 
 $m(i)\leq l$ such that $pd_{i}(\rho)=\sig_{n_{m(i)}+1}$ and this $m(i)$ is the minimal 
 $m$ for which $i_{m}<i$.

It remains to determine the o.d.'s $rg_{i}(\rho)$ for $N-1\neq i=i_{m_{k}}\in In(\rho)$. 
As in the case $N=4$ there are two cases to consider. 
First suppose there is a $p<n$ such that
\benu
\item $\rho\prec_{i}\sig_{p+1}\prec_{i}\sig_{n_{m_{k}}+1}=pd_{i_{m_{k-1}}}(\rho)=pd_{i+1}(\rho)$ and
\item $i\in In(\sig_{p+1})$.
\eenu
Then pick the minimal $p$ satisfying these two conditions, vz. the uppermost rule  
$(c)_{\sig_{p+1}}\, J_{p}$ below the merging point $(\Sig_{i})^{pd_{i+1}(\rho)}\, K_{m_{k}}$ 
with $\sig_{p+1}\prec_{i}pd_{i+1}(\rho) \spand i\in In(\sig_{p+1})$.
Then set 
\bdes
\item[Case 1] $rg_{i}(\rho):=\sig_{q}:=rg_{i}(\sig_{p+1})$.
\edes
Otherwise set
\bdes
\item[Case 2] $rg_{i}(\rho)=pd_{i}(\rho)$.
\edes

In general we have the following fact. 
\bprp\label{prp:N5.1}
Let $\calc=J_{0},\ldots,J_{n-1}$ be a chain leading to a $(c)^{\sig_{n-1}}_{\sig_{n}}\, J_{n-1}$. Each $J_{p}$ is a rule $(c)^{\sig_{p}}_{\sig_{p+1}}$ with $\sig_{0}=\pi$. Suppose that the chain passes through the left side of a $(\Sig_{j})^{\sig_{p}}\, K$ for a $p$ with $0<p<n$ and a $j\geq i$ so that $J_{p-1}$ is in the left upper part of $K$ and $J_{p}$ is below $K$. Then $\sig_{n}\prec_{i}\sig_{p}$ and if further $N-1\neq i\in In(\sig_{n})$, then $\sig_{q}=rg_{i}(\sig_{n})\preceq_{i}\sig_{p}$, cf. the figure in Proposition \ref{prp:N4.1}.
\eprp
Let us expalin this Proposition \ref{prp:N5.1} using the new chain 
$\calc_{\rho}=J_{0},\ldots,J_{n-1},J_{n}$ leading to $(c)^{\sig}_{\rho}\, J_{n}$, cf. $Fig.4$. 
When a $(\Sig_{j+1})^{\sig_{t}}\, K^{j+1}\, (0<t\leq n)$ in the new chain $\calc_{\rho}$ 
is to be resolved,  a $(\Sig_{j})^{\sig_{s}}\, K^{j}$ is introduced at a point below $K^{j+1}$. 
The point and $s\geq t$ is determined as the lowest position as far as we can lower a rule 
$(\Sig_{j})^{\sig_{t}}$, cf. Definition \ref{df:5res} in Section \ref{sec:TNc}. 
For example when $K^{j+1}$ is the rule $(\Sig_{i_{m}})^{\sig_{n_{m}+1}}\, K_{m}$ in $Fig.4$, let $m_{1}$ denote the minimal $m_{1}$ such that $i_{m_{1}}<i_{m}$ and 
we introduce a new $(\Sig_{i_{m}-1})^{\sig_{n_{m_{1}+1}}}\, (s=n_{m_{1}+1})$ 
between the rules $(c)_{\sig_{n_{m_{1}+1}}}\, J_{n_{m_{1}}}$ and 
$(\Sig_{i_{m_{1}}})^{\sig_{n_{m_{1}+1}}}\, K_{m_{1}}$. 
Observe that the new $(\Sig_{i_{m}-1})$ together with 
$(\Sig_{i_{m_{2}}})\, K_{m_{2}}\, (m<m_{2}<m_{1})$ by inversion will be merging points 
for the next chain leading to a $(c)^{\rho}$. 

Let us consider the case when the $(\Sig_{i_{m}-1})^{\sig_{n_{m_{1}+1}}}$ is the rule 
$(\Sig_{j})^{\sig_{p}}\, K$ in the Proposition \ref{prp:N5.1}: 
$j=i_{m}-1\spand p=n_{m_{1}+1}$. Also put $pd_{i}(\rho)=\sig_{n_{m(i)}+1}$, where $m(i)$ denotes the minimal $m(i)$ such that $i_{m(i)}<i$. 
Then $i\leq j=i_{m}-1$. 
By Proposition \ref{prp:N5.2} below we see that $i_{m}<i_{m_{3}}$ for any $m_{3}<m$, i.e., any merging point $(\Sig_{i_{m_{3}}})\, K_{m_{3}}$ above $(\Sig_{i_{m}})\, K^{j+1}=K_{m}$ 
has larger index since we are assuming that $K_{m}$ is to be resolved. 
Therefore $m(i)\geq m_{1}$, i.e., the merging point 
$(\Sig_{i_{m(i)}})^{\sig_{n_{m(i)+1}}}\, K_{m(i)}$ determining $pd_{i}(\rho)$ is equal to or 
below the merging point $(\Sig_{i_{m_{1}}})^{\sig_{n_{m_{1}+1}}}\, K_{m_{1}}$. 
In the former case we have $pd_{i}(\rho)=\sig_{n_{m(i)}+1}=\sig_{n_{m_{1}}+1}=\sig_{p}$ 
and hence $\rho\prec_{i}\sig_{p}$. In the latter case we have $i_{m_{3}}\geq i$ for 
$m_{1}\leq m_{3}<m(i)$. Thus we see $\rho\prec_{i}\sig_{p}$ inductively. 
This shows the first half of Proposition \ref{prp:N5.1}.

Now assume $N-1\neq i\in In(\sig_{n})$ and show $rg_{i}(\rho)\preceq_{i}\sig_{p}$. 
Consider the {\bf Case 1}, vz. $\sig_{q}=rg_{i}(\rho)\neq pd_{i}(\rho)$. 
Let $p_{0}$ denote the minimal $p_{0}$ such that 
$\rho\prec_{i}\sig_{p_{0}+1}\prec_{i}pd_{i+1}(\rho)$ and $i\in In(\sig_{p_{0}+1})$. 
By the definition we have $\sig_{q}=rg_{i}(\rho)=rg_{i}(\sig_{p_{0}+1})$.

Let $m(i+1)<m(i)$ denote the number such that $pd_{i+1}(\rho)=\sig_{n_{m(i+1)}+1}$. 
Then by $i\in In(\rho)$ we have $i_{m(i+1)}=i\leq i_{m_{1}}$, i.e., 
$m_{1}\leq m(i+1)<m(i)$ and hence $pd_{i+1}(\rho)\leq\sig_{n_{m_{1}}+1}=\sig_{p}$. 
On the other we have $\rho\prec_{i}\sig_{q}=rg_{i}(\rho)$ by the definition and $\rho\prec_{i}\sig_{p}$ by the first half of the Proposition \ref{prp:N5.1}. 
Hence it suffices to show $\sig_{q}\leq\sig_{p}$ since the set $\{\tau: \rho\prec_{i}\tau\}$ is 
linearly ordered by $\prec_{i}$. 
Now we see $\sig_{q}=rg_{i}(\rho)=rg_{i}(\sig_{p_{0}+1})\preceq_{i}pd_{i+1}(\rho)$ 
inductively, i.e., by using Proposition \ref{prp:N5.1} for smaller parts. 
Thus we get $\sig_{q}\preceq_{i}pd_{i+1}(\rho)\leq\sig_{p}$. 
This shows the second half of Proposition \ref{prp:N5.1}.

Further we have the following fact. 
\bprp\label{prp:N5.2}
Let $\calc=J_{0},\ldots,J_{n-1}$ be a chain leading to a $(c)^{\sig_{n-1}}_{\sig_{n}}\, J_{n-1}$. Each $J_{k}$ is a rule $(c)^{\sig_{k}}_{\sig_{k+1}}$ with $\sig_{0}=\pi$. 
Suppose that the chain $\calc$ passes through the left side of a 
$(\Sig_{j})^{\sig_{p}}\, K^{lw}$ for a $p$ with $0<p<n$ so that $J_{p-1}$ is 
in the left upper part of $K^{lw}$ and $J_{p}$ is below $K^{lw}$. 
Let $\cald=I_{0},\ldots,I_{m-1}\, (m\geq p)$ be a chain leading to a 
$(c)^{\sig_{m-1}}_{\sig_{m}}\, I_{m-1}$. 
Each $I_{k}$ is a rule $(c)^{\tau_{k}}_{\tau_{k+1}}$ 
such that $\tau_{k}=\sig_{k}$ for $0\leq k<\min\{n,m\}$. 
Suppose that the chain $\cald$ passes through the left side of a 
$(\Sig_{i})^{\sig_{k}}\, K^{up}$ for a $k$ with $0<k<p$ so that 
$I_{k-1}$ is in the left upper part of $K^{up}$ and $I_{k}$ is below $K^{up}$. 
Further assume the rule $(c)_{\sig_{p}}\, I_{p-1}$ is in the right upper part of $(\Sig_{j})^{\sig_{p}}\, K^{lw}$ and $i\leq j$. 

Then the upper $K^{up}$ foreruns the lower $K^{lw}$, i.e., analyses of $K^{up}$ have to precede ones of $K^{lw}$.
\eprp

Let us explain Proposition \ref{prp:N5.2} by referring $Fig.4$: 
$\calc$ is the new chain $\calc_{\rho}$, $K^{lw}$ is the new 
$(\Sig_{i_{m}-1})^{\sig_{n_{m_{1}+1}}}$ which is resulted from 
$(\Sig_{i_{m}})^{\sig_{n_{m}+1}}\, K_{m}$ with $m=l-1$, i.e., 
the resolved rule $K_{l-1}$ is the lowest merging point. 
Then $K^{lw}$ is a $(\Sig_{i_{m}-1})^{\sig}$ with $m_{1}=l$. 
Further $\cald$ is the chain 
$\calc_{n_{m+1}}=J^{m+1}_{0},\ldots,J^{m+1}_{n_{m}},J_{n_{m}+1},\ldots,J_{n_{m+1}}$ leading to the last member $(c)_{\sig}\, J_{n-1}\, (n-1=n_{m+1}=n_{l})$ of the series $\calr$. Then the last member $(c)_{\sig}\, J_{n-1}$ is in the right upper part of 
$(\Sig_{i_{m}-1})^{\sig}\, K^{lw}$. Let $I$ be a rule $(\Sig_{i+1})^{\tau}$ 
such that the chain $\cald$ passes through its right side. 
Suppose the rule $I$ in the chain $\cald$ is resolved and produces a 
$(\Sig_{i})^{\sig_{k}}\, K^{up}$ for a $k$ with $0<k<n$ so that the chain $\cald$ passes through the left side of $K^{up}$. 

{\footnotesize
\[
\deduce
{\hskip-0.3cm
 \infer[(c)^{\sig}_{\rho}\, J_{n}]{\Gam_{n}'}
 {
  \infer*[\calc]{\Gam_{n}}
  {}
 }
}
{\hskip-1.5cm
 \infer[(\Sig_{i_{m}-1})^{\sig}\, K^{lw}]{\Phi,\Psi}
 {
  \deduce
  {\hskip0.8cm
   \infer*{\Phi,\lnot C^{\sig}}
   {
    \infer[(c)^{\sig_{n-1}}_{\sig}\, J_{n-1}]{\Gam_{n-1}',\lnot C^{\sig}} 
    {
     \infer*[\cald,\calc]{\Gam_{n-1},\lnot C^{\sig_{n-1}}}
     {}
     }
    }
   }
  {\hskip-1.0cm
   \infer[(\Sig_{i_{m}})^{\sig_{n_{m}+1}}\, K_{m}]{\Phi_{m},\Psi_{m},\lnot C^{\sig_{n_{m}+1}}}
   {\hskip2.0cm
    \infer*[\cald]{\Phi_{m},C_{m}^{\sig_{n_{m}+1}}}
    {
     \infer[(\Sig_{i+1})^{\tau}]{\Pi,\Lam,C_{m}^{\tau}}
     {
      \infer*{\Pi,\lnot B}{}
     &
      \infer*[\cald]{B,\Lam,C_{m}^{\tau}}{}
      }
     }
   &
   \hskip-1.0cm
    \infer*[\calc]{\lnot C_{m}^{\sig_{n_{m}+1}},\Psi_{m},\lnot C^{\sig_{n_{m}+1}}}{}
   }
  }
 &
  \deduce
   {\hskip0.0cm
    \infer*{C^{\sig},\Psi}
    {
     \infer[(c)^{\sig_{n-1}}_{\sig}\, J_{n-1}]{C^{\sig},\Gam_{n-1}'} 
     {
      \infer*[\cald]{C^{\sig_{n-1}},\Gam_{n-1}}
      {
       \infer{C^{\sig_{n_{m}+1}},\Phi_{m},\Psi_{m}}
       {}
      }
     }
    }
   }
   {\hskip-0.9cm
    \infer*{C^{\sig_{n_{m}+1}},\Phi_{m}}
    {
     \infer[(\Sig_{i+1})^{\tau}\, I]{C^{\tau},\Pi,\Lam}
     {
      \infer*{\Pi,\lnot B}{}
     &
      \infer*[\cald]{B,C^{\tau},\Lam}{}
     }
    }
   }
 }
}
\]
}
\hspace*{\fill} $Fig.5$
\\
where $\lnot A_{m}\equiv C_{m}^{\sig_{n_{m}+1}}\equiv \fal x<\sig_{n_{m}+1}C_{0}(x)$ and $C^{\sig_{n_{m}+1}}\equiv C_{0}(\alp)$ for an $\alp<\sig_{n_{m}+1}$.

\[
\hskip3.0cm
\deduce
{\hskip-2.3cm
 \infer[(c)^{\sig}_{\rho}\, J_{n}]{\Gam_{n}'}
 {
  \infer*[\calc]{\Gam_{n}}
  {}
 }
}
{\hskip-3.5cm
 \infer[(\Sig_{i_{m}-1})^{\sig}\, K^{lw}]{\Phi,\Psi}
 {
  \infer*{\Phi,\lnot C^{\sig}}
  {
   \infer[(c)^{\sig_{n-1}}_{\sig}\, J_{n-1}]{\Gam_{n-1}',\lnot C^{\sig}} 
   {
    \infer*[\cald,\calc]{\Gam_{n-1},\lnot C^{\sig_{n-1}}}
    {}
    }
   }
 &
  \deduce
   {\hskip0.0cm
    \infer*{C^{\sig},\Psi}
    {
     \infer[(c)^{\sig_{n-1}}_{\sig}\, J_{n-1}]{C^{\sig},\Gam_{n-1}'} 
     {
      \infer*[\cald]{C^{\sig_{n-1}},\Gam_{n-1}}
      {}
     }
    }
   }
   {\hskip-1.9cm
    \infer[(\Sig_{i})^{\sig_{k}}\, K^{up}]{C^{\sig_{k}},\Lam_{1},\Pi_{1}}
    {
     \infer*[\cald]{C^{\sig_{k}},\Lam_{1},B_{1}^{\sig_{k}}}
     {
      \infer[(\Sig_{i+1})^{\tau}\, I]{C^{\tau},\Pi,\Lam,B_{1}^{\tau}}
      {
       \infer*{\Pi,\lnot B}{}
      &
       \infer*[\cald]{B,C^{\tau},\Lam,B_{1}^{\tau}}{}
       }
      }
     &
      \infer*{\lnot B_{1}^{\sig_{k}},\Pi_{1}}
      {
       \infer{\lnot B_{1}^{\tau},\Pi,\Lam}
       {
        \infer*{\Pi,\lnot B_{1}^{\tau}}{}
        }
       }
     }
    }
 }
}
\]
\hspace*{\fill} $Fig.6$

We show, in $Fig.6$, no ancestor of the right cut formula $C^{\sig}$ of $K^{lw}$ is in the right upper part of $K^{up}$ in order to see that $K^{up}$ foreruns $K^{lw}$. It sufices to see that, in $Fig.5$, no ancestor of the right cut formula $C^{\sig}$ of $K^{lw}$ is in the left upper part of the resolved rule $(\Sig_{i+1})^{\tau}\, I$. Any ancestor of the right cut formula $C^{\sig}$ of $K^{lw}$ comes from the left cut formula $\lnot A_{m}\equiv C_{m}^{\sig_{n_{m}+1}}$ of $(\Sig_{i_{m}})^{\sig_{n_{m}+1}}\, K_{m}$ and any ancestor of the latter is in the chain $\cald$, which in turn passes through the right side of $(\Sig_{i+1})^{\tau}\, I$. Thus any ancestor of the right cut formula $C^{\sig}$ of $K^{lw}$ is in the right upper part of $I$ in $Fig.5$, a fortiori, in the left upper part of $K^{up}$ in $Fig.6$. This shows Proposition \ref{prp:N5.2}.

For full statements and proofs of Propositions \ref{prp:N5.1}, \ref{prp:N5.2}, see Lemmata \ref{lem:5.3.22.A},
the proviso {\bf (uplw)} in Definition \ref{df:5prf} in Section \ref{sec:TNc} and the case {\bf M7.2} in Section \ref{sec:ml5}.

From Propositions \ref{prp:N5.1}, \ref{prp:N5.2} we see that the conditions $(\cald^{Q}.1)$ for $Od(\Pi_{N})$ in \cite{WienpiN} or Section \ref{sec:WienpiN} are enjoyed with respect to the $Q$ part of $\rho$ as for the case $N=4$. A set-theoretic meaning and a wellfoundedness proof of $Od(\Pi_{N})$ are derived from these conditions on o.d.'s as we saw in \cite{LMPSPeking} and \cite{WienpiN}.

Consider a rule $(\Sig_{j})$ in the chain $\calc_{\rho}$ for $j\geq i\in In(\rho)$ which is below 
$(\Sig_{i_{m(i)}})^{pd_{i}(\rho)}\, K_{m(i)}\, (i_{m(i)}<i)$. 
Then from Proposition \ref{prp:N5.2} we see that analyses of such a 
$(\Sig_{j})$ have to follow ones of the rule
 $(\Sig_{i_{m(i)}})^{pd_{i}(\rho)}\, K_{m(i)}$. 
 Thus when such a reversal happens, the lower rule with greater indices $(j>i_{m(i)})$ is dead and we can ignore it. 
 The o.d. $pd_{i}(\rho)$ and the rule $(c)_{pd_{i}(\rho)}\, J_{n_{m(i)}}$ is the predecessor of the o.d. $\rho$ and the rule $(c)_{\rho}$ 
 with respect to $i$: any member $(c)_{\kap}$ of the chain
  $\calc_{\rho}$ with $\rho<\kap<pd_{i}(\rho)$ is irrelevant to the fact that $pd_{i}(\rho)$ and $rg_{i}(\rho)$ are iterated $\Pi_{i}$-reflecting. But the member may be relevant to $\Pi_{j}$-reflection for $j<i$. 
This motivates the definitions of $In(\rho)$ and $pd_{i}(\rho)$.
A series $\kap_{n}\prec_{i}\kap_{n-1}\prec_{i}\cdots\prec_{i}\kap_{0}$ expesses a possible stepping down for the fact that $\kap_{0}$ is an iterated $\Pi_{i}$-reflecting ordinal. 
Degrees of iterations are measured by an ordinal 
$\nu<\kap^{+}$ with $\kap=rg_{i}(\kap_{0}), \nu=st_{i}(\kap_{0})$ 
(and by predecessors of $rg_{i}(\kap_{0})$) as we saw in  \cite{LMPSPeking} and \cite{WienpiN}. 
Therefore we seach only for o.d.'s $\sig_{p+1}$ with $\rho\prec_{i}\sig_{p+1}$ in determing the o.d. $rg_{i}(\rho)=rg_{i}(\sig_{p+1})$.

In the {\bf Case 1} the reason why we chose $\sig_{q}$ as the uppermost one is explained by Propositions \ref{prp:N5.1}, \ref{prp:N5.2} as in the case $N=4$.
\\

Now details follow.

\section{The theory $\mbox{{\rm T}}_{N}$ for $\Pi_{N}$-reflecting ordinals}\label{sec:5}
In this section a theory $\mbox{{\rm T}}_{N}$ of $\Pi_{N}$-reflecting ordinals is defined.

Let $T_{0}$ denote the base theory defined in \cite{ptMahlo}. 
$\calL_{1}$ denotes the language of $T_{0}$. 
Recall that $\calL_{1}=\calL_{0}\cup\{R^{\cala},R^{\cala}_{<}:\cala \mbox{ is a } \Del_{0} \mbox{ formula in } \calL_{0}\cup\{X\}\}$ 
with $\calL_{0}=\{0, 1, +, -, \cdot, q, r, \max, j,( )_0 ,( )_1,=,<\}$. $R^{\cala},R^{\cala}_{<}$ are predicate constants for inductively defined predicates. 
The axioms and inference rules in $T_{0}$ are designed for this language $\calL_{1}$.

The {\it language\/}  $\calL(\mbox{{\rm T}}_{N})$ of the theory $\mbox{{\rm T}}_{N}$ is defined to be $\calL_1\cup\{\Ome\}$ with an individual constant $\Ome$.

The {\it axioms\/} of $\mbox{{\rm T}}_{N}$ are the same as for the theory $T_{3}$ in \cite{ptpi3}, i.e., 
are obtained from those of $T_{22}$ in \cite{ptMahlo} by deleting the axiom
$\Gam,Ad(\Ome).$ Thus the axioms $\Gam,\Lam_{f}$ for the closure of $\Ome$ under the function $f$ in $\calL_{0}$ are included as mathematical axiom in $\mbox{{\rm T}}_{N}$.

The {\it inference rules\/} in $\mbox{{\rm T}}_{N}$ are obtained from $T_{0}$ by adding the following rules $(\Pi_{N}\mbox{{\rm -rfl}})$ and $(\Pi^\Ome_2\mbox{{\rm -rfl}})$.
\[
\infer[(\Pi_{N}\mbox{{\rm -rfl}})]{\Gam}{\Gam,A & \lnot \exi z(t_0<z \land A^z),\Gam }
\]
where $A\equiv\fal x_{N}\exi x_{N-1}\cdots Q x_{1}B(x_{N},x_{N-1},\ldots,x_{1},t_{0})$ is a $\Pi_{N}$ formula. 
\[
\infer[(\Pi_2^\Ome\mbox{{\rm -rfl}})]{\Gam}{\Gam,A^{\Ome} & \lnot \exi z(t<z<\Ome\land A^z),\Gam & \Gam,t<\Ome}
\]
where $A\equiv\fal x\exi yB(x,y,t)$ is a $\Pi_{2}$ formula. 

Concepts related to proof figures, principal or auxiliary formulae, pure variable condition, branch, etc. are defined exactly as in Section 2 of \cite{ptMahlo}.

\section{The system $Od(\Pi_{N})$ of ordinal diagrams}\label{sec:WienpiN}

In this section first let us recall briefly the system $Od(\Pi_{N})$ of ordinal diagrams (abbreviated by o.d.'s) in \cite{WienpiN}.

Let $0,\vphi,\Omega,+,\pi$ and $d$ be distinct symbols. Each o.d. in $Od(\Pi_{N})$ is a finite sequence of these symbols. $\vphi$ is the Veblen function. $\Ome$ denotes the first recursively regular ordinal $\ome^{CK}_{1}$ and $\pi$ the first $\Pi_{N}$-reflecting ordinal.

The set $Od(\Pi_{N})$ is classified into subsets $R, SC, P$ according to the intended meanings of o.d.'s. $P$ denotes the set of additive principal numbers, $SC$ the set of strongly critical numbers and $R$ the set of recursively regular ordinals (less than or equal to $\pi$).
If $\pi>\sig\in R$, then $\sig^{+}$ denotes the next recursively regular diagram to $\sig$.

Recall that $K\alp$ denotes the finite set of o.d.'s defined as follows.
 \begin{enumerate}
 \item $K0 =\emptyset$.
 \item $K (\alpha_1+\cdots +\alpha_n)=\bigcup \{K \alpha_i:1\leq i\leq n \}$
 \item $K \vphi\alpha\bet=K \alpha\cup K \bet$
 \item $K\alp=\{\alp\}$ otherwise, i.e., $\alp\in SC $.
 \eenu

\bdf\label{def:DB}
\benu
\item $\cald_{\sig}(\alp)\incl\cald_{\sig}$.
 \benu
 \item $\cald_{\sig}(\alp)=\emptyset$ {\rm if} $\alp\in\{0,\Ome,\pi\}$.
 \item $\cald_{\sig}(\alp)=\cald_{\sig}(K\alp)$ {\rm if} $\alp\not\in SC$.
 \item {\rm If} $\alp\in\cald_{\tau}$,
 \[\cald_{\sig}(\alp)=\left\{
 \begin{array}{ll}
 \cald_{\sig}(\{\tau\}\cup c(\alp)) & \mbox{{\rm if} } \tau>\sig \\
 \{\alp\}\cup\cald_{\sig}(c(\alp)) & \mbox{{\rm if} } \tau=\sig \\
 \cald_{\sig}(\tau) & \mbox{{\rm if} } \tau<\sig
 \end{array}
 \right.
 \]
 \eenu
\item $\calb_{\sig}(\alp)=\max\{b(\bet):\bet\in\cald_{\sig}(\alp)\}$.
\item $\calb_{>\sig}(\alp)=\max\{\calb_{\tau}(\alp):\tau>\sig\}$.
\eenu
\edf

For an o.d. $\alp$ set
\[
\alp^{+}=\min\{\sig\in R\cup\{\infty\}:\alp<\sig\}
.\]

For $\sig\in R$, $\cald_{\sig}\incl SC$ denotes the set of o.d.'s of the form 
$\rho=d^{q}_{\sig}\alp$ with a (possibly empty) list $q$, where the following condition has to be met:
\beqn\label{eq:Odmu}
\calb_{>\sig}(\{\sig,\alp\}\cup q)<\alp 
\eeqn
$\alp$ is the {\it body\/} of $d^{q}_{\sig}\alp$. 

If $q$ is not empty, then $d_{\sig}^{q}\alp\in\cald^{Q}$ by definition.
Its $Q$ {\it part\/} $Q(d^{q}_{\sig}\alp)=q=\ovl{\nu\kap\tau j}$ denotes a sequence of quadruples $\nu_m\kap_m\tau_m j_m$ of length $l+1 \, (0\leq l)$ such that  
\benu
\item 
$2\leq j_0<j_1<\cdots <j_l=N-1,$
\item 
$\kap_l=\pi, \kap_m\in R\mid\pi \, (m<l) \spand \sig\preceq\kap_m \, (m\leq l),$
\item 
$\nu_{l}\in Od(\Pi_{N})$,
 \beqn\label{eq:stboundN-1}
 \sig=\pi \Rarw \nu_{l}\leq\alp
 \eeqn
 {\rm and}
 \beqn\label{eq:stbound}
 m<l \Rarw \nu_m<\kap_{m}^{+}
 \eeqn
\item 
$\tau_0=\sig, \tau_m\in\{\pi\}\cup\cald^{Q}, \sig\preceq \tau_m\, (m\leq l)$
  {\rm and}
 \beqn\label{eq:5pred}
 \tau_{l}=\pi \Rarw \sig=\pi
 \eeqn
\eenu
From $q=Q(\rho)$ define
\benu
\item $in_j(\rho)=st_j(\rho)rg_j(\rho)$ (a pair) and $pd_j(\rho)$: Given $j$ 
with $2\leq j<N$, put $m=\min\{m\leq l:j\leq j_m\}$.
\item $pd_j(\rho)=\tau_m$.
\item  $\exi m\leq l(j=j_m)$: Then $st_j(\rho)=\nu_m, \, rg_j(\rho)=\kap_m$.
\item Otherwise: $in_j(\rho)=in_j(pd_j(\rho))=in_j(\tau_m)$. 
If $in_j(\tau_m)=\emptyset$, then set $st_j(\rho)\uarw, rg_j(\rho)\uarw$.
\item $In(\rho)=\{j_m:m\leq l\}$.
\eenu
 
Observe that
\beqn\label{eq:Qrest}
\pi<\bet\in q=Q(\rho)\Rarw \bet=\nu_{l}=st_{N-1}(\rho)
\eeqn

The relation $\alp\prec_{i}\bet$ is the transitive closure of the relation 
$pd_{i}(\alp)=\bet$.

In \cite{WienpiN} we impose several conditions on a diagram of the form 
$\rho=d^{q}_{\sig}\alp$ to be in $Od(\Pi_{N})$. 
For $\alp\in Od(\Pi_{N}), q\incl Od(\Pi_{N}) \, \& \,  \sig\in R\setm\{\Ome\}$, 
$\rho=d^{q}_{\sig}\alp\in Od(\Pi_{N})$ if the following conditions are fulfilled
besides (\ref{eq:Odmu});
\bdes
\item[$(\cald^{Q}.1)$] Assume $i\in In(\rho)$. Put $\kap=rg_i(\rho)$. Then \label{cnd:stn}
 \bdes
 \item[$(\cald^{Q}.11)$] $in_i(rg_i(\rho))=in_i(pd_{i+1}(\rho))$, $rg_i(\rho)\preceq_{i}pd_{i+1}(\rho)$ and 
 $pd_i(\rho)\neq pd_{i+1}(\rho)$ if $i<N-1$. \\
 Also $pd_i(\rho)\preceq_{i}rg_i(\rho)$ for any $i$.\label{cnd:stn.1}
 \item[$(\cald^{Q}.12)$] One of the following holds: \label{cnd:stn.2}
  \bdes
  \item[$(\cald^{Q}.12.1)$] 
  $rg_i(\rho)=pd_i(\rho)
  \spand \calb_{>\kap}(st_i(\rho))<b(\alp_{1})$ {\rm with} $\rho\preceq\alp_{1}\in\cald_{\kap}$. \label{cnd:stn.21}
  
  \item[$(\cald^{Q}.12.2)$]  
  $rg_i(\rho)=rg_i(pd_i(\rho)) \spand st_i(\rho)<st_i(pd_i(\rho))$.
 \label{cnd:stn.22}
 
  \item[$(\cald^{Q}.12.3)$]  
  $rg_i(pd_i(\rho))\prec_i\kap \spand \\
\fal \tau(rg_i(pd_i(\rho))\preceq _i\tau\prec _i\kap\to rg_i(\tau)\preceq _i\kap) \spand 
st_i(\rho)<st_i(\sig_1)$ with
 \[\sig_1=\min\{\sig_1:rg_i(\sig_1)=\kap\spand pd_i(\rho)\prec _i\sig_1\prec _i\kap\}\]
 and such a $\sig_1$ exists. \label{cnd:stn.23}
  \edes
 \edes
\item [$(\cald^{Q}.2)$]
   \beqn \label{cnd:Kst}
   \fal\kap\leq rg_i(\rho)(K_{\kap}st_i(\rho)<\rho)
   \eeqn
    {\rm for} $i\in In(\rho)$.

\edes

We set $Q(d_{\sig}\alp)=\emptyset$, i.e., $d_{\sig}^{\emptyset}\alp=d_{\sig}\alp$.

The order relation $\alp<\bet$ on $\cald_{\sig}$ is defined through finite sets $K_{\tau}\alp$ for $\tau\in R, \alp\in Od(\Pi_{N})$, and the latter is defined through the relation $\alp\prec\bet$, which is the transitive closure of the relation $\alp\in\cald_{\bet}$. Thus $\alp\prec_{2}\bet\Lrarw\alp\prec\bet$.


For $\rho=d^{q}_{\tau}\alp$ $c(\rho)=\{\tau,\alp\}\cup q$ and
\[
K_{\sig}\rho=\left\{
        \begin{array}{ll}
         K_{\sig}(\{\tau\}\cup c(\rho))=K_{\sig}\{\tau,\alp\}\cup q, & \sig<\tau \\
         K_{\sig}\tau, & \tau<\sig \, \& \, \tau\not\preceq\sig
        \end{array}
        \right.
        \]


The following Proposition \ref{lem:5.4} is shown in \cite{WienpiN}.

\bprp\label{lem:5.4}
\benu
\item 
The finite set $\{\tau:\sig\prec_{i}\tau\}$ is linearly ordered by $\prec_{i}$.
\label{lem:5.3.1}
\\
In the following assume $\kap=rg_i(\rho)\darw$.
\item 
$\rho\prec _{i}rg_i(\rho)$.
\label{lem:5.4.1}

\item 
$\rho\prec_i\sig\prec_i\tau \spand in_i(\rho)=in_i(\tau) \Rarw in_i(\rho)=in_i(\sig)$.
\label{lem:5.4.6}

\item 
$\rho\prec_i\tau\prec_{i}rg_i(\rho) \Rarw rg_i(\tau)\preceq_{i}rg_i(\rho)$.
\label{lem:5.4.8}
\eenu
\eprp

\bdf\label{df:ksige-1}
{\rm For o.d.'s} $\alp,\sig$ {\rm with} $\sig\in R$,
\[\calk_{\sig}(\alp):=\max K_{\sig}\alp.\]
\edf

The following lemmata are seen as in \cite{ptMahlo}.




\blem\label{lem:rgbnd}
Suppose $\calb_{>\kap}(\alp_{i})<\alp_{i}$ for $i=0,1$, and $\alp_{0}<\alp_{1}$.
Then
\[
\tau>\kap\Rarw d_{\tau}\alp_{i}\in Od(\Pi_{N})\spand d_{\tau}\alp_{0}<d_{\tau}\alp_{1}
.\]
\elem

\blem\label{lem:OdbndB}
For $\alp,\bet,\sig\in Od(\Pi_{N})$ with $\sig\in R|\pi$ assume $\fal\tau<\pi[\calb_{\tau}(\bet)\leq\calb_{\tau}(\alp)]$, and put $\gam=\max\{\calb_{\pi}(\bet),\calb_{>\sig}(\{\sig,\alp\})\}+\ome^{\bet}$. Then
$\calb_{>\sig}(\{\sig,\gam,\gam+\calk_{\sig}(\alp)\})<\gam$, and hence
(\ref{eq:Odmu}) is fulfilled for
$d_{\sig}\gam, d_{\sig}(\gam+\calk_{\sig}(\alp))\in Od(\Pi_{N})$.
\elem


\section{The system $\mbox{{\rm T}}_{Nc}$}\label{sec:TNc}
In this section we extend $\mbox{{\rm T}}_{N}$ to a formal system $\mbox{{\rm T}}_{Nc}$. 
The {\it universe\/} $\pi(\mbox{{\rm T}}_{N})$ of the theory $\mbox{{\rm T}}_{N}$ is defined to be the o.d. 
$\pi\in Od(\Pi_{N})$. 
The language is expanded so that individual constants $c_{\alp}$ for o.d.'s 
$\alp\in Od(\Pi_{N})\mid\pi$ are included. 
Inference rules $(c)^{\sig}$ are added. 
To each proof $P$ in $\mbox{{\rm T}}_{Nc}$ an o.d. $o(P)\in Od(\Pi_{N})\mid\Ome$ is attached. 
{\it Chains\/} are defined to be a consecutive sequence of rules $(c)$. 
{\it Proofs\/} in $\mbox{{\rm T}}_{Nc}$ defined in Definition \ref{df:5prf} are proof figures 
enjoying some provisos and obtained from given proofs in $\mbox{{\rm T}}_{N}$ 
by operating rewriting steps. 
Some lemmata for proofs are established. 
These are needed to verify that rewrited proof figures enjoy these provisos.

The {\it language\/} $\calL_{Nc}$ of $\mbox{{\rm T}}_{Nc}$ is obtained from the language 
$\calL(\mbox{{\rm T}}_{N})$ by adding individual constants $c_{\alp}$ for each o.d. 
$\alp\in Od(\Pi_{N})$ such that $1<\alp<\pi \spand \alp\neq\Ome$. 
We identify the constant $c_\alp$ with the o.d. $\alp$. 

In what follows $A,B,\ldots$ denote formulae in $\calL_{Nc}$ and $\Gam,\Del,\ldots$ sequents in $\calL_{Nc}$.

The {\it axioms\/} of $\mbox{{\rm T}}_{Nc}$ are obtained from those of $\mbox{{\rm T}}_{N}$ as in \cite{ptMahlo}. 


Complexity measures $\mbox{deg}(A), \mbox{rk}(A)$ of formulae $A$ are defined as in \cite{ptMahlo} by replacing the universe $\pi(\mbox{T}_{22})=\mu$ by $\pi(\mbox{T}_{N})=\pi$.



Also the sets $\Del^{\sig}_{0},\Sig_{i}^{\sig}$ of formulae are defined as in \cite{ptMahlo}.
Recall that for a bounded formula $A$ and a multiplicative principal number 
$\alp\leq\pi$, we have $A\in\Del^{\alp} \Lrarw \mbox{deg}(A)<\alp$.

\bdf\label{df:deg3}
\[
\mbox{{\rm deg}}_{N}(A):=
\left\{
 \begin{array}{ll}
 \mbox{{\rm deg}}(A)+N-1 & \mbox{{\rm if} } A \mbox{ {\rm is a bounded formula}} \\
 \mbox{{\rm deg}}(A) & \mbox{{\rm otherwise}} 
 \end{array}
 \right.
\]
\edf
Note that
\[
\mbox{deg}_{N}(A)\not\in\{\alp+i: i<N-1, \alp<\pi \mbox{ is a limit o.d.}\}
\]

The {\it inference rules\/} of $\mbox{{\rm T}}_{Nc}$ are obtained from those of 
$\mbox{{\rm T}}_{N}$ by adding the following rules $(h)^{\alp}\, (\alp\in\{\alp:\pi\leq\alp<\pi+\ome\}\cup\{0,\Ome\})$, 
$(c\Pi_{2})^{\Ome}_{\alp_{1}}$, 
$(c\Sig_{1})^{\Ome}_{\alp_{1}}$, 
$(c\Pi_{N})^{\sig}_{\tau}$, 
$(c\Sig_{N-1})^{\sig}_{\tau}$ for each 
$\sig\in R\incl Od(\Pi_{N})\spand \sig\neq\Ome$ and 
$(\Sig_{i})^{\sig}$ for each  $\sig\in R\incl Od(\Pi_{N})\spand \sig\not\in\{\Ome,\pi\}$ and $i=1,2,\ldots,N$. 
The rule $(h)^{\alp}$, $(c\Pi_2)^\Ome_{\alp_1}$ and 
$(c\Sig_1)^\Ome_{\alp_1}$ are the same as in \cite{ptpi3}.
We write $(w)$ for $(h)^{0}$.

\benu



 
 
\item
\[\infer[(c\Pi_{N})^{\sig}_{\tau}]{\Gam,A^{\tau}}{\Gam,A^{\sig}}\]
where 
 \benu
 \item 
 $A\equiv\fal x_{N}\exi x_{N-1}\cdots Q x_{1}B$ is a $\Pi_{N}$-sentence with a $\Del^{\tau}$-matrix $B$,
 \item 
 $\tau\in\cald_{\sig}$ with the {\it body\/} $\alp=b(\tau)$ of the rule and
 \item 
 the formula $A^{\tau}$ in the lowersequent is the {\it principal formula\/} of the rule and the formula $A^{\sig}$ in the uppersequent is the {\it auxiliary formula\/} of the rule, resp. Each formula in  $\Gam$ is a {\it side formula\/} of the rule.
 \eenu
 
\item
\[\infer[(c\Sig_{N-1})^{\sig}_{\tau}]{\Gam,\Lam^{\tau}}{\Gam,\Lam^{\sig}}\]
where
 \benu
 \item $\Lam$ is a nonempty set of unbounded $\Pi_{N}$-sentences with 
 $\Del^{\tau}$-matrices.
 
 \item $\tau\in\cald_{\sig}$ with the {\it body\/} $\alp=b(\tau)$ of the rule and
 \item each formula in $\Gam$ is a {\it side formula\/} of the rule.
 \eenu
 
\item
\[\infer[(\Sig_{i})^{\sig}]{\Gam,\Lam}{\Gam,\lnot A^{\sig} & A^{\sig},\Lam}\]
where $1\leq i\leq N$ and $A^{\sig}$ is a genuine $\Sig_{i}^{\sig}$-sentence, i.e., $A^{\sig}\in\Sig_{i}^{\sig}$ and $A^{\sig}\not\in\Pi^{\sig}_{i-1}\cup\Sig^{\sig}_{i-1}$.

$A^{\sig}$ [$\neg A^{\sig}$] is said to be the {\it right\/} [{\it left\/}] {\it cut formula\/} of the rule $(\Sig_i)^{\sig}$, resp.
\eenu


The rules $(c\Pi_2)^{\Ome}$ and $(c\Pi_{N})$ are {\it basic rules\/} but not the rules
 $(h)^{\alp}$, $(c\Sig_{1})^{\Ome}$, $(c\Sig_{N-1})^{\sig}$ and $(\Sig_i)^{\sig}$. 

A {\it preproof\/} in $\mbox{{\rm T}}_{Nc}$ is a proof in $\mbox{{\rm T}}_{Nc}$ in the sense of \cite{ptMahlo}, i.e., a proof tree built from axioms and inference rules in $\mbox{{\rm T}}_{Nc}$.
The underlying tree $\mbox{Tree}(P)$ of a preproof $P$ is a tree of finite sequences of natural numbers
such that each occurrence of a sequent or an inference rule receives a finite sequence.
The root (empty sequence) $(\,)$ is attached to the endsequent,
and in an inference rule
\[
\infer[(r)\, a*(0)]{a:\Gam}{a*(0,0):\Lam_0 & \cdots & a*(0,n):\Lam_n}
\]
where $(r)$ is the name of the inference rule.
Finite sequences are denoted by Roman letters $a,b,c,\ldots, I,J,K,\ldots$.
Roman capitals $I,J,K,\ldots$ denote exclusively inference nodes.
We will identify the attached sequence $a$ with the occurrence of a sequent or an inference rule.





Let $P$ be a preproof and $\gam<\pi+\ome$ an o.d. in $Od(\Pi_{N})$. 
For each sequent $a:\Gam\, (a\in\mbox{Tree}(P))$, we assign the {\it height\/} $\mbox{h}_{\gam}(a;P)<\pi+\ome$ of the node $a$ {\it with the base height\/} $\gam$ in $P$ as in \cite{ptMahlo} except we replace $\pi(\mbox{T}_{22})=\mu$ by $\pi(\mbox{T}_{N})=\pi$ and replace $\mbox{deg}(A)$ by $\mbox{deg}_{N}(A)$.

Then the {\it height\/} $\mbox{{\rm h}}(a;P)$ of $a$ in $P$ is defined to be the height with the base height $\gam=0$:
\[
\mbox{{\rm h}}(a;P):=\mbox{{\rm h}}_{0}(a;P)
.\]

A pair $(P,\gam)$ of a preproof $P$ and an o.d. $\gam$ is said to be {\it height regulated\/}
if it enjoys the conditions in \cite{ptMahlo}, or equivalently in \cite{ptpi3}, Definition 5.4.
For the rules $(\Sig_i)^{\sig}$, this requires the condition:
If $a:\Gam$ is the lowersequent of a rule $(\Sig_i)^{\sig}\, a*(0)\, (1\leq i\leq N)$ in $P$, then 
$\mbox{{\rm h}}_{\gam}(a;P)\leq\sig+i-2$ if $i=N-1,N$.
Otherwise
$\mbox{{\rm h}}_{\gam}(a;P)\leq\sig+i-1$.
 
Therefore for the uppersequent $a*(0,k):\Lam$ of a $(\Sig_i)^{\sig}$ we have
 $\mbox{{\rm h}}_{\gam}(a*(0,k);P)=\sig+i-1$. 
 Note that this implies that there is no nested rules $(\Sig_{i})^{\sig}$, i.e., there is no $(\Sig_{i})^{\sig}$  
 below any $(\Sig_{i})^{\sig}$ for $i\geq N-1$.

A preproof is height regulated iff $(P,0)$ is height regulated.

Let $P$ be a preproof and $\gam<\pi+\ome$.
Assume that $(P,\gam)$ is height regulated. 
Then the o.d $o_{\gam}(a;P)\in O(\Pi_{N})$ assigned to each node $a$ in the underlying tree $\mbox{{\rm Tree}}(P)$ of $P$ is defined exactly as in \cite{ptpi3}.

Furthermore for $\tau\in R\cap Od(\Pi_{N})$,
o.d.'s $B_{\tau,\gam}(a;P), Bk_{\tau,\gam}(a;P)\in O(\Pi_{N})$ are assigned 
to each sequent node $a$ such that $\mbox{{\rm h}}_{\gam}(a;P)\leq\tau\in R$ as in \cite{ptpi3}.
Namely

\beqnarrs
B_{\tau,\gam}(a;P) & := &
\left\{
\begin{array}{l}
\pi\cdot o_{\gam}(a;P) \\
 \mbox{ {\rm if }} \mbox{{\rm h}}_{\gam}(a;P)=\tau=\pi \\
\max\{\calb_{\pi}(o_{\gam}(a;P)),\calb_{>\tau}(\{\tau\}\cup (a;P))\}+\ome^{o_{\gam}(a;P)} \\
 \mbox{ {\rm if }} \mbox{{\rm h}}_{\gam}(a;P)<\pi
\end{array}
\right.
\\
Bk_{\tau,\gam}(a;P) & := & B_{\tau,\gam}(a;P)+\calk_{\tau}(a;P)
\eeqnarrs
$B_{\tau}(a;P)$ [$Bk_{\tau}(a;P)$] denotes $B_{\tau,0}(a;P)$ [$Bk_{\tau,0}(a;P)$], resp.

Then propositions and lemmata (Rank Lemma 7.3, Inversion Lemma 7.9, etc.) in Section 9 of \cite{ptMahlo} 
and Replacement Lemma 5.15 in \cite{ptpi3} hold also for $\mbox{T}_{Nc}$.

Lemma \ref{lem:OdbndB} yields $o_{\gam}(a;P)\in Od(\Pi_{N})$ for each node $a\in\mbox{Tree}(P)$ if 
$(P,\gam)$ is height regulated and $\gam<\pi+\ome$.

\bdf\label{df:5branch}
{\rm Let} $\calt$ {\rm be a branch in a preproof} $P$ {\rm and} $J$ {\rm a rule} $(\Sig_i)^{\sig}$.
\benu
\item Left branch: $\calt$ {\rm is a} left branch {\rm of} $J$ {\rm if}
 \benu
 \item $\calt$ {\rm starts with a lowermost sequent} $\Gam$ {\rm such that} $h(\Gam)\geq\pi$, 
 \item {\rm each sequent in} $\calt$ {\rm contains an ancestor of the left cut formula of} $J$ {\rm and}
 \item $\calt$ {\rm ends with the left uppersequent of} $J${\rm.}
 \eenu
\item Right branch: $\calt$ {\rm is a} right branch {\rm of} $J$ {\rm if}
 \benu
 \item $\calt$ {\rm starts with a lowermost sequent} $\Gam$ {\rm such that} $\Gam$ {\rm is a lowersequent of a basic rule whose principal formula is an ancestor of the right cut formula of} $J$ {\rm and} 
 \item $\calt$ {\rm ends with the right uppersequent of} $J${\rm.}
 \eenu
\eenu
\edf


{\it Chains\/} in a preproof are defined as in Definition 6.1 of \cite{ptpi3} 
when we replace $((c\Pi_{3}),(\Sig_{3}))$, $((c\Sig_{2}), (\Sig_{2}))$ by
$((c\Pi_{N}),(\Sig_{N}))$, $((c\Sig_{N-1}), (\Sig_{N-1}))$.
For definitions related to chains such as {\it starting with\/}, {\it top\/}, {\it branch\/} of a chain,
{\it passing through\/}, see Definition 6.1 of \cite{ptpi3}.
Also {\it rope sequence\/} of a rule, the {\it end\/} of a rope sequence and the {\it bar\/}
of a rule are defined as in Definition 6.2 of \cite{ptpi3}.
Moreover a {\it chain analysis\/} for a preproof together with the {\it bottom\/} of a rule
is defined as in Definition 6.3 of \cite{ptpi3}.

\bdf\label{df:5qorign}$Q$ part of a chain and the $i$-origin.
\benu
\item {\rm Let} $\calc=J_0,J_0',\ldots,J_n,J_n'$ {\rm be a chain starting with a} $(c)_{\sig}\, J_n${\rm . Put}
 \benu
 \item $In(\calc):=In(J_n):=In(\sig)$.
 \item $in_{i}(\calc):=in_{i}(J_n):=in_{i}(\sig)$ {\rm for} $2\leq i<N$.
 \item $st_{i}(\calc):=st_{i}(J_n):=st_{i}(\sig), \: rg_{i}(\calc):=rg_{i}(J_n):=rg_{i}(\sig)$\\
 {\rm where} $st_{i}(\calc)\uparrow \spand rg_{i}(\calc)\uparrow$ {\rm if} $st_{i}(\sig)\uparrow \spand rg_{i}(\sig)\uparrow$.
 \item $J_k$ {\rm is the} $i$-origin {\rm of the chain} $\calc$ {\rm or the rule} $J_n$ {\rm if} $J_k$ {\rm is a rule} $(c)^{\kap}$ {\rm with} $\kap=rg_{i}(\sig)\downarrow$.
 \item $J_k$ {\rm is the} $i$-predecessor {\rm of} $J_n${\rm , denoted by} $J_k=pd_{i}(J_n)$ {\rm or} \\
 $i$-predecessor {\rm of the chain} $\calc${\rm , denoted by} $J_k=pd_{i}(\calc)$ {\rm if} $J_k$ {\rm is a rule} $(c)_{\rho}$ {\rm with} $\rho=pd_{i}(\sig)$. 
 \eenu
\eenu
\edf

\bdf\label{df:5knot}Knot and rope.\\
{\rm Assume that a chain analysis for a preproof} $P$ {\rm is given and by a chain we mean a chain in the chain analysis.}
\benu
\item \label{df:5knot.1} $i$-knot: {\rm Let} $K$ {\rm be a rule} $(\Sig_i)^{\sig}\, (1\leq i\leq N-2)$. {\rm We say that} $K$ {\rm is an} $i$-knot {\rm if there are an uppermost rule} $(c)^{\sig}\, J_{lw}$ {\rm below} $K$ {\rm and a chain} $\calc$ {\rm such that} $J_{lw}$ {\rm is a member of} $\calc$ {\rm and} $\calc$ {\rm passes through the left side of} $K$.

{\rm The rule} $J_{lw}$ {\rm is said to be the} lower rule {\rm of the} $i${\rm -knot} $K$. {\rm The member} $(c)_{\sig}\, J_{ul}$ {\rm of the chain} $\calc$ {\rm is the} upper left rule {\rm of} $K$ {\rm and a rule} $(c)_{\sig}\, J_{ur}$ {\rm which is above the right uppersequent of} $K$ {\rm is an} upper right rule {\rm of} $K$ {\rm if such a rule} $(c)_{\sig}\, J_{ur}$ {\rm exists.}
\[
\infer*{}
{
 \infer[\mbox{{\rm uppermost}}\: (c)^{\sig}\, J_{lw}\in\calc]{\Del'}
 {
  \infer*{\Del}
  {
   \infer[(\Sig_i)^{\sig}\, K]{\Gam,\Lam}
   {
    \infer*[\calc]{\Gam,\lnot A^{\sig}}{}
   &
    \infer*{A^{\sig},\Lam}{}
   }
  }
 }
}
\]
\item\label{df:5knot.2} {\rm A rule is a} knot {\rm if it is an} $i${\rm -knot for some} $i> 1$.\\
{\bf Remark}{\rm . Note that a} $1${\rm -knot} $(\Sig_1)$ {\rm is not a knot by definition.}
\item\label{df:5knot.3} {\rm Let} $K$ {\rm be a knot,} $J_{lw}$ {\rm the lower rule of} $K$ {\rm and} $J_{ur}$ {\rm an upper right rule of} $K${\rm . Then we say that} $K$ {\rm is a} knot of $J_{ur}$ and $J_{lw}$.
\item {\rm Let} $\calc_n=J_0,\ldots,J_n$ {\rm be a chain starting with} $J_n$ {\rm and} $K$ {\rm a knot.} $K$ {\rm is a} knot for the chain $\calc_n$ {\rm or} the rule $J_n$ {\rm if}
 \benu
 \item {\rm the lower rule} $J_{lw}$ {\rm of} $K$ {\rm is a member} $J_k\, (k<n)$ {\rm of} $\calc_n$,
 \item $\calc_n$ {\rm passes through the right side of} $K${\rm , and}
 \item {\rm for any} $k<n$ {\rm the chain} $\calc_k$ {\rm starting with} $J_k$ {\rm does not pass through the right side of} $K$.
 \eenu
{\rm The knot} $K$ {\rm is a merging rule of the chain} $\calc_n$ {\rm and the chain} $\calc_k$ {\rm starting with the lower rule} $J_{lw}=J_k$.
\[
\deduce
{\hskip4.6cm
\infer[J_n]{\Del_n'}
{
 \infer*{\Del_n}
 {
  \infer[\mbox{{\rm uppermost }} (c)^{\sig}\, J_{lw}=J_k\in\calc_n]{\Del_k'}
  {
   \infer*{\Del_k}
   {}
  }
 }
}
}
{\hskip-0.5cm
   \infer[(\Sig_i)^{\sig}\, K]{\Gam,\Lam}
   {
    \infer*[\calc_k]{\Gam,\lnot A^{\sig}}{}
   &
    \infer*[\calc_n]{A^{\sig},\Lam}{}
   }
 }
\]
\item \label{df:5knot.4} {\rm A series} $\calr_{J_0}=J_0,\ldots,J_{n-1}\, (n\geq 1)$ {\rm of rules} $(c)$ {\rm is said to be the} rope starting with $J_0$ {\rm if there is an increasing sequence of numbers (uniquely determined)}
\beqn\label{eqn:5rope}
0\leq n_0<n_1<\cdots<n_l=n-1\, (l\geq 0)
\eeqn
{\rm for which the following hold:}
 \benu
 \item {\rm each} $J_{n_m}$ {\rm is the bottom of} $J_{n_{m-1}+1}$ {\rm for} $m\leq l\, (n_{-1}=-1)$,
 \item {\rm there is an uppermost knot} $K_m$ {\rm such that} $J_{n_m}$ {\rm is an upper right rule and} $J_{n_{m}+1}$ {\rm is the lower rule of} $K_m$ {\rm for} $m<l${\rm , and}
 \item {\rm there is no knot whose upper right rule is} $J_{n_l}=J_{n-1}$.
 \eenu
 {\rm We say that the rule} $J_{n-1}$ {\rm is the} edge of the rope $\calr_{J_0}$ {\rm or} the edge of the rule $J_0$.
 \\
{\rm For a rope the increasing sequence of numbers (\ref{eqn:5rope}) is called the} knotting numbers {\rm of the rope.}\\
{\bf  Remark}. {\rm These knots} $K_m$ {\rm are uniquely determined for a} proof {\rm defined below.}
\item\label{df:5knot.5} {\rm Let} $K_{-1}$ {\rm be an} $i_{-1}${\rm -knot} $(i_{-1}\geq 1)$ {\rm and} $J_0$ {\rm the lower rule of} $K_{-1}$. {\rm The} left rope ${}_{K_{-1}}\calr$ of $K_{-1}$ {\rm is inductively defined as follows:}
 \benu
 \item {\rm Pick the lowermost rule} $(c)\, J_{n_0}$ {\rm such that the chain} $\calc$ {\rm starting with} $J_{n_0}$ {\rm passes through the left side of the} $i_{-1}${\rm -knot} $K_{-1}$ {\rm and} $J_0$ {\rm is a member of} $\calc$. {\rm Let} ${}_{0}\calr=I_0,\ldots,I_q$ {\rm be the part of the chain} $\calc$ {\rm with} $J_0=I_0\spand J_{n_0}=I_q$.
 \item {\rm If there exists an uppermost knot} $K_{0}$ {\rm such that} $J_{n_0}$ {\rm is an upper right rule of} $K_0${\rm , then} ${}_{K_{-1}}\calr$ {\rm is defined to be a concatenation :}
\[{}_{K_{-1}}\calr={}_{0}\calr^{\frown}{}_{K_{0}}\calr\]
{\rm where} ${}_{K_{0}}\calr$ {\rm denotes the left rope of} $K_0$.
 \item {\rm Otherwise. Set:}
 \[{}_{K_{-1}}\calr={}_{0}\calr\]
 \eenu
{\rm Therefore for the left rope} ${}_{K_{-1}}\calr=J_0,\ldots,J_{n-1}$ {\rm of} $K_{-1}$ {\rm there exists a uniquely determined increasing sequence of numbers (\ref{eqn:5rope}) such that:}
 \benu
 \item {\rm each} $J_{n_m}$ {\rm is the lowermost rule} $(c)$ {\rm such that the chain} $\calc$ {\rm starting with} $J_{n_m}$ {\rm passes through the left side of the} $i_{m-1}${\rm -knot} $K_{m-1}$ {\rm and} $J_{n_{m-1}+1}$ {\rm is a member of} $\calc\, (n_{-1}=-1)$ {\rm for} $m\leq l$,
  \item {\rm there is an} $i_m${\rm -knot} $K_m\, (i_m>1)$ {\rm such that} $J_{n_m}$ {\rm is an upper right rule and} $J_{n_{m}+1}$ {\rm is the lower rule of} $K_m$ {\rm for} $m<l${\rm , and}
 \item {\rm there is no knot whose upper right rule is} $J_{n_l}=J_{n-1}${\rm . (} $K_{-1}$ {\rm is the} $i_{-1}${\rm -knot whose lower rule is} $J_0${\rm .)}
 \eenu
{\rm These numbers (\ref{eqn:5rope}) is called the} knotting numbers {\rm of the left rope and each knot} $K_{m}\, (m<l)$ {\rm a} knot for the left rope.

{\rm By the} left rope ${}_{J_0}\calr$ {\rm of the lower rule} $J_0$ {\rm of} $K_{-1}$ {\rm we mean the left rope} ${}_{K_{-1}}\calr$ {\rm of} $K_{-1}$.
\eenu
\edf

When a rule $(\Sig_{i+1})^{\sig}\, K \, (0<i<N)$ is resolved, we introduce a new rule \\
$(\Sig_{i})^{\sig_{n_{m(i+1)}+1}}$ at a sequent $\Phi$, which is defined to be the {\it resolvent} of $K$ and a $\sig_{n_{m(i+1)}+1}\preceq\sig$ defined as follows.

\bdf\label{df:5res}Resolvent\\
{\rm Let} $K$ {\rm be a rule} $(\Sig_{i+1})^{\sig}\, (0<i<N)$. 
{\rm The} resolvent {\rm of the rule} $K$ {\rm is a sequent} $a:\Phi$ {\rm defined as follows: let} $K'$ {\rm denote the lowermost rule} $(\Sig_{i+1})^{\sig}$ {\rm below or equal to} $K$ {\rm and} $b:\Psi$ {\rm the lowersequent of} $K'$.
\bdes
\item[Case 1]{\rm The case when there exists an} $(i+1)${\rm -knot} $(\Sig_{i+1})^{\sig}$ {\rm which is between an uppersequent of} $K$ {\rm and} $b:\Psi${\rm : Pick the uppermost such knot} $(\Sig_{i+1})^{\sig}\, K_{-1}$ {\rm and let} ${}_{K_{-1}}\calr=J_0,\ldots,J_{n-1}$ {\rm denote the left rope of} $K_{-1}$. {\rm Each} $J_p$ {\rm is a rule} $(c)^{\sig_p}_{\sig_{p+1}}$. {\rm Let} 
\beqn
\renewcommand{\theequation}{\ref{eqn:5rope}}
0\leq n_0<n_1<\cdots<n_l=n-1\, (l\geq 0)
\eeqn
\addtocounter{equation}{-1}
{\rm be the knotting numbers of the left rope} ${}_{K_{-1}}\calr$ {\rm and} $K_m$ {\rm an} $i_m${\rm -knot} $(\Sig_{i_m})^{\sig_{n_{m}+1}}$ {\rm of} $J_{n_m}$ {\rm and} $J_{n_{m}+1}$ {\rm for} $m<l$. {\rm Put}
\beqn\label{eqn:m(i+1)}
m(i+1)=\max\{m:0\leq m\leq l\spand \fal p\in [0,m)(i+1\leq i_p)\}
\eeqn
{\rm Then the resolvent} $a:\Phi$ {\rm is defined to be the uppermost sequent} $a:\Phi$ {\rm below} $J_{n_{m(i+1)}}$ {\rm such that} $h(a;P)<\sig_{n_{m(i+1)}+1}+i$.
\item[Case 2] {\rm Otherwise: Then the resolvent} $a\Phi$ {\rm is defined to be the sequent} $b:\Psi$.
\edes
\edf

\bdf\label{df:5forun}
{\rm Let} $J$ {\rm and} $J'$ {\rm be rules in a preproof such that both} $J$ {\rm and} $J'$ {\rm are one of rules} $(\Sig_i)\, (1\leq i\leq N-1)$ {\rm and} $J$ {\rm is above the right uppersequent of} $J'$. {\rm We say that} $J$ foreruns $J'$ {\rm if any right branch} $\calt$ {\rm of} $J'$ is left to $J${\rm , i.e., there exists a merging rule} $K$ {\rm such that} $\calt$ {\rm passes through the left side of} $K$ {\rm and the right uppersequent of} $K$ {\rm is equal to or below the right uppersequent of} $J$.
\[
\infer[(\Sig_i)^{\sig_0}\, J']{\Gam_0,\Lam_0}
{
 \infer*{\Gam_0,\lnot A^{\sig_0}}{}
&
 \infer*{A^{\sig_0},\Lam_0}
 {
  \infer[K]{\Gam_1,\Lam_1}
  {
   \infer*[\mbox{{\rm right branch }}\calt \mbox{{\rm of }} J']{\Gam_1,\lnot B}
   {A^{\sig_1}}
  &
   \infer*{B,\Lam_1}
   {
    \infer[(\Sig_j)\, J]{\Gam_2,\Lam_2}
    {
     \Gam_2,\lnot C
    &
     C,\Lam_2
     }
    }
   }
  }
 }
\]
\edf

If $J$ foreruns $J'$, then resolving steps of $J$ precede ones of $J'$. In other words we have to resolve $J$ in advance in order to resolve $J'$.


\bdf\label{df:5reach}
{\rm Let} $\calr=J_0,\ldots,J_{n-1}$ {\rm denote a series of rules} $(c)${\rm . Each} $J_p$ {\rm is a rule} $(c)^{\sig_p}_{\sig_{p+1}}${\rm Assume that} $J_0$ {\rm is above a rule} $(\Sig_i)^{\sig}\, I$ {\rm and} $\sig=\sig_p$ {\rm for some} $p$ {\rm with} $0<p\leq n${\rm . Then we say that the series} $\calr$ reaches to {\rm the rule} $I$.
\edf

In a {\it proof} defined in the next definition, if a series $\calr=J_0,\ldots,J_{n-1}$ reaches to the rule $(\Sig_i)^{\sig}\, I$, then either $\calr$ passes through $I$ in case $p<n$, or the subscript $\sig_{n}$ of the last rule $(c)^{\sig_{n-1}}_{\sig_{n}}\, J_{n-1}$ is equal to $\sig$, i.e., $J_{n-1}$ is a lowermost rule $(c)$ above $I$.

\bdf\label{df:5prf} Proof\\
{\rm Let} $P$ {\rm be a preproof. Assume a chain analysis for} $P$ {\rm is given. The preproof} $P$ together with the chain analysis {\rm is said to be a} proof {\rm in} $\mbox{{\rm T}}_{Nc}$ {\rm if it satisfies the following conditions
besides the conditions {\bf (pure)}, {\bf (h-reg)}, {\bf (c:side)}, {\bf (c:bound)}, {\bf (next)}, {\bf (h:bound)}, {\bf (ch:pass)}
(a chain passes through only rules $(c), (h),(\Sig_{i})\,(i<N)$), 
{\bf (ch:left)}, which are the same as in \cite{ptpi3}:}
\bdes

\item[(st:bound)] 
{\rm Let} $\calc$ {\rm be a 
chain, } $i\in In(\calc)$ {\rm and}  $a:\Gam$ {\rm be the uppersequent of the}
$i${\rm -origin of the chain} $\calc${\rm .}
 \bdes
 \item[(st:bound1)]
 {\rm Let} $i=N-1${\rm . Then}
 \[
 o(a;P)\leq st_{N-1}(\calc)
 .\]
 \item[(st:bound2)]
 {\rm Let} $i<N-1$ {\rm and} $\kap=rg_{i}(\calc)${\rm . Then for an} $\alp$
\[
st_{i}(\calc)=d_{\kap^{+}}\alp
\]
{\rm and}
\[
B_{\kap}(a;P)\leq \alp
.\]
 \edes

\item[(ch:link)] Linking chains: {\rm Let} $\calc=J_0,J_0',\ldots,J_n,J_n'$ {\rm and} $\cald=I_0,I_0',\ldots,I_m,I_m'$ {\rm be chains such that} $J_i$ {\rm is a rule} $(c)^{\tau_{i}}_{\tau_{i+1}}$ {\rm and} $I_i$ {\rm a rule} $(c)^{\sig_{i}}_{\sig_{i+1}}$. {\rm Assume that branchs of these chains intersect. Then one of the following three types must occur
(Cf. \cite{ptpi3} for {\bf Type1 (segment)} and {\bf Type2 (jump)}):}
 \bdes
 \item[Type1 (segment)]{\rm: One is a part of the other, i.e.,} 
 \[ n\leq m\spand J_i=I_i \]
{\rm or vice versa.}
 \edes
 
{\rm Assume that there exists a merging rule} $K$ {\rm such that} $\calc$ {\rm passes through the left side of} $K$ {\rm and} $\cald$ {\rm the right side of} $K$. {\rm Then by} {\bf (ch:left)} {\rm the merging rule} $K$ {\rm is a} $(\Sig_l)^{\tau_{j}}$ {\rm for some} $j\leq n$ {\rm and some} $l$ {\rm with} $1\leq l\leq N-2$. 

 \bdes 
 \item[Type2 (jump)]{\rm : The case when there is an} $i\leq m$ {\rm so that}
  \benu
  \item $J_{j-1}'$ {\rm is above} $K$ {\rm and} $J_j$ {\rm is below} $K$,
  \item $I_i$ {\rm is above} $K$,
  \item $I_i'$ {\rm is below} $J_n'$ {\rm and}
  \item $\sig_{i+1}<\tau_{n+1}$.
  \eenu
  
 \item[Type3 (merge)]{\rm : The case when} $\tau_{j}=\sig_{j}${\rm . Then it must be the case:}
  \benu
  \item $l>1$,
  \item $I_{j-1}'$ {\rm and} $J_{j-1}'$ {\rm are rules} $(\Sig_{N-1})^{\tau_{j}}$ {\rm above} $K${\rm , and}
  \item $n<m\spand J_{j+k}=I_{j+k}\spand J_{j+k}'=I_{j+k}'$ {\rm for any} $k$ {\rm with}\\
   $j\leq j+k\leq n$.
  \eenu
 {\rm That is to say,} $\calc$ {\rm and} $\cald$ {\rm share the part from} $J_{j}=I_{j}$ {\rm to} $J_{n}=I_{n}$ {\rm , the right chain} $\cald$ {\rm has to be longer} $n<m$ {\rm than the left chain} $\calc$ {\rm and the merging rule} $K$ {\rm is not a rule} $(\Sig_{1})$.
 \edes
{\rm If} {\bf Type2 (jump)} {\rm or} {\bf Type3 (merge)} {\rm occurs for chains} $\calc$ {\rm and} $\cald${\rm , then we say that} $\cald$ foreruns $\calc${\rm , since the resolving of the chain} $\cald$ {\rm precedes the resolving of the chain} $\calc$.
\[
\deduce
{\hskip0.3cm
 \infer[(c)^{\sig_{m}}_{\sig_{m+1}}\, I_m']{\Gam_m'}
 {
  \infer*{\Gam_m}
  {
   \infer[(\Sig_{N-1})^{\tau_{n+1}}\, J_n'=(\Sig_{N-1})^{\sig_{n+1}}\, I_n']{\Phi_n,\Psi_n}
   {
    \infer*{\Phi_n,\lnot A_n^{\tau_{n+1}}}{}
    &
    \infer*{A_n^{\tau_{n+1}},\Psi_n}
    {
     \infer[(c\Sig_{N-1})^{\tau_n}_{\tau_{n+1}}\, J_n=(c\Sig_{N-1})^{\sig_n}_{\sig_{n+1}}\, I_n]{\Gam_n'}
     {}
     }
    }
   }
  }
 }
 {
  \deduce
  {\hskip2.8cm
    \infer*{\Gam_n}
    {
     \infer[(c\Sig_{N-1})^{\tau_j}_{\tau_{j+1}}\, J_j=(c\Sig_{N-1})^{\sig_j}_{\sig_{j+1}}\, I_j]{\Gam_j'}
        {
         \infer*{\Gam_j}
         {}
        }
     }
    }
    {\hskip-3.0cm
      \infer[(\Sig_l)^{\tau_{j}}\, K]{\Phi,\Psi}
      {
       \infer*[\calc]{\Phi,\lnot A^{\tau_{j}}}
       {
        \infer[(\Sig_{N-1})^{\tau_{j}}\, J_{j-1}']{\Phi_{j-1},\Psi_{j-1}}
        {
         \infer*{\Phi_{j-1},\lnot A_{j-1}^{\tau_{j}}}{}
        &
         \infer*{ A_{j-1}^{\tau_{j}},\Psi_{j-1}}
          {}
         }
        }
       &
        \infer*[\cald]{ A^{\tau_{j}},\Psi}
       {
        \infer[(\Sig_{N-1})^{\sig_{j}}\, I_{j-1}']{\Pi,\Del}
        {
         \infer*{\Pi,\lnot B^{\sig_{j}}}{}
        &
         \infer*{B^{\sig_{j}},\Del}{}
        }
       }
      }
     }
   }
\hskip-2.0cm Type3
\]


\item[(ch:Qpt)] 
{\rm Let} $\calc=J_0,\ldots,J_n$ {\rm be a chain with a} $(c)^{\sig_p}_{\sig_{p+1}}\, J_p \, (p\leq n)$ 
{\rm and put} $\rho=\sig_{n+1}$. {\rm Then by} {\bf (ch:link)} 
{\rm there exists a uniquely detremined increasing sequence of numbers}
\beqn
\renewcommand{\theequation}{\ref{eqn:5rope}}
0\leq n_0<n_1<\cdots<n_l=n-1\, (l\geq 0)
\eeqn
\addtocounter{equation}{-1}
{\rm such that for each} $m<l$ {\rm there exists an} $i_m${\rm -knot} 
$(\Sig_{i_m})^{\sig_{n_{m}+1}}\, K_m \, (2\leq i_{m}\leq N-2)$ {\rm for the chain} $\calc$. 
{\rm (The} $i_m${\rm -knot} $K_m$ {\rm is the merging rule of the chain} $\calc$ 
{\rm and the chain starting with the rule} $J_{n_{m}+1}${\rm, cf.} {\bf Type3 (merge)}{\rm .) 
These numbers are called the} knotting numbers {\rm of the chain} $\calc$. 

{\rm Then} $pd_i(\rho), In(\rho), rg_i(\rho)$ {\rm have to be determined as follows:}
 \benu
 \item {\rm For} $2\leq i<N$, 
 \[pd_i(\rho)=\sig_{n_{m(i)}+1}\]
  {\rm with}
  \beqn\label{eqn:m(i)}
  m(i)=\max\{m:0\leq m\leq l\spand \fal p\in [0,m)(i\leq i_p)\}
  \eeqn
 {\rm that is to say,} 
 \[
 J_{n_{m(i)}}=pd_i(J_n)
 .\]
 \item {\rm For}  $2\leq i<N-1$
 \begin{eqnarray*}
 i\in In(\calc)=In(\rho) & \Lrarw & \exi p\in [0,m(i))(i_p=i) \\
 & \Lrarw & \exi p\in [0,l)(i_p=i\spand \fal q<p(i_q>i)) \\
 & \Lrarw & m(i)>m(i+1)=\min\{m<l: i_{m}=i\}
 \end{eqnarray*}
 {\rm And by the definition} $N-1\in In(\calc)=In(\rho)$.
 
 \item {\rm For} $i\in In(\calc) \spand i\neq N-1$,
  \benu
  \item {\rm The case when there exists a} $q$ {\rm such that} 
  \beqn\label{eqn:5rg}
  \exi p[n_{m(i)}\geq p\geq q>n_{m(i+1)}\spand \rho\prec_{i}\sig_{p+1}\spand \sig_{q}=rg_i(\sig_{p+1})]
  \eeqn
  {\rm Then} 
  \[rg_i(\rho)=\sig_q\]
  {\rm where} $q$ {\rm denotes the minimal} $q$ {\rm satisfying (\ref{eqn:5rg}).}
  \item {\rm Otherwise.}
  \[rg_i(\rho)=pd_i(\rho)=\sig_{n_{m(i)}+1}\]
  \eenu
 \eenu

\item[(lbranch)] 
{\rm Any left branch of a} $(\Sig_i)$ {\rm is the rightmost one in the left upper part of the} $(\Sig_i)$.
 
\item[(forerun)]
{\rm Let} $J^{lw}$ {\rm be a rule} $(\Sig_j)^{\sig}${\rm . Let} $\calr_{J_0}=J_0,\ldots,J_{n-1}$ {\rm denote the rope starting with a} $(c)\, J_0${\rm . Assume that} $J_0$ {\rm is above the right uppersequent of} $J^{lw}$ {\rm and the series} $\calr_{J_0}$ {\rm reaches to the rule} $J^{lw}${\rm . Then there is no merging rule} $K$, cf. {\rm the figure below, such that}
 \benu
 \item {\rm the chain} $\calc_0$ {\rm starting with} $J_0$ {\rm passes through the right side of} $K${\rm , and}
 \item {\rm a right branch} $\calt$ {\rm of} $J^{lw}$ {\rm passes through the left side of} $K$.
 \eenu
\[
\infer[(\Sig_j)^{\sig}\, J^{lw}]{\Phi,\Psi}
{
 \Phi,\lnot A
&
 \infer*[\calr_{J_0}]{A,\Psi}
 {
  \infer[(c)\, J_0]{\Gam_0'}
  {
   \infer*{\Gam_0}
   {
    \infer[K]{\Gam,\Lam}
    {
     \infer*[\calt]{\Gam,\lnot B}{}
    &
     \infer*[\calc_0]{B,\Lam}{}
    }
   }
  }
 }
}
\]

\item[(uplw)]
{\rm Let} $J^{lw}$ {\rm be a rule} $(\Sig_j)^{\sig}$ {\rm and} $J^{up}$ {\rm an} $i${\rm -knot} $(\Sig_i)^{\sig_0}\, (1\leq i,j\leq N)$. {\rm Let} $J_0$ {\rm denote the lower rule of} $J^{up}${\rm . Assume that the left rope} ${}_{J^{up}}\calr=J_0,\ldots,J_{n-1}$ {\rm of} $J^{up}$ {\rm reaches to the rule} $J^{lw}${\rm . Then}
\item[(uplwl)] {\rm if} $J^{up}$ {\rm is above the left uppersequent of} $J^{lw}${\rm , then} $j<i<N$.
 \[
 \infer[(\Sig_j)^{\sig}\, J^{lw}]{\Phi,\Psi}
 {
  \infer*[{}_{J^{up}}\calr]{\Phi,\lnot A}
  {
   \infer[(c)^{\sig_0}\, J_0]{\Gam_0'}
   {
    \infer*{\Gam_0}
    {
     \infer[(\Sig_i)^{\sig_0}\, J^{up}]{\Gam,\Lam}
     {
      \infer*[\calc_0]{\Gam,\lnot B}{}
     &
      \infer*{B,\Lam}{}
     }
    }
   }
  }
 &
  \infer*{A,\Psi}{}
 }
 \Longrightarrow j<i
 \]
 {\rm where} $\calc_0$ {\rm denotes the chain starting with} $J_0${\rm , and}
\item[(uplwr)] {\rm if} $J^{up}$ {\rm is above the right uppersequent of} $J^{lw}$ {\rm and} $i\leq j\leq N${\rm , then the rule} $(\Sig_i)^{\sig_0}\, J^{up}$ {\rm foreruns the rule} $(\Sig_j)^{\sig}\, J^{lw}$, cf. {\rm Proposition \ref{prp:N5.2} in Subsection \ref{subsec:5aQpart}}.
 
{\rm In other words if there exists a right branch} $\calt$ {\rm of} $J^{lw}$ {\rm as shown in the following figure, then} $j<i$.
\[
\deduce
{\hskip0.0cm
 \infer[(\Sig_j)^{\sig}\, J^{lw}]{\Phi,\Psi}
 {
  \infer*{\Phi,\lnot A}{}
 &
  \infer*[{}_{J^{up}}\calr]{A,\Psi}
  {
   \infer[{}^{\exi}K]{\Pi,\Del}
   {
    \infer*{\Pi,\lnot C}
    {
     \infer[(c)^{\sig_0}\, J_0]{\Gam_0'}
     {}
    }
   &
    \infer*[\calt]{C,\Del}{}
   }
  }
 }
}
{\hskip-1.5cm
 \infer*{\Gam_0}
    {
     \infer[(\Sig_i)^{\sig_0}\, J^{up}]{\Gam,\Lam}
     {
      \infer*[\calc_0]{\Gam,\lnot B}{}
     &
      \infer*[\calt]{B,\Lam}{}
     }
    }
 }
 \]   
\edes
\edf 

\noindent
{\bf Decipherment}. These provisos for a preproof to be a proof are obtained by inspection to rewrited proof figures. We decipher only additional provisos from \cite{ptpi3}.  
\bdes
\item[(ch:link)] Now a new type of linking chains, {\bf Type3 (merge)} enters, cf. Subsection \ref{subsec:5amerge}.

 For a chain $\cald=I_0,I_0',\ldots,I_m,I_m'$ and a member $I_n\, (n<m)$ of $\cald$ let $\calc=J_0,J_0',\ldots,J_n,J_n'$ denote the chain starting with $J_n=I_n$. Then there are two possibilities:
 \bdes
 \item[Type1 (segment)]$\calc$ is a part $I_0,\ldots,I_n$ of $\cald$ and hence the tops $I_0$ and $J_0$ are identical.
 \item[Type3 (merge)] The branch of $\calc$ is left to the branch of $\cald$.
 \edes
 
\item[(st:bound), (ch:Qpt)] By these provisos we see that an o.d. $\rho$ is in $Od(\Pi_N)$ for a newly introduced rule $(c)_{\rho}$, cf. Propositions \ref{prp:N5.1}, \ref{prp:N5.2} in Subsection \ref{subsec:5aQpart}, Lemma \ref{lem:5.3.240} below and the case {\bf M5.2} in the next Section \ref{sec:ml5}.

\item[(uplwl)] By the proviso we see that a preproof $P'$ which is resulted from a proof $P$
 is again a proof with respect to the proviso {\bf (ch:Qpt)}, cf. Lemma \ref{lem:5.3.22.A}.\ref{lem:5.3.22.A0}.

\item[(uplwr), (forerun), (lbranch)]
By these provisos we see that a preproof $P'$ which is resulted from a proof $P$ by resolving a rule $(\Sig_{i+1})$ is again a proof with respect to the provisos {\bf (forerun)} and {\bf (uplw)}, cf. Proposition \ref{prp:N5.2} in Subsection \ref{subsec:5aQpart}, the case {\bf M7.2} in the next Section \ref{sec:ml5}, Lemma \ref{lem:5.3.19} and Lemma \ref{lem:5.3.25}.
\edes

In the following any sequent and any rule are in a fixed proof.

As in the previous paper \cite{ptpi3} we have the following lemmata. 
Lemma \ref{lem:5.3.16} follows from the provisos {\bf (h-reg)} and {\bf (ch:link)} in Definition \ref{df:5prf}, 
, Lemma \ref{lem:5.3.18} from {\bf (h-reg)} and {\bf (c:bound1)}
and Lemma \ref{lem:(bar)} from {\bf (h-reg)}.

\blem\label{lem:5.3.16}
Let $J$ be a rule $(c)_\sig$ and $J'$ the trace $(\Sig_{N-1})^\sig$ of $J$. 
Let $J_1$ be a rule $(c)^\sig$ below $J'$.
If there exists a chain $\calc$ to which both $J$ and $J_1$ belong, 
then $J_1$ is the uppermost rule $(c)^{\sig}$ below $J$ and there is no rule $(c)$ between $J'$ and $J_1$.
\elem


\blem\label{lem:5.3.18}
Let $J_{top}$ be a rule $(c)^{\pi}$. Let $\Phi$ denote the bar of $J_{top}$. Assume that the branch $\calt$ from $J_{top}$ to $\Phi$ is the rightmost one in the upper part of $\Phi$. Then no chain passes through $\Phi$.
\elem

\blem\label{lem:(bar)}
Let $J$ be a rule $(c)$ and $b:\Phi$ the bar of the rule $J$.
Then there is no $(cut)\, I$ with $b\subset I\subset J$ nor a right uppersequent of a $(\Sig_{N})\, I$ 
with $b\subset I*(1)\incl J$ between $J$ and $b:\Phi$.
\elem

The following lemma is used to show that a preproof $P'$ which results from a proof $P$ by resolving a rule $(\Sig_{j})\, J^{lw}$ is again a proof with respect to the proviso {\bf (uplwl)}, cf. the Claim \ref{clm:5M8.2.3} in the case {\bf M7} in the next subsection.

\blem\label{lem:5.3.25}
Let $J^{lw}$ be a rule $(\Sig_j)$. Assume that there exists a right branch $\calt$ of $J^{lw}$ such that $\calt$ is the rightmost one in the upper part of $J^{lw}$. Then there is no $i$-knot $(\Sig_i)\, J^{up}$ above the right uppersequent of $J^{lw}$ such that $i\leq j$ and the left rope ${}_{J^{up}}\calr$ of $J^{up}$ and $J^{up}$ reaches to $J^{lw}$.
\elem
\bprf
Suppose such a rule $J^{up}$ exists. By {\bf (uplwr)} the rule $J^{up}$ foreruns $J^{lw}$. Thus the branch $\calt$ would not be the rightmost one.
\[
\infer[(\Sig_j)\, J^{lw}]{\Gam,\Lam}
{
 \Gam,\lnot A
&
 \infer*{A,\Lam}
 {
  \infer[(\Sig_i)\, J^{up}]{\Psi,\Phi}
  {
   \Psi,\lnot B
  &
   \infer*[\calt]{B,\Phi}{}
  }
 }
}
\]
\eprf

The following lemma is used to show that a preproof $P'$ which results from a proof $P$ by resolving a $(\Sig_{i+1})$ is a proof with respect to the proviso {\bf (uplwr)}, and to show a newly introduced rule $(\Sig_{i})$ in such a $P'$ does not split any chain, cf. the Claim \ref{clm:5M8.2.3} in the case {\bf M7}.

\blem\label{lem:5.3.19}
Let $J$ be a rule $(\Sig_{i+1})^{\sig_0}\, (0<i<N)$ and $b: \Phi$ the resolvent of $J$. 
Assume that the branch $\calt$ from $J$ to $b$ is the rightmost one in the upper part of $b$. 
Then every chain passing through $b$ passes through the right side of $J$.
\elem
\bprf
Let $a*(0)$ denote the lowermost rule $(\Sig_{i+1})^{\sig_{0}}$ below or equal to $J$,
 and $a: \Psi$ the lowersequent of $a*(0)$. 
The sequent $a: \Psi$ is the uppermost sequent below $J$ such that $h(a;P)<\sig_{0}+i$ by {\bf (h-reg)}. \\
{\bf Case 2}. $b=a$: 
If a chain passes through $a$ and a left side of a $(\Sig_{i+1})^{\sig_{0}}\, K_{-1}$
with $a\subset K_{-1}\subseteq J$,
then the chain would produce an $(i+1)$-knot $K_{-1}$.
\\
{\bf Case 1}. Otherwise:
Then there exists an $(i+1)$-knot $(\Sig_{i+1})^{\sig_{0}}$ with $a\subset K_{-1}\subseteq J$.
Let $(\Sig_{i+1})^{\sig}\, K_{-1}$ denote the uppermost such knot and 
${}_{K_{-1}}\calr=J_0,\ldots, J_{n-1}$ the left rope of $K_{-1}$. 
Each $J_p$ is a rule $(c)^{\sig_p}_{\sig_{p+1}}$. 
Let 
\beqn
\renewcommand{\theequation}{\ref{eqn:5rope}}
0\leq n_0<n_1<\cdots<n_l=n-1\, (l\geq 0)
\eeqn
\addtocounter{equation}{-1}
 be the knotting numbers of the left rope ${}_{K_{-1}}\calr$ and $K_m$ an $i_m$-knot 
 $(\Sig_{i_m})^{\sig_{n_{m}+1}}$ of $J_{n_m}$ and $J_{n_{m}+1}$ for $m<l$. 
 Put
\beqn
\renewcommand{\theequation}{\ref{eqn:m(i+1)}}
m(i+1)=\max\{m:0\leq m\leq l\spand \fal p\in [0,m)(i+1\leq i_p)\}
\eeqn
\addtocounter{equation}{-1}
 Then the resolvent $b: \Phi$ is the uppermost sequent $b:\Phi$ below $J_{n_{m(i+1)}}$ such that 
 \[
 h(b;P)<\sig_{n_{m(i+1)}+1}+i
 .\]
 Put 
 \[
m=m(+1),  \sig=\sig_{n_{m}+1}.
 \]
 Assume that there is a chain $\calc$ passing through $b$.
 As in {\bf Case 2} it suffices to show that the chain $\calc$ passes through the right side of $K_{-1}$.
Assume that this is not the case.
 Let $(c)^{\rho}_{\rho'}\, K$ denote the lowermost member of $\calc$ which is above $b$.
 \bclm\label{clm:5.3.19}
$K$ is on the branch $\calt$.
\eclm
{\bf Proof} of the Claim \ref{clm:5.3.19}. 
Assume that this is not the case.
Then we see that there exists a merging rule $(\Sig_{j})^{\rho'}\, I$ and a member $(c)^{\rho'}\, K'$ of $\calc$ such that the chain $\calc$ passes through the left side of $I$. 
$K'\subset b\subset I$ and hence $h(K';P)=\rho'\leq\sig$. 
We see $\rho'=\sig$ from {\bf (h-reg)}. 

Suppose $m=l$. 
Then by the definition of the left rope ${}_{K_{-1}}\calr$, the rule $(\Sig_{j})^{\rho'}\, I$ is not a knot, i.e.,
$j=1$.
But then $h(I*(1);P)=h(K';P)=\sig$, and hence $I\subset b$. A contradiction.
Therefore $m<l$ and $i_{m}\leq i$.
This means $K_{m}*(1)\incl b\subset I$.
On the other hand we have $1\leq i<j$ by $b\subset I$, and
by Lemma \ref{lem:5.3.16}
$K'$ is the uppermost rule $(c)^{\sig}$ below $(\Sig_{j})^{\sig}\, I$.
Therefore $(\Sig_{j})^{\rho'}\, I$ would be a knot below $J_{n_{m}}$.
On the other side $K_{m}$ is the uppermost knot below $J_{n_{m}}$.
This is a contradiction.

\[
\deduce
{\hskip4.1cm
\infer[(c)^{\sig}\, K'\in\calc]{\Del'}
{
 \infer*{\Del}
 {
    \infer[(\Sig_{i_{m}})^{\sig}\, K_{m}]{\Gam_1,\Lam_1}
     {
      \Gam_1,\lnot A_1
     &
      \infer*{A_1,\Lam_1}
       {
        \infer*{b: \Phi}
      {}
     }
    }
  }
 }
}
{\hskip3.0cm
   \infer[(\Sig_{j})^{\sig}\,I]{\Gam_0,\Lam_0}
   {
    \infer*[\calc]{\Gam_0,\lnot A_0}
    {
     \infer[(c)^{\rho}_{\sig}\, K\in\calc]{\Pi'}{\Pi}
     }
   &
    \infer*{A_0,\Lam_0}
    {
     J
     }
   }
 }
\]
\eprf
\\

\noindent
Then as in the proof of Lemma 7.13 of \cite{ptpi3} we see that
$K=J_{n_{m}}$, i.e., $(c)^{\rho}_{\rho'} \, K$ and $(c)^{\sig_{n_{m}}}_{\sig} \, J_{n_{m}}$
coincide.
Consider the chain $\calc_{m}$ starting with $(c)^{\sig_{n_{m}}}_{\sig} \, J_{n_{m}}$.
Then by {\bf (ch:link)} either $\calc_{m}$ is a segment of $\calc$ by {\bf Type1(segment)},
or $\calc$ foreruns $\calc_{m}$ by {\bf Type3(merge)}.
Since $(c)^{\sig_{n_{m}}}_{\sig} \, J_{n_{m}}$ is the lowest one such that
$\calc_{m}$ passes through the left side of $K_{m-1}$ and $J_{n_{m-1}+1}$
is a member of $\calc_{m}$,
{\bf Type1(segment)} does not occur.
In {\bf Type3(merge)} $K_{m-1}$ has to be the merging rule of $\calc_{m}$ and $\calc$
since, again, $(c)^{\sig_{n_{m}}}_{\sig} \, J_{n_{m}}$ is the lowest one,
and the branch $\calt$ is the rightmost one.
Therefore $\calc$ passes through the right side of $K_{m-1}$.
If $m=0$, then we are done.
Otherwise we see the chain $\calc$ and the chain $\calc_{m-1}$ starting with 
$(c)^{\sig_{n_{m-1}}}_{\sig} \, J_{n_{m-1}}$ has to share the rule 
$(c)^{\sig_{n_{m-1}}}_{\sig} \, J_{n_{m-1}}$.
As above we see that $\calc$ passes through the right side of $K_{m-2}$, and so forth.
\eprf
\\

\blem\label{lem:5.3.20}
Let $\calc= J_0,\ldots, J_n$ be a chain with rules $(c)^{\sig_p}_{\sig_{p+1}}\, J_{p}$ for $p\leq n$, and $(\Sig_{j})^{\sig_{p}}\, K\, (p<n)$ a rule such that $\calc$ passes through the right side of $K$ 
and the chain $\calc_{p}$ stating with $J_{p}$ passes through the left side of $K$. 
Further let $\calr={}_{K}\calr=J_{p},\ldots, J_{q-1}\, (q\leq n)$ denote the left rope of the $j$-knot $K$. 
Then the chain $\calc_{q}$ starting with $J_{q}$ is a part of the chain $\calc=\calc_{n}$, $\calc_{q}\incl\calc$. Therefore any knot for the chain $\calc$ is below $J_{q}$ and $q<n$, 
and in particular, if $K$ is a knot for the chain $\calc$, then $b=n$, cf. Definition \ref{df:5knot}.\ref{df:5knot.3}.
\elem
\bprf 
Suppose $q<n$.
By the Definiton \ref{df:5knot}.\ref{df:5knot.5} there is no knot of $J_{q-1}$ and $J_q$. 
Let $I_q$ denote a knot such that the chain $\calc_{q-1}$ starting with $J_{q-1}$ passes through 
the left side of $I_q$. 
$c\incl b$.
From the definition of a left rope we see that the chain $\calc_q$ starting with $J_q$ does not pass through the left side of the knot $I_q$.  Therefore by {\bf (ch:link)} {\bf Type1 (segment)} the chain $\calc_q$ must be a part of the chain $\calc$, $\calc_q\subset\calc$, i.e., the top of the chain $\calc_q$ is the top $J_0$ of $\calc$.
\[
\deduce
{\hskip1.0cm
 \infer[J_n]{\Gam_n'}
 {
  \infer*[\calc]{\Gam_n}
  {
   \infer[(c)^{\sig_{q}}\, J_{q}]{\Gam_{q}'}
   {}
  }
 }
}
{
\deduce
 {\hskip1.4cm
  \infer*[\calc_q]{\Gam_{q}}
  {
   \infer[(c)_{\sig_{q}}\, J_{q-1}]{\Gam_{q-1}'}
   {
    \infer*[\calr]{\Gam_{q-1}}
    {}
   }
  }
 }
 {
  \deduce
   {\hskip0.4cm
    \infer[I_q]{\Phi_{q-1},\Psi_{q-1}}
    {
     \infer*[\calc_{q-1}]{\Phi_{q-1},\lnot A_{q-1}}{}
    &
     \infer*[\calr]{A_{q-1},\Psi_{q-1}}
     {
      \infer[(c)^{\sig_{p}}\, J_{p}]{\Gam_{p}'}
      {}
     }
    }
   }
   {\hskip2.0cm
    \infer*{\Gam_{p}}
    {
     \infer[(\Sig_{j})^{\sig_{p}}\, K]{\Phi_{p},\Psi_{p}}
     {
      \infer*[\calc_{p}]{\Phi_{p},\lnot A_{p}}{}
     &
      \infer*[\calc_q\subset\calc]{A_{p},\Psi_{p}}{}
     }
    }
   }
  }
 }
\]
\eprf

The following Lemma \ref{lem:5.3.22.A}
is a preparation for Lemma \ref{lem:5.3.240}. 
From the Lemma \ref{lem:5.3.240} 
we see that an o.d. $\rho$ is in $Od(\Pi_N)$ for a newly introduced rule $(c)_{\rho}$, cf. the case {\bf M5.2} in the next subsection.

In the following Lemma \ref{lem:5.3.22.A}, 
$J$ denotes a rule $(c)_{\rho}$ and $\calc=J_0,\ldots,J_n$ the chain starting with $J_{n}=J$. 
Each $J_{p}$ is a rule $(c)^{\sig_p}_{\sig_{p+1}}$ for $p\leq n$ with $\sig_{n+1}=\rho$. 

$K$ denotes a rule $(\Sig_j)^{\sig_a}\, (j\leq N-2, 0<a\leq n)$ such that the chain $\calc$ passes through $K$.
If $\calc$ passes through the left side of $K$, then $j\leq N-2$ holds by {\bf (ch:left)}.

$J_{a-1}$ denotes the lowermost member $(c)_{\sig_a}$ of $\calc$ above $K$,
 $K\subset J_{a-1}$.

Let
\beqn
\renewcommand{\theequation}{\ref{eqn:5rope}}
0\leq n_0<n_1<\cdots<n_l=n-1\, (l\geq 0)
\eeqn
\addtocounter{equation}{-1}
be the knotting numbers of the chain $\calc$, cf. {\bf (ch:Qpt)}, 
 and $K_m$ an $i_m$-knot $(\Sig_{i_m})^{\sig_{n_{m}+1}}$ of 
 $J_{n_{m}}$ and $J_{n_{m}+1}$ for $m<l$. 
 Let $m(i)$ denote the number
\beqn
\renewcommand{\theequation}{\ref{eqn:m(i)}}
m(i)=\max\{m:0\leq m\leq l\spand \fal p\in [0,m)(i\leq i_p)\}
\eeqn
\addtocounter{equation}{-1}

\blem\label{lem:5.3.22.A} (cf. Proposition \ref{prp:N5.1} in Subsection \ref{subsec:5aQpart}.)\\
 \benu
 \item\label{clm:5.3.22.A.3}
Let $m\leq m(i)$. Then
\[
i\leq i_{m-1}\spand \sig_{n_{m}+1}\prec_i\sig_{n_{m-1}+1}
.\]

 \item\label{lem:5.3.22.A0}
Assume that $\calc$ passes through the left side of the rule $K$, i.e., $K*(0)\subset J_{a-1}$.
Then $J_{a-1}$ is the upper left rule of $K$. 
Let $i\leq j$.
\benu
\item \label{lem:5.3.22.A1} $\rho\prec_{i}\sig_a$, \\
and hence
\item the $i$-predecessor of $J$ is equal to or below $J_{a-1}$, and\label{lem:5.3.22.A2}
\item if $K_p$ is an $i_p$-knot $(\Sig_{i_p})$ for the chain $\calc$ above $K$, then $j<i_p$.\label{lem:5.3.22.A3}
\eenu
\[
\deduce
{\hskip1.0cm
 \infer[(c)^{\sig_n}_{\rho}\, J_{n}=J]{\Gam_n'}
 {
  \infer*{\Gam_n}
  {
   \infer[(\Sig_j)^{\sig_a}\, K]{\Phi,\Psi}
   {
    \infer*[\calc]{\Phi,\lnot A}
    {
     \infer[(c)_{\sig_a}\, J_{a-1}]{\Gam_{a-1}'}
     {}
    }
   &
    \infer*{A,\Psi}{}
   }
  }
 }
}
{\hskip-3.3cm
 \infer*{\Gam_{a-1}}
 {
  \infer[(\Sig_{i_p})^{\sig_{n_{p}+1}}\, K_p]{\Phi_p,\Psi_p}
  {
   \infer*{\Phi_p,\lnot A_p}{}
  &
   \infer*[\calc]{A_p,\Psi_p}{}
  }
 }
}
\Longrightarrow \rho\prec_{j}\sig_a \spand j<i_p
\]

\item\label{lem:5.3.26-1}
Let $J_{b-1}$ be a member of $\calc$ such that $\rho\prec_i\sig_{b}$ for an $i$ with $2\leq i\leq N-2$. 
Let $\calc_{b-1}$ denote the chain starting with $J_{b-1}$. 
Assume that the chain $\calc_{b-1}$ intersects $\calc$ of {\bf Type3 (merge)} in {\bf (ch:link)} and 
$(\Sig_j)\, K$ is the merging rule of $\calc_{b-1}$ and $\calc$. 

Then $i\leq j$.
\[
\infer[(c)_{\sig_{n+1}}\, J_n]{\Gam_n'}
{
 \infer*[\calc]{\Gam_n}
 {
  \infer[(c)_{\sig_b}\, J_{b-1}]{\Gam_{b-1}'}
  {
   \infer*{\Gam_{b-1}}
   {
    \infer[(\Sig_j)\, K]{\Phi,\Psi}
    {
     \infer*[\calc_{b-1}]{\Phi,\lnot A}{}
    &
     \infer*[\calc]{A,\Psi}{}
    }
   }
  }
 }
}
\spand \sig_{n+1}\prec_i\sig_{b} \Longrightarrow i\leq j
\]

\item\label{lem:5.3.22} 
Assume that $\calc$ passes through the left side of the rule $K$, i.e., $K*(0)\subset J_{a-1}$.
Let $i\leq j$.

Assume that the $i$-origin $J_{q}$ of $\calc$ is not below $K$, i.e., 
$\sig_q=rg_i(\rho)\darw\Rarw q<a$. 
Then
\[
\fal b\in(a,n+1]\{\rho\preceq_{i}\sig_{b}\prec_{i}\sig_{a}\to i\not\in In(\sig_{b})\}
\]
and hence
\beqnarrs
&& \fal b\in(a,n+1]\{\rho\preceq_{i}\sig_{b}\prec_{i}\sig_{a}\to in_i(J)=in_{i}(J_{b-1})=in_i(J_{a-1})\: 
\mbox{, i.e., } \\
&& in_{i}(\rho)=in_{i}(\sig_{b})=in_{i}(\sig_{a})\}
\eeqnarrs
In particular by Lemma \ref{lem:5.3.22.A}.\ref{lem:5.3.22.A0} we have
\[
\rho\prec_{i}\sig_{a}\spand in_i(J)=in_i(J_{a-1})\: \mbox{, i.e., } \: in_{i}(\rho)=in_{i}(\sig_{a})
.\]

\item\label{lem:5.3.26}

Assume that $\calc$ passes through the left side of the rule $K$. 
Let $J_{b-1}$ be a member of $\calc$ such that $J_{b-1}$ is below $K$, i.e., $a<b$,
 and assume that $\sig_{n+1}\preceq_i\sig:=\sig_{b}$ for an $i\leq j$. 
 If $\sig_q=rg_i(\sig)\darw \Rarw q<a$, then 
\[
\fal d\in(a,b]\{\sig\preceq_{i}\sig_{d}\prec_{i}\sig_{a}\to i\not\in In(\sig_{d})\}
\]
and
\[
\sig\prec_i\sig_a
.\]
Hence
\[
\fal d\in(a,b]\{\sig\preceq_{i}\sig_{d}\prec_{i}\sig_{a}\to in_{i}(\sig_{d})=in_{i}(\sig_{a})\}
\spand in_i(\sig)=in_i(\sig_{a})
.\]
The following figure depicts the case $\sig_{q}=rg_{i}(\sig)\darw$:
\[
\deduce
{\hskip0.0cm
 \infer[(c)_{\sig_{n+1}}\, J_n]{\Gam_n'}
 {
  \infer*[\calc]{\Gam_n}
  {
   \infer[(c)_{\sig}\, J_{b-1}]{\Gam_{b-1}'}
   {
    \infer*{\Gam_{b-1}}
   {}
   }
  }
 }
}
{\hskip-1.4cm
 \infer[(\Sig_j)^{\sig_a}\, K]{\Phi,\Psi}
 {
  \infer*[\calc]{\Phi,\lnot A}
  {
   \infer[(c)^{rg_i(\sig)}\, J_q]{\Gam_q'}
   {
    \infer*[\calc]{\Gam_q}{}
   }
  }
 &
  \infer*{A,\Psi}{}
 }
}
\]

\item\label{lem:5.3.23}

Assume that the chain $\calc$ passes through the left side of the rule $K$. 
For an $i\leq j$ assume that there exists a $q$ such that 
\[ 
\exi p[n\geq p\geq q\geq a\spand \rho\preceq_{i}\sig_{p+1}\spand \sig_{q}=rg_i(\sig_{p+1})]
.\]
Pick the minimal such $q_0$ and put $\kap=\sig_{q_0}$. 
Then
\benu
\item $\fal d\in(a,q_{0}]\{\sig_{q_{0}}\preceq_{i}\sig_{d}\prec_{i}\sig_{a}\to i\not\in In(\sig_{d})\}$ and
$in_i(J_{a-1})=in_i(J_{q_{0}-1})$, i.e., $in_i(\sig_a)=in_i(\kap)$ and $\kap\preceq_i\sig_a$.
\label{lem:5.3.23.1}
\item $\fal t[\rho\preceq_i\sig_t\prec_i\kap \Rarw rg_i(\sig_t)\preceq_i\kap]$.
\label{lem:5.3.23.2}
\eenu

\item\label{lem:5.3.27}
Assume that $\calc$ passes through the left side of the rule $K$. 
Let $J_{b-1}$ be a member of $\calc$ such that $J_{b-1}$ is below $K$, i.e., $a< b$ and 
$\rho\preceq_i\sig:=\sig_{b}$ for an $i\leq j$. 
Suppose $rg_i(\sig)\darw$ and put $\sig_q=rg_i(\sig)$. 
If the member $(c)^{\sig_{q}}\, J_{q}$ is below $K$, i.e., $a\leq q$, then
for $st_i(\sig)=d_{\sig_{q}^{+}}\alp$, cf. {\bf (st:bound)},
\[
B_{\sig_{q}}(c;P)\leq\alp
\]
 for the uppersequent $c:\Gam_q$ of the rule $J_q$.

\eenu
\elem
\bprf
First we show Lemmata \ref{lem:5.3.22.A}.\ref{clm:5.3.22.A.3} and \ref{lem:5.3.22.A}.\ref{lem:5.3.22.A0}
simultaneously by induction on the number of sequents between $K$ and $J$. 
\\

\noindent
{\bf Proof} of
Lemma \ref{lem:5.3.22.A}.\ref{clm:5.3.22.A.3}. 

By the definition of the number $m(i)$ we have $i\leq i_{m-1}$. 
Since the chain $\calc_{n_m}$ starting with $J_{n_{m}}$ passes through 
the left side of the $i_{m-1}$-knot $K_{m-1}$,
we have the assertion $\sig_{n_{m}+1}\prec_i\sig_{n_{m-1}+1}$ by IH 
on Lemma \ref{lem:5.3.22.A}.\ref{lem:5.3.22.A0}.
\[
\infer[(c)_{\sig_{n_{m}+1}}\, J_{n_{m}}]{\Gam_{n_{m}}'}
{
 \infer*{\Gam_{n_{m}}}
 {
  \infer[(\Sig_{i_{m-1}})^{\sig_{n_{m-1}+1}}\, K_{m-1}]{\Phi_{m-1},\Psi_{m-1}}
  {
   \infer*[\calc_{n_m}]{\Phi_{m-1},\lnot A_{m-1}}{}
  &
  A_{m-1},\Psi_{m-1}
  }
 }
}
\]
This shows Lemma \ref{lem:5.3.22.A}.\ref{clm:5.3.22.A.3}.
\eprf
\\

\noindent
{\bf Proof} of
Lemma \ref{lem:5.3.22.A}.\ref{lem:5.3.22.A0}.

By {\bf (ch:Qpt)} we have
 \[
 pd_i(\rho)=\sig_{n_{m(i)}+1}\: \mbox{{\rm and }}\: J_{n_{m(i)}}=pd_i(J)
 .\]
 
\bclm\label{clm:5.3.22.A.1}
$a\leq n_{m(i)}+1$, i.e., $J_{n_{m(i)}+1}\subset K$.
\eclm
{\bf Proof} of Claim \ref{clm:5.3.22.A.1}. 
If $m(i)=l$, then $a\leq n=n_{l}+1$. 
Assume 
\[
m(i)<l\neq 0 \spand a>n_{m(i)}+1
.\]
Then the $i_{m(i)}$-knot $(\Sig_{i_{m(i)}})^{\sig_{n_{m(i)}+1}}\, K_{m(i)}$ is 
above the left uppersequent of $K$, $K*(0)\subset K_{m(i)}$,
and $j\geq i>i_{m(i)}$. 
Consider the left rope ${}_{K_{m(i)}}\calr=J_{n_{m(i)}+1},\ldots, J_{b-1}$ of the knot 
$K_{m(i)}$ for the chain $\calc$. 
Then by Lemma \ref{lem:5.3.20} we have $b=n$. 
Therefore ${}_{K_{m(i)}}\calr$ reaches to the rule $K$.
Thus by {\bf (uplwl)} we have $i\leq j<i_{m(i)}$. 
This is a contradiction.
\eprf
\\

\noindent
By the Claim \ref{clm:5.3.22.A.1} we have Lemma \ref{lem:5.3.22.A}.\ref{lem:5.3.22.A2}.
\\
{\bf Case 1} $a=n_{m(i)}+1$: 
This means that the $i$-predecessor $J_{n_{m(i)}}$ of $J$ is the rule $J_{a-1}$, 
and $pd_i(\rho)=\sig_a$.
\\
{\bf Case 2} $a<n_{m(i)}+1$: This means that $J_{n_{m(i)}} \subset K$.
Put
\beqn\label{eqn:5.3.22A.m1}
m_1=\min\{m\leq m(i):a<n_m+1\}
\eeqn
Then $J_{n_{m_1}}$ is the uppermost rule $J_{n_{m}}$ below $K$. 
The chain $\calc_{n_{m_1}}$ starting with $J_{n_{m_1}}$ passes through the left side of the 
knot $(\Sig_{i_{m_{1}-1}})^{\sig_{n_{m_{1}-1}+1}}\, K_{m_{1}-1}$. 
If $K\subset K_{m_{1}-1}$,
then $\calc_{n_{m_1}}$ passes through the left side of $K$.
\[
\infer[(c)_{\sig_{n_{m_1}+1}}\, J_{n_{m_1}}]{\Gam_{n_{m_1}}'}
{
 \infer*{\Gam_{n_{m_1}}}
 {
  \infer[(\Sig_j)^{\sig_a}\, K]{\Phi,\Psi}
  {
   \infer*[\calc_{n_{m_1}}]{\Phi,\lnot A}{}
  &
   \infer*{A,\Psi}{}
  }
 }
}
\]
And by the minimality of $m_{1}$, if $K_{m_{1}-1}\subset K$,
then $J_{a-1}=J_{n_{m_{1}-1}}$,
i.e., $a=n_{m_{1}-1}+1$. 
\[
\deduce
{\hskip0.1cm
 \infer[(c)_{\sig_{n_{m_1}+1}}\, J_{n_{m_1}}]{\Gam_{n_{m_1}}'}
 {
  \infer*[\calc_{n_{m_1}}]{\Gam_{n_{m_1}}}
  {
   \infer[(\Sig_{i_{m_{1}-1}})^{\sig_{a}}\, K_{m_{1}-1}]{\Phi_{m_{1}-1},\Psi_{m_{1}-1}}
   {
    \infer*[\calc_{n_{m_1}}]{\Phi_{m_{1}-1},\lnot A_{m_{1}-1}}{}
   &
    \infer*[\calc]{A_{m_{1}-1},\Psi_{m_{1}-1}}{}
   }
  }
 }
}
{\hskip-0.0cm
 \infer[(\Sig_j)^{\sig_a}\, K]{\Phi,\Psi}
 {
  \infer*[\calc]{\Phi,\lnot A}
  {
   \infer[J_{a-1}=J_{n_{m_{1}-1}}]{\Gam_{a-1}'}
   {
    \infer*{\Gam_{a-1}}{}
   }
  }
 &
  \infer*{A,\Psi}{}
 }
}
\]

By Lemma \ref{lem:5.3.22.A}.\ref{clm:5.3.22.A.3} we have $pd_i(\rho)=\sig_{n_{m(i)}+1}\preceq_i\sig_{n_{m_1}+1}$. 
Once again by IH we have $\sig_{n_{m_1}+1}\preceq_i\sig_a$. 
Thus we have shown Lemma \ref{lem:5.3.22.A}.\ref{lem:5.3.22.A1}, $\rho\prec_{i}\sig_a$.
\eprf
\\

\noindent
{\bf Proof} of
Lemma \ref{lem:5.3.22.A}.\ref{lem:5.3.22.A3}. $j<i_p$: 
This is seen from {\bf (uplwl)} as in the proof of the Claim \ref{clm:5.3.22.A.1} since in this case we have $l\neq 0$.

A proof of Lemma \ref{lem:5.3.22.A}.\ref{lem:5.3.22.A0} is completed.
\eprf
\\

\noindent
{\bf Proof} of 
Lemma \ref{lem:5.3.22.A}.\ref{lem:5.3.26-1}.

The chain $\calc_{b-1}$ passes through the left side of $K$ and $\calc$ the right side of $K$. 
By {\bf (ch:Qpt)} we have
 \[
 pd_i(\sig_{n+1})=\sig_{n_{m(i)}+1}\: \mbox{{\rm and }}\:
 J_{n_{m(i)}}=pd_i(J_{n}).
 \]
If $K$ is a $K_m$ for an $m<l$, then the assertion $i\leq j=i_m$ follows from (\ref{eqn:m(i)}) 
since $b-1\leq n_{m(i)}$ by  $\sig_{n+1}\prec_{i}\sig_{b}$, and hence $m<m(i)$. 
 
Otherwise let $m\leq m(i)$ denote the number such that $n_m>b-1>n_{m-1}$, i.e., 
$J_{b-1}$ is between $K_{m-1}$ and $J_{n_{m}}$. 
Then $K$ is below $K_{m-1}$ and the rule $K$ is the merging rule of $\calc_{b-1}$ and 
the chain $\calc_{n_m}$ starting with $J_{n_{m}}$, i.e., 
$\calc_{n_m}$ passes through the right side of $K$.
\[
\deduce
{\hskip0.8cm
 \infer[(c)_{\sig_{n_{m}+1}}\, J_{n_m}]{\Gam_{n_m}'}
 {
  \infer*{\Gam_{n_m}}
  {
   \infer[(c)_{\sig_b}\, J_{b-1}]{\Gam_{b-1}'}
   {
    \infer*{\Gam_{b-1}}
    {}
   }
  }
 }
}
{\hskip-1.3cm
 \infer[(\Sig_j)\, K]{\Phi,\Psi}
 {
  \infer*[\calc_{b-1}]{\Phi,\lnot A}{}
 &
  \infer*[\calc_{n_m},\calc]{A,\Psi}
  {
   \infer[K_{m-1}]{\Phi_{m-1},\Psi_{m-1}}
   {
    \infer*[\calc_{n_m}]{\Phi_{m-1},\lnot A_{m-1}}{}
   &
    \infer*[\calc]{A_{m-1},\Psi_{m-1}}{}
   }
  }
 }
}
\]
 By IH it suffices to show that $\sig_{n_{m}+1}\prec_i\sig_{b}$ and this follows from 
\beqn\label{eqn:5.3.26-1}
\sig_{n+1}\prec_i\sig_{n_{m}+1}
\eeqn
 since the set $\{\tau:\sig_{n+1}\prec_i\tau\}$ is linearly ordered by $\prec_i$, 
 Proposition \ref{lem:5.4}.\ref{lem:5.3.1}.
Now (\ref{eqn:5.3.26-1}) follows from (\ref{eqn:m(i)}) and Lemma \ref{lem:5.3.22.A}.\ref{lem:5.3.22.A1}, i.e,
\[
\sig_{n+1}\prec_i pd_i(\sig_{n+1})=\sig_{n_{m(i)}+1}\prec_i\cdots\prec_i\sig_{n_{m-1}+1}\prec_i\sig_{n_{m}+1}
.\]
This shows Lemma \ref{lem:5.3.22.A}.\ref{lem:5.3.26-1}.
\eprf
\\

\noindent
{\bf Proof} of
Lemma \ref{lem:5.3.22.A}.\ref{lem:5.3.22} by induction on the number of sequents between $K$ and $J$.

By the Claim \ref{clm:5.3.22.A.1} we have $a\leq n_{m(i)}+1$.
\\
{\bf Case 1} $a=n_{m(i)}+1$: 
This means that the $i$-predecessor $J_{n_{m(i)}}$ of $J$ is the rule $J_{a-1}$ and $pd_i(\rho)=\sig_a$. 
By Lemma \ref{lem:5.3.22.A}.\ref{lem:5.3.22.A3} we have $i_p>j\geq i$ for any $p<m(i)$. 
On the other side by {\bf (ch:Qpt)} 
\beqn\label{eqn:5.3.22.In}
i\in In(\calc)=In(\rho)  \Lrarw  \exi p\in [0,m(i))(i_p=i)
\eeqn
Hence $i\not\in In(\rho)$. 
Thus $in_i(\sig_a)=in_i(pd_i(\rho))=in_i(\rho)$.
\\

\noindent
{\bf Case 2} $a<n_{m(i)}+1$: 
This means that $J_{n_{m(i)}}$ is below $K$. 
Let $m_1$ denote the number (\ref{eqn:5.3.22A.m1}) defined in the proof of 
Lemma \ref{lem:5.3.22.A}.\ref{lem:5.3.22.A0}.
\bclm\label{clm:5.3.22}
For each $m\in (m_1,m(i)]$ the $i$-origin of $J_{n_{m}}$ is not below $K_{m-1}$, 
$
i<i_{m-1}
$, 
$
i\not\in In(\rho)
$,
 $\sig_{n_{m}+1}\prec_{i}\sig_{n_{m-1}+1}$ 
 and 
 $\fal b\in(n_{m-1}+1,n_{m}+1]\{\sig_{n_{m}+1}\preceq_{i}\sig_{b}\prec_{i}\sig_{n_{m-1}+1}\to i\not\in In(\sig_{b})\}$
and
\beqnarr
&& \fal b\in(n_{m-1}+1,n_{m}+1]\{\sig_{n_{m}+1}\preceq_{i}\sig_{b}\prec_{i}\sig_{n_{m-1}+1}\to i\not\in In(\sig_{b})\}
\nonumber \\
&& \fal b\in(n_{m-1}+1,n_{m}+1]\{\sig_{n_{m}+1}\preceq_{i}\sig_{b}\prec_{i}\sig_{n_{m-1}+1} \to \nonumber \\
&& in_i(J_{n_{m}})=in_{i}(J_{b-1})=in_i(J_{n_{m-1}})\: \mbox{, i.e., } \nonumber \\
&& in_i(\sig_{n_{m}+1})=in_{i}(\sig_{b})=in_i(\sig_{n_{m-1}+1})\}
\label{eqn:5.3.22}
\eeqnarr
\eclm
{\bf Proof} of the Claim \ref{clm:5.3.22}. 
First we show $i<i_{m-1}$. 
By Lemma \ref{lem:5.3.22.A}.\ref{clm:5.3.22.A.3} we have $i\leq i_{m-1}$. 
Assume $i=i_{m-1}$ for some $m\in (m_1,m(i)]$. 
Pick the minimal such $m_2$. 
Then by {\bf (ch:Qpt)}, (\ref{eqn:5.3.22.In}) we have $i\in In(\rho)$ and hence $rg_i(\rho)\darw$. 
By Lemma \ref{lem:5.3.22.A}.\ref{lem:5.3.22.A3} we have 
\beqn\label{eqn:5.3.22.ab}
p<m_1 \Rarw i_p>j\geq i
\eeqn
Here $p<m_1$ means that $K_p$ is above $K$. 
Thus by {\bf (ch:Qpt)}
\[
m(i+1)=\max\{m:0\leq m\leq l\spand \fal p\in [0,m)(i+1\leq i_p)\}=m_{2}-1\geq m_{1}
,\]
i.e.,
\[
J_{n_{m_{2}-1}}=pd_{i+1}(J)\spand \sig_{n_{m_{2}-1}+1}=pd_{i+1}(\rho)\spand n_{m_{2}-1}\geq n_{m_{1}}\geq a
.\]
On the other hand by {\bf (ch:Qpt)} we have for the $i$-origin $J_{q}$ of $\calc$, i.e., $\sig_q=rg_i(\rho)$,
\[
n_{m_{2}-1}=n_{m(i+1)}<q\leq n_{m(i)}+1
.\]
Thus $J_{q}$ is below $J_{n_{m_{2}-1}}$ and hence by $a\leq n_{m_{2}-1}$ the $i$-origin $J_{q}$ is below $K$. 
This contradicts our hypothesis. 
\[
\deduce
{\hskip1.0cm
 \infer[(c)_{pd_{i}(\rho)}\, J_{n_{m(i)}}]{\Gam_{n_{m(i)}}'}
 {
  \infer*{\Gam_{n_{m(i)}}}
  {
   \infer[(c)^{rg_i(\rho)}\, J_q]{\Gam_q'}
   {
    \infer*{\Gam_q}
    {}
   }
  }
 }
}
{\hskip-2.5cm
 \infer[(\Sig_i)\, K_{m_{2}-1}]{\Phi_{m_{2}-1},\Psi_{m_{2}-1}}
 {
  \Phi_{m_{2}-1},\lnot A_{m_{2}-1}
 &
  \infer*{A_{m_{2}-1},\Psi_{m_{2}-1}}
  {
   \infer[(c)_{pd_{i+1}(\rho)}\, J_{n_{m_{2}-1}}]{\Gam_{n_{m_{2}-1}}'}
   {
    \infer*{\Gam_{n_{m_{2}-1}}}
    {
     \infer[K]{\Phi,\Psi}
     {
      \Phi,\lnot A
     &
      A,\Psi
     }
    }
   }
  }
 }
}
\]
 Thus we have shown $i<i_{m-1}$ for any $m\in (m_1,m(i)]$. 
 From this, (\ref{eqn:5.3.22.ab}) and (\ref{eqn:5.3.22.In}) we see $i\not\in In(\rho)$ and hence 
\[
in_i(\rho)=in_i(pd_i(\rho))=in_i(\sig_{n_{m(i)}+1})
.\]
In particular, if $rg_i(\rho)\darw$, then $\sig_q=rg_i(\rho)=rg_i(\sig_{n_{m(i)}+1})$, i.e.,
 the $i$-origin of $J_{n_{m(i)}}$ equals to the $i$-origin $J_{q}$ of $J$. 
Therefore the $i$-origin of $J_{n_{m(i)}}$ is not below $K$ and a fortiori not below the 
 $i_{m(i)-1}$-knot $K_{m(i)-1}$. 
 Also the chain starting with $J_{n_{m(i)}}$ passes through the left side of $K_{m(i)-1}$ 
 and $i\leq i_{m(i)-1}$. 
 Thus by IH we have (\ref{eqn:5.3.22}) for $m=m(i)$.
We see similarly that for each $m\in (m_1,m(i)]$ the $i$-origin of $J_{n_{m}}$ is not below 
$K_{m-1}$ and (\ref{eqn:5.3.22}). 

Thus we have shown the Claim \ref{clm:5.3.22}.
\eprf
\\

\noindent
From Claim \ref{clm:5.3.22} we see
\[
\fal b\in(n_{m_{1}}+1,n+1]\{\rho\preceq_{i}\sig_{b}\prec_{i}\sig_{n_{m_{1}}+1}\to i\not\in In(\sig_{b})\}
\]
and hence
\beqnarrs
&& \fal b\in(n_{m_{1}}+1,n+1]\{\rho\preceq_{i}\sig_{b}\prec_{i}\sig_{n_{m_{1}}+1}\to \\
&& in_i(J)=in_{i}(J_{b-1})=in_i(J_{n_{m_{1}}}) \mbox{, i.e., } \: in_{i}(\rho)=in_{i}(\sig_{b})=in_{i}(\sig_{n_{m_{1}}+1})\}.
\eeqnarrs
Further
\[
\fal m\in(m_{1},m(i)]\{\rho\prec_{i}\sig_{n_{m}+1}\prec_{i}\sig_{n_{m-1}+1}\preceq_{i}\sig_{n_{m_{1}}+1}\}
.\]
Once more by IH we  have, cf. Figures in the proof of Lemma \ref{lem:5.3.22.A}.\ref{lem:5.3.22.A0}, {\bf Case 2}., 
$\sig_{n_{m_{1}}+1}\preceq_{i}\sig_{a}$ and
\[\fal b\in(a,n_{m_{1}}+1]\{\sig_{n_{m_{1}}+1}\preceq_{i}\sig_{b}\prec_{i}\sig_{a}\to i\not\in In(\sig_{b})\}\]
and hence
\beqnarrs
&& \fal b\in(a,n_{m_{1}}+1]\{\sig_{n_{m_{1}}+1}\preceq_{i}\sig_{b}\prec_{i}\sig_{a}\to \\
&& in_i(J_{n_{m_{1}}})=in_{i}(J_{b-1})=in_i(J_{a-1})\mbox{, i.e., } \: in_{i}(\sig_{n_{m_{1}}+1})=in_{i}(\sig_{b})=in_{i}(\sig_{a})\}
\eeqnarrs
Thus we have shown Lemma \ref{lem:5.3.22.A}.\ref{lem:5.3.22}.
\eprf
\\

\noindent
{\bf Proof} of 
Lemma \ref{lem:5.3.22.A}.\ref{lem:5.3.26} by induction on the number of sequents between $K$ and $J_{b-1}$. 

Let $\calc_b=I_0,\ldots,I_{b-1}$ denote the chain starting with $J_{b-1}=I_{b-1}$. 
Each rule $I_p$ is again a rule $(c)^{\sig_p}_{\sig_{p+1}}$. 
Chains $\calc_{b}$ and  $\calc$ intersect in a way described as {\bf Type1 (segment)} or {\bf Type3 (merge)} 
in {\bf (ch:link)}. 
If the chain $\calc_b$ passes through the left side of $K$, then the $i$-origin $I_q$ of $\calc_b$ 
is above $K$ if it exists, and hence the assertion follows from Lemma  \ref{lem:5.3.22.A}.\ref{lem:5.3.22}.

 Otherwise there exists a merging rule $(\Sig_l)^{\sig_c}\, I$ below $K$ such that 
 the chain $\calc_b$ passes through the left side of $I$ and $\calc$ the right side of $I$.
\[
\deduce
{\hskip0.0cm
 \infer[(c)_{\sig}\, J_{b-1}]{\Gam_{b-1}'}
 {
  \infer*{\Gam_{b-1}}
  {
   \infer[(\Sig_l)^{\sig_c}\, I]{\Pi,\Del}
   {
    \infer*[\calc_b]{\Pi,\lnot B}{}
   &
    \infer*[\calc]{B,\Del}
    {
     \infer[(c)_{\sig_c}\, J_{c-1}]{\Gam_{c-1}'}
     {}
    }
   }
  }
 }
}
{\hskip-0.3cm
 \infer*{\Gam_{c-1}}
 {
  \infer[K]{\Phi,\Psi}
  {
   \infer*[\calc]{\Phi,\lnot A}{}
  &
   \infer*{A,\Psi}{}
  }
 }
}
\]
Then by Lemma \ref{lem:5.3.22.A}.\ref{lem:5.3.26-1} we have $i\leq l$. 
The $i$-origin $I_q$ of $\calc_b$ is not below $I$. 
Therefore by Lemma  \ref{lem:5.3.22.A}.\ref{lem:5.3.22} we have
\[
\fal d\in(c,b]\{\sig\preceq_{i}\sig_{d}\prec_{i}\sig_{c}\to i\not\in In(\sig_{d})\}
\]
and
\[
\sig\prec_i\sig_{c}
.\]
Hence
\[
\fal d\in(c,b]\{\sig\preceq_{i}\sig_{d}\prec_{i}\sig_{c}\to in_{i}(\sig_{d})=in_{i}(\sig_{c})\}
.\]
In particular
\beqn\label{eqn:5.3.26}
in_i(\sig_c)=in_i(\sig)\spand \sig\prec_i\sig_c
\eeqn 
Now consider the member $J_{c-1}$ of $\calc$. 
$J_{c-1}$ is again below $K$, $\sig_{n+1}\preceq_i\sig_{c}$ and 
$rg_i(\sig_c)\simeq rg_i(\sig)$ by (\ref{eqn:5.3.26}). 
Thus by IH we have 
\[
\fal d\in(a,c]\{\sig_{c}\preceq_{i}\sig_{d}\prec_{i}\sig_{a}\to i\not\in In(\sig_{d})\}
\]
and
\[
\sig_{c}\prec_i\sig_{a}
.\]
Therefore
\[
\fal d\in(a,c]\{\sig_{c}\preceq_{i}\sig_{d}\prec_{i}\sig_{a}\to in_{i}(\sig_{d})=in_{i}(\sig_{a})\}
.\]
This shows Lemma \ref{lem:5.3.22.A}.\ref{lem:5.3.26}.
\eprf
\\

\noindent
{\bf Proof} of
Lemma \ref{lem:5.3.22.A}.\ref{lem:5.3.23}.

Pick a $p_0$ so that $n\geq p_0\geq q_0$, $\rho\preceq_{i}\sig_{p_{0}+1}$ and 
$rg_i(\sig_{p_{0}+1})=\sig_{q_0}=\kap$.
\[
\deduce
{\hskip1.7cm
 \infer[(c)_{\rho}\, J_n]{\Gam_n'}
 {
  \infer*{\Gam_n}
  {
   \infer[(c)_{\sig_{p_{0}+1}}\, J_{p_0}]{\Gam_{p_0}'}
   {
    \infer*{\Gam_{p_0}}
    {
     \infer[(c)^{\kap}\, J_{q_0}\, (rg_i(\sig_{p_{0}+1})=\kap)]{\Gam_{q_0}'}
     {}
    }
   }
  }
 }
}
{\hskip-2.2cm
 \infer*{\Gam_{q_0}}
 {
  \infer[(\Sig_j)^{\sig_a}\, K]{\Phi,\Psi}
  {
   \infer*[\calc]{\Phi,\lnot A}
   {
    \infer[(c)_{\sig_a}\, J_{a-1}]{\Gam_{a-1}'}
    {
     \infer*{\Gam_{a-1}}{}
    }
   }
  &
   \infer*{A,\Psi}{}
  }
 }
}
\]
Lemma \ref{lem:5.3.22.A}.\ref{lem:5.3.23.1}. 
First note that $\rho\preceq_{i}\sig_{p_{0}+1}\prec_{i}rg_i(\sig_{p_{0}+1})=\kap$ 
by Proposition \ref{lem:5.4}.\ref{lem:5.4.1} (or by the proviso {\bf (ch:Qpt)}) and hence $\rho\prec_{i}\kap$.
Thus the assertion follows from Lemma \ref{lem:5.3.22.A}.\ref{lem:5.3.26} and the minimality of $q_0$.
\\

\noindent
Lemma \ref{lem:5.3.22.A}.\ref{lem:5.3.23.2}. 
Suppose $rg_i(\sig_t)\not\preceq_i\kap$ for a $t$ with $\rho\preceq_i\sig_t\prec_i\kap$. 
Put $\sig_b=rg_i(\sig_t)$. 
Then by Propositions \ref{lem:5.4}.\ref{lem:5.3.1} and \ref{lem:5.4}.\ref{lem:5.4.1} 
we have $\kap=\sig_{q_0}\prec_{i}\sig_b$ and $b<q_0<t\spand q_0\geq a$.
Hence by the minimality of $q_0$ we have $b<a$.
\[
\deduce
{\hskip0.0cm
 \infer[(c)_{\sig_{t}}\, J_{t-1}]{\Gam_{t-1}'}
 {
  \infer*{\Gam_{t-1}}
  {
   \infer[(c)^{\kap}\, J_{q_0}]{\Gam_{q_0}'}
   {
    \infer*{\Gam_{q_0}}
   {}
   }
  }
 }
}
{\hskip-1.6cm
 \infer[(\Sig_j)^{\sig_a}\, K]{\Phi,\Psi}
 {
  \infer*[\calc]{\Phi,\lnot A}
  {
   \infer[(c)^{\sig_{b}}\, J_{b}\, (\sig_{b}=rg_i(\sig_t))]{\Gam_{b}'}
   {
    \infer*{\Gam_{b}}{}
   }
  }
 &
  \infer*{A,\Psi}{}
 }
}
\]    
 Thus by Lemma \ref{lem:5.3.22.A}.\ref{lem:5.3.26} we have
\[
in_i(\sig_a)=in_i(\sig_t)
.\] 
From this and Lemma \ref{lem:5.3.22.A}.\ref{lem:5.3.23.1} we have 
\beqn\label{eqn:5.3.23}
in_i(\sig_a)=in_i(\sig_t)=in_i(\kap)\spand \sig_t\prec_i\kap\preceq_{i}\sig_a
\eeqn
{\bf Case 1} $t\leq p_0$: 
Then $\sig_{p_{0}+1}\prec_i\sig_t\prec_i\kap=rg_i(\sig_{p_{0}+1})$ by Proposition \ref{lem:5.4}.\ref{lem:5.3.1}. 
By Proposition \ref{lem:5.4}.\ref{lem:5.4.8} we would have $\sig_b=rg_i(\sig_t)\preceq_i\kap$. 
Thus this is not the case.

Alternatively we can handle this case without appealing Proposition \ref{lem:5.4}.\ref{lem:5.4.8} as follows. 
Let $p_{0}$ denote the minimal $p_0$ such that 
\[
n\geq p_0\geq q_0\spand \rho\preceq_{i}\sig_{p_{0}+1}\spand rg_i(\sig_{p_{0}+1})=\sig_{q_0}=\kap
.\]
 Then by {\bf (ch:Qpt)} we have $\kap=rg_i(\sig_{p_{0}+1})=pd_{i}(\sig_{p_{0}+1})$ 
 and hence this is not the case, i.e., $p_{0}<t$.
\\

\noindent
{\bf Case 2} $p_0<t$: 
Then $\sig_t\preceq_i\sig_{p_{0}+1}\prec_i\kap$. 
By (\ref{eqn:5.3.23}) and Proposition \ref{lem:5.4}.\ref{lem:5.4.6}, or by Lemma \ref{lem:5.3.22.A}.\ref{lem:5.3.26} 
we would have $in_i(\sig_{p_{0}+1})=in_i(\kap)$. 
In particular $\kap=rg_i(\sig_{p_{0}+1})=rg_i(\kap)$ but $rg_i(\kap)$ is a proper subdiagram of $\kap$. 
This is a contradiction.

This shows Lemma \ref{lem:5.3.22.A}.\ref{lem:5.3.23.2}.
\eprf
\\

\noindent
{\bf Proof} of Lemma \ref{lem:5.3.22.A}.\ref{lem:5.3.27} by induction on the number of sequents between $K$ and $J_{b-1}$. 
\\

\noindent
{\bf Case 1} 
$J_{q}$ is the $i$-origin of $J_{b-1}$, i.e., $J_q$ is a member of the chain starting with $(c)_{\sig}\, J_{b-1}$: 
By the proviso {\bf (st:bound)}
we can assume $i\not\in In(\sig)$. 
Then $in_i(\sig)=in_i(pd_i(\sig))=in_i(J_p)$ with $J_p=pd_i(J_{b-1})\spand a\leq b-1<p$ 
by Lemma \ref{lem:5.3.22.A}.\ref{lem:5.3.22.A0}.
In particular $rg_i(\sig)=rg_i(pd_i(\sig))$. 
IH and $st_i(pd_i(\sig))=st_i(\sig)$ yields the lemma.
\\

\noindent
{\bf Case 2} Otherwise: 
First note that $\sig_{n+1}\neq\sig$ and $\sig_{n+1}\prec_i\sig$. 
By {\bf (ch:Qpt)} we have
 \[
 pd_i(\sig_{n+1})=\sig_{n_{m(i)}+1}\: \mbox{{\rm and }}\:J_{n_{m(i)}}=pd_i(J_n)
 .\]
Also $\sig_{n_{m(i)}+1}\preceq_i\sig$ and hence $\sig_{n_{m(i)}+1}\leq \sig$. 
Let $m_1$ denote the number such that
\[
m_1=\min\{m:\sig_{n_{m}+1}\leq\sig\}\leq m(i)
.\]
Then the rule $(c)_{\sig}\, J_{b-1}$ is a member of the chain $\calc_{n_{m_1}}$
 starting with $J_{n_{m_1}}$ and $J_{b-1}$ is below $(\Sig_{i_{m_{1}-1}})\, K_{m_{1}-1}$. 
Also the chain $\calc_{n_{m_1}}$ passes thorugh the left side of the knot $K_{m_{1}-1}$.
By $m_1\leq m(i)$ and Lemma \ref{lem:5.3.22.A}.\ref{clm:5.3.22.A.3} we have
$\sig_{n_{m(i)}+1}\preceq_i\sig_{n_{m_1}}$ and hence 
\beqn\label{eqn:5.3.27}
i\leq i_{m_{1}-1}\spand  \sig_{n_{m_1}}\preceq_i\sig.
\eeqn
{\bf Case 2.1} 
$J_q$ is below $K_{m_{1}-1}$, i.e., $\sig_q\leq\sig_{n_{m_{1}-1}+1}$, i.e., $n_{m_{1}-1}<q$: 
By IH and (\ref{eqn:5.3.27}) we get the assertion. 
\\

\noindent
{\bf Case 2.2} Otherwise: 
By Lemma \ref{lem:5.3.22.A}.\ref{lem:5.3.26} and (\ref{eqn:5.3.27}) we have
\[
in_i(\sig)=in_i(\sig_{n_{m_1}+1})\spand \sig\prec_i\sig_{n_{m_1}+1}.
\]
Hence $st_i(\sig)=st_i(\sig_{n_{m_1}+1})\spand rg_i(\sig)=rg_i(\sig_{n_{m_1}+1})$. 
IH and $st_i(\sig_{n_{m_1}+1})=st_i(\sig)$ yield the lemma.
\eprf

\blem\label{lem:5.3.240}
Let $\calr=J_{0},\ldots, J_{n-1}$ denote the rope starting with a top $(c)^{\pi}\, J_0$. 
Each $J_{p}$ is a rule $(c)^{\sig_p}_{\sig_{p+1}}$. 
Let 
\beqn
\renewcommand{\theequation}{\ref{eqn:5rope}}
0\leq n_0<n_1<\cdots<n_l=n-1\, (l\geq 0)
\eeqn
\addtocounter{equation}{-1}
be the knotting numbers of the rope $\calr$, and $K_m$ an $i_m$-knot $(\Sig_{i_m})^{\sig_{n_{m}+1}}$ 
of $J_{n_{m}}:$ and $J_{n_{m}+1}$ for $m<l$. 
For $2\leq i<N$ let $m(i)$ denote the number
\beqn
\renewcommand{\theequation}{\ref{eqn:m(i)}}
m(i)=\max\{m:0\leq m\leq l\spand \fal p\in [0,m)(i\leq i_p)\}
\eeqn
\addtocounter{equation}{-1}
Note that $i_{m}\leq N-2$ by {\bf (ch:left)}.
Also put {\rm (cf. {\bf (ch:Qpt)})}
\benu
\item 
\[
pd_i=\sig_{n_{m(i)}+1}
.\]

\item 
\[
i\in In  \Lrarw  \exi p\in [0,m(i))(i_p=i)
.\]

\item For $i\in In$ $( i\neq N-1)$,
  \bdes
  \item[Case 1]  
  The case when there exists a $q$ such that 
  \beqn\label{eqn:5rg24}
  \exi p[n_{m(i)}\geq p\geq q>n_{m(i+1)}\spand pd_i\preceq_{i}\sig_{p+1}\spand \sig_{q}=rg_i(\sig_{p+1})]
  \eeqn
  Then 
  \[
  rg_i=\sig_q
  \]
   where $q$ denotes the minimal $q$ satisfying (\ref{eqn:5rg24}).
   
  \item[Case 2] 
  Otherwise.
  \[
  rg_i=pd_i=\sig_{n_{m(i)}+1}
  .\]
  \edes
\eenu

\benu
\item\label{lem:5.3.24}
For each $i\in In$ we have
\benu
\item 
$in_i(rg_i)=in_i(pd_{i+1})\spand pd_{i}\preceq_{i}rg_i\preceq_{i}pd_{i+1}\spand pd_i\neq pd_{i+1}$.
\label{lem:5.3.24.1}

\item 
$\fal t[rg_i(pd_i)\preceq_i\sig_t\prec_{i}rg_i\Rarw rg_i(\sig_t)\preceq_{i}rg_i]$.
\label{lem:5.3.24.2}

\item 
Either $rg_{i}=pd_{i}$ or $rg_{i}(pd_{i})\preceq_{i}rg_{i}$.
 \label{lem:5.3.24.3}
\eenu

\item\label{lem:5.3.28}
Assume $i\in In\spand \sig_q:=rg_i\neq pd_i$, i.e., {\bf Case 1} occurs. 
Then
\[
B_{\sig_{q}}(c;P)\leq\alp
\]
for the uppersequent $c:\Gam_q$ of the rule $J_q$, $st_i(\sig_{p+1})=d_{\sig_{q}^{+}}\alp$ and
$p$ denotes a number such that
\[
n_{m(i)}\geq p\geq q>n_{m(i+1)}\spand pd_i=\sig_{n_{m(i)}+1}\preceq_{i}\sig_{p+1}\spand \sig_{q}=rg_i(\sig_{p+1})
.\]
\eenu
\elem
\bprf 

Lemma \ref{lem:5.3.240}.\ref{lem:5.3.24}.

Let $i\in In$, and put $\sig_{q_0}=rg_i \spand \sig_{p_0}=pd_i\spand \sig_r=pd_{i+1}$. 
By the definition we have 
$p_0=n_{m(i)}+1\spand r=n_{m(i+1)}+1$, $m(i)>m(i+1)\spand i_{m(i+1)}=i$,
$p_0\leq q_0\leq r$ and $\sig_{p_0}\preceq_{i}\sig_{q_0}$. 
Also 
\[
\fal p\in [m(i+1),m(i))(i\leq i_p)
.\]
From this and Lemma \ref{lem:5.3.22.A}.\ref{clm:5.3.22.A.3} we see
\beqn\label{eqn:5.3.24.1}
\fal p\in [m(i+1),m(i))(\sig_{n_{p+1}+1}\prec_i\sig_{n_{p}+1})
\eeqn
On the other hand we have by the definition of $rg_{i}=\sig_{q_{0}}$
\beqn\label{eq:5.3.24.rg}
\lnot\exi q<q_{0}\exi p[p_{0}-1\leq p\leq q>r-1\spand pd_{i}\preceq_{i}\sig_{p+1}\spand \sig_{q}=rg_{i}(\sig_{p})]
\eeqn
{\bf Case 2}. 
Then $pd_i=rg_i$, i.e., $p_{0}=q_{0}$,
 and Lemma \ref{lem:5.3.240}.\ref{lem:5.3.24.2} vacuously holds.
Lemma \ref{lem:5.3.240}.\ref{lem:5.3.24.1}, $in_{i}(\sig_{q_{0}})=in_{i}(\sig_{r})\spand \sig_{q_{0}}\preceq_{i}\sig_{r}$,
 follows from (\ref{eqn:5.3.24.1}) and Lemma  \ref{lem:5.3.22.A}.\ref{lem:5.3.22} with (\ref{eq:5.3.24.rg}).
\\

\noindent
{\bf Case 1}. 
Let $m$ denote the number such that 
\beqn\label{eqn:5.3.24.Cs1}
m(i)\geq m>m(i+1)\spand n_m\geq q_0>n_{m-1}
\eeqn
i.e., the rule $(c)^{\sig_{q_0}}\, J_{q_0}$ is a member of the chain $\calc_{n_m}$ starting with $J_{n_m}$.
\bclm\label{clm:5.2.34.1}
Let $p_1$ denote the minimal $p_{1}$ such that 
$\sig_{p_{0}}\preceq_{i}\sig_{p_{1}+1}$ and $\sig_{q_0}=rg_i(\sig_{p_{1}+1})$. 
Then $p_1\leq n_m\spand \sig_{n_{m}+1}\preceq_i\sig_{p_{1}+1}$.
\[
\deduce
{\hskip0.0cm
 \infer[(c)_{\sig_{n_{m}+1}}\, J_{n_m}]{\Gam_{n_m}'}
 {
  \infer*{\Gam_{n_m}}
  {
   \infer[(c)_{\sig_{p_{1}+1}}\, J_{p_1}]{\Gam_{p_1}'}
   {
    \infer*{\Gam_{p_1}}
    {
     \infer[(c)^{rg_i(\sig_{p_{1}+1})}\, J_{q_0}]{\Gam_{q_0}'}
     {}
    }
   }
  }
 }
}
{\hskip-2.4cm
 \infer*{\Gam_{q_0}}
 {
  \infer[K_{m-1}]{\Phi_{m-1},\Psi_{m-1}}
  {
   \infer*[\calc_{n_m}]{\Phi_{m-1},\lnot A_{m-1}}{}
  &
   \infer*{A_{m-1},\Psi_{m-1}}{}
  }
 }
}
\]
\eclm
{\bf Proof} of Claim \ref{clm:5.2.34.1}. 
Let $m_1$ denote the number such that 
\[
m(i)\geq m_1>m(i+1)\spand n_{m_1}\geq p_1>n_{m_{1}-1}
.\]
 Then by (\ref{eqn:5.3.24.1}), $pd_i\preceq_i\sig_{n_{m_1}+1}$ and $pd_i\preceq_{i}\sig_{p_{1}+1}$ 
 we have $\sig_{n_{m_1}+1}\preceq_i\sig_{p_{1}+1}$. 
 It remains to show $m=m_1$. Assume $m<m_1$. 
 Then by Lemma \ref{lem:5.3.22.A}.\ref{lem:5.3.26} and $q_0<n_{m_{1}-1}+1$ we would have $in_i(\sig_{n_{m_{1}-1}+1})=in_i(\sig_{p_{1}+1})$ and hence 
 $rg_i(\sig_{n_{m_{1}-1}+1})=rg_i(\sig_{p_{1}+1})=\sig_{q_0}$. 
 This contradicts the minimality of $p_1$ by (\ref{eqn:5.3.24.1}).

\[
\deduce
{\hskip-2.3cm
 \infer[J_{n_{m_1}}]{\Gam_{n_{m_1}}'}
 {
  \infer*{\Gam_{n_{m_1}}}
  {
   \infer[(c)_{\sig_{p_{1}+1}}\, J_{p_1}]{\Gam_{p_1}'}
     {}
  }
 }
}
{
\deduce
 {\hskip-0.4cm
  \infer*{\Gam_{p_1}}
  {
   \infer[(\Sig_{i_{m_{1}-1}})^{\sig_{n_{m_{1}-1}}}\, K_{m_{1}-1}]{\Phi_{m_{1}-1},\Psi_{m_{1}-1}}
   {
    \infer*[\calc_{n_{m_{1}}}]{\Phi_{m_{1}-1},\lnot A_{m_{1}-1}}{}
   &
    \infer*{A_{m_{1}-1},\Psi_{m_{1}-1}}
    {}
   }
  }
 }
 {
 \deduce
  {\hskip1.2cm
   \infer[J_{n_m}]{\Gam_{n_m}'}
   {
    \infer*{\Gam_{n_m}}
    {
     \infer[(c)^{rg_i(\sig_{p_{1}+1})}\, J_{q_0}]{\Gam_{q_0}'}
     {}
    }
   }
  }
  {\hskip-1.2cm
   \infer*{\Gam_{q_0}}
   {
    \infer[K_{m-1}]{\Phi_{m-1},\Psi_{m-1}}
    {
     \infer*[\calc_{n_m}]{\Phi_{m-1},\lnot A_{m-1}}{}
    &
     \infer*{A_{m-1},\Psi_{m-1}}{}
    }
   }
  }
 }
}
\]
This shows Claim \ref{clm:5.2.34.1}.
\eprf
\\

\noindent
 By the minimality of $q_0$ and the Claim \ref{clm:5.2.34.1} $q_0$ is the minimal $q$ such that 
\[ 
\exi p[n_{m}\geq p\geq q\geq n_{m-1}+1\spand \sig_{n_{m}+1}\preceq_{i}\sig_{p+1}\spand \sig_{q}=rg_i(\sig_{p+1})]
.\]
Hence by Lemma \ref{lem:5.3.22.A}.\ref{lem:5.3.23} we have 
\beqn\label{eqn:5.3.24.2}
in_i(\sig_{n_{m-1}+1})=in_i(\sig_{q_0})\spand \sig_{q_0}\preceq_i\sig_{n_{m-1}+1}
\eeqn
and
\beqn\label{eqn:5.3.24.3}
\fal t[\sig_{n_{m}+1}\preceq_i\sig_t\prec_i \sig_{q_0} \Rarw rg_i(\sig_t)\preceq_i\sig_{q_0}]
\eeqn
Lemma \ref{lem:5.3.240}.\ref{lem:5.3.24.1}. 
By (\ref{eqn:5.3.24.2}) it suffices to show that 
\[
in_i(\sig_{n_{m-1}+1})=in_i(pd_{i+1})\spand \sig_{n_{m-1}+1}\preceq_{i}pd_{i+1}
.\]
This follows from (\ref{eqn:5.3.24.1}) and Lemma  \ref{lem:5.3.22.A}.\ref{lem:5.3.22} with (\ref{eq:5.3.24.rg}).
\\

\noindent
Lemma \ref{lem:5.3.240}.\ref{lem:5.3.24.2} and \ref{lem:5.3.240}.\ref{lem:5.3.24.3}.
In view of (\ref{eqn:5.3.24.3}) it suffices to show the
\bclm\label{clm:5.2.34.2}
$\sig_{p_{0}}\preceq_i\sig_t\prec_i\sig_{n_{m}+1} \Rarw rg_i(\sig_t)\preceq_i\sig_{q_0}$.
\eclm
{\bf Proof} of Claim \ref{clm:5.2.34.2} by induction on $t$ with $p_0>t>n_{m}+1$.  

Let $m_2\geq m$ denote the number such that $n_{m_{2}+1}\geq t>n_{m_{2}}$. 
Then the chain $\calc_{n_{m_{2}+1}}$ starting with  $J_{n_{m_{2}+1}}$ 
passes through the left side of the rule $(\Sig_{i_{m_2}})^{\sig_{n_{m_{2}}+1}}\, K_{m_2}$.
\[
\infer[(c)_{\sig_{n_{m_{2}+1}+1}}\, J_{n_{m_{2}+1}}]{\Gam_{n_{m_{2}+1}}'}
{
 \infer*{\Gam_{n_{m_{2}+1}}}
 {
  \infer[(c)_{\sig_t}\, J_{t-1}]{\Gam_{t-1}'}
  {
   \infer*{\Gam_{t-1}}
   {
    \infer[(\Sig_{i_{m_2}})^{\sig_{n_{m_{2}}+1}}\, K_{m_2}]{\Phi_{m_2},\Psi_{m_2}}
    {
     \infer*[\calc_{n_{m_{2}+1}}]{\Phi_{m_2},\lnot A_{m_2}}{}
    &
     \infer*{A_{m_2},\Psi_{m_2}}{}
    }
   }
  }
 }
}
\]
We have 
\beqn\label{eqn:5.3.24.4}
\sig_{n_{m_{2}+1}+1}\preceq_i\sig_t\prec_i\sig_{n_{m_2}+1}
\eeqn
 by (\ref{eqn:5.3.24.1}). 
 Put $\sig_b=rg_i(\sig_t)$.
 It suffices to show $b\geq q_{0}$.

First consider the case when $b\leq n_{m_2}$. 
Then by (\ref{eqn:5.3.24.4}) and Lemma \ref{lem:5.3.22.A}.\ref{lem:5.3.26} we have $in_i(\sig_t)=in_i(\sig_{n_{m_2}+1})$ and $rg_i(\sig_t)=rg_i(\sig_{n_{m_2}+1})$. 
Thus IH when $m_{2}>m$ and (\ref{eqn:5.3.24.3}) when $m_{2}=m$ yield $b\geq q_{0}$.

Next suppose $b>n_{m_2}$. 
Let $q_1\leq b$ denote the minimal $q\leq b$ such that
\[
\exi p[n_{m_{2}+1}\geq p\geq q\geq n_{m_{2}}+1 \spand \sig_{n_{m_{2}+1}}\preceq_{i}\sig_{p+1}\spand \sig_{q}=rg_i(\sig_{p+1})]
.\]
The pair $(p,q)=(t-1,b)$ enjoys this condition.

Then by Lemma \ref{lem:5.3.22.A}.\ref{lem:5.3.23} we have $\sig_{q_1}\preceq_i\sig_{n_{m_2}+1}$.
Thus $\sig_b\preceq_i\sig_{q_1}\preceq_{i}\sig_{n_{m_2}+1}\preceq_{i}\sig_{n_{m}+1}\preceq_i\sig_{q_0}$.
This shows Claim \ref{clm:5.2.34.2}.
\eprf
\\

\noindent
Lemma \ref{lem:5.3.240}.\ref{lem:5.3.28}.

First observe that as in (\ref{eqn:5.3.27}) in the proof of Lemma  \ref{lem:5.3.22.A}.\ref{lem:5.3.27},
\beqn\label{eqn:5.3.28}
\fal m\leq m(i)[\sig_{n_{m(i)}+1}\preceq_i\sig_{n_{m}+1}]
\eeqn
Put
\beqnarrs
m_1 & = & \min\{m:p\leq n_m\}
\\
m_2 & = & \min\{m:q\leq n_m\}.
\eeqnarrs
Then the rule $J_p$ [$J_q$] is a member of the chain $\calc_{n_{m_1}}$ [$\calc_{n_{m_2}}$] 
starting with $J_{n_{m_1}}$ [starting with $J_{n_{m_2}}$], resp. and $m(i+1)<m_2\leq m_1\leq m(i)$.
\bclm\label{clm:5.3.28}
(Cf. Claim \ref{clm:5.2.34.1}.)
There exists a $p_0$ such that
\[
in_i(\sig_{p+1})=in_i(\sig_{p_{0}+1})\spand \sig_{n_{m_2}+1}\preceq_i\sig_{p_{0}+1}\spand 
n_{m_2}\geq p_0>n_{m_{2}-1},
\]
i.e., the rule $J_{p_0}$ is a member of the chain $\calc_{n_{m_2}}$ starting with $J_{n_{m_2}}$.
\eclm
{\bf Proof} of Claim \ref{clm:5.3.28}. 
\benu
\item 
The case $m_1=m_2$: 
By (\ref{eqn:5.3.28}) and $\sig_{n_{m(i)}+1}\preceq_{i}\sig_{p+1}$ we have 
\beqn\label{eqn:5.3.28.2}
\sig_{n_{m_1}+1}\preceq_i\sig_{p+1}.
\eeqn
 Set $p_0=p$.
\item 
The case $m_2<m_1$: 
By (\ref{eqn:5.3.28.2}) and Lemma  \ref{lem:5.3.22.A}.\ref{lem:5.3.26} we have
\[
in_i(\sig_{p+1})=in_i(\sig_{n_{m_{1}-1}+1})=\cdots=in_i(\sig_{n_{m_2}+1})
\]
and
\[
\sig_{p+1}\prec_i\sig_{n_{m_{1}-1}+1}\prec_{i}\cdots\prec_i\sig_{n_{m_2}+1}
.\]
Set $p_0=n_{m_2}$.
\eenu
This shows the Claim \ref{clm:5.3.28}.
\eprf
\\

\noindent
By the Claim \ref{clm:5.3.28} and Lemma  \ref{lem:5.3.22.A}.\ref{lem:5.3.27} we conclude $rg_i(\sig_{p+1})=rg_i(\sig_{p_{0}+1})$ and 
\[
B_{\sig_{q}}(c;P)\leq b(st_i(\sig_{p_{0}+1}))=b(st_i(\sig_{p+1}))=\alp
.\]
\eprf

\bmlem\label{mlem:5}
If $P$ is a proof, then the endsequent of $P$ is true.
\emlem

In the next section we prove the Main Lemma \ref{mlem:5} by a transfinite induction on $o(P)\in Od(\Pi_N)\mid\Ome$. 

Assuming the Main Lemma  \ref{mlem:5} we see Theorem \ref{th:5} as in \cite{ptpi3}, i.e., attach $(h)^{\pi}$, 
$(c\Pi_2)^{\Ome}$ and $(h)^{\Ome}$ as last rules to a proof $P_0$ of $A^{\Ome}$ in $\mbox{T}_{N}$. 
\[
\infer[(h)^{\Ome}]{A^{\alp}}
{
 \infer[(c\Pi_{2})^{\Ome}_{\alp}]{A^{\alp}}
 {
  \infer[(h)^{\pi}]{A^{\Ome}}
  {
   \infer*[P_{0}]{A^{\Ome}}{}
  }
 }
}
\msfiv P
\]


\section{Proof of Main Lemma}\label{sec:ml5}
Throughout this section $P$ denotes a proof with a chain analysis in $\mbox{{\rm T}}_{Nc}$ 
and $r:\Gam_{rdx}$ the redex of $P$.
\\

\noindent
{\bf M1}. 
The case when $r:\Gam_{rdx}$ is a lowersequent of an explicit basic rule $J$.\\
{\bf M2}. 
The case when $r:\Gam_{rdx}$ is a lowersequent of an $(ind) \, J$.
\\
{\bf M3}. 
The case when the redex $r:\Gam_{rdx}$ is an axiom. 
\\
These are treated as in \cite{ptMahlo}, \cite{ptpi3}.
\\

\noindent
By virtue of {\bf M1-3} we can assume that the redex $r:\Gam_{rdx}$ of $P$ is a lowersequent of a rule $J=r*(0)$ such that $J$ is one of the rules $(\Pi_2^\Ome\mbox{{\rm -rfl}})$, $(\Pi_N\mbox{{\rm -rfl}})$ or an implicit basic rule.
\\

\noindent
{\bf M4}. 
$J$ is a $(\Pi_2^\Ome\mbox{{\rm -rfl}})$. 
As in \cite{ptMahlo} introduce a $(c)^{\Ome}_{d_{\Ome}\alp}$ and a $(\mbox{{\rm cut}})$.
\\

\noindent
{\bf M5}. 
$J$ is a $(\Pi_N\mbox{{\rm -rfl}})$.\\
{\bf M5.1}. 
There is no rule $(c)^{\pi}$ below $J$.
\[ 
\infer[(h)^{\pi}]{a_{0}:\Lam}
{
  \infer*{a:\Phi}
  {
   \infer[J]{r:\Gam}
   {
    \infer*{\Gam,A}{}
   &
    \infer*{\lnot\exi z(t<z\land A^z),\Gam}{}
   }
  }
 }
\msten  P\]
where $a:\Phi$ denotes the uppermost sequent below $J$ such that 
$h(a;P)=\pi$.
The sequent $a_{0}:\Lam$ is the lowersequent of the lowermost $(h)^{\pi}$.

Let $P'$ be the following:

\[\infer[(\Sig_{N-1})^{\sig}\, J_{0}']{a_{0}:\Lam}
 {
  \infer[(h)^{\pi}]{\Lam,A^{\sig}}
   {
    \infer[(c\Pi_{N})^{\pi}_{\sig}\, J_{0}]{\Phi,A^{\sig}}
    {
      \infer*{a_{1}:\Phi,A}
      {
       \infer[(w)]{\Gam,A}{\infer*{\Gam}{}}
      }
     }
    }
  &
   \infer[(h)^{\pi}]{\lnot A^{\sig},\Lam}
    {
     \infer[(w)]{\lnot  A^{\sig},\Phi}
      {
       \infer*{\lnot A^{\sig},\Phi}
        {
          \infer[(w)]{\lnot  A^{\sig},\Gam}
          {\infer*[z:=\sig]{\lnot  A^{\sig},\Gam }{}}
         }
       }
     }
  }
\msten P'\]
where the o.d. $\sig$ in the new $(c\Pi_{N})^{\pi}_{\sig}\, J_{0}$ is defined to be
\[
\sig=d_{\pi}^{q}\alp \mbox{ with } q=\nu\pi\pi N-1, \nu=o(a_{1};P') \mbox{ and }
\alp=\pi\cdot o(a_{1};P')+\calk_{\pi}(a;P)
\]
Namely $In(\sig)=\{N-1\}$, $st_{N-1}(\sig)=\nu$ and $pd_{N-1}(\sig)=rg_{N-1}(\sig)=\pi$.

Then as in \cite{ptpi3} we see that $\Phi\incl\Del^{\sig}$, $\alp<Bk_{\pi}(a;P) \spand \sig<o(a_{0};P)$,
$\sig\in Od(\Pi_{N})$ and $o(a_{0};P')<o(a_{0};P)$.
Hence $o(P')<o(P)$. Moreover in $P$, no chain passes through $a_{0}:\Lam$, 
and the new $(\Sig_{N})^{\sig}\,J_{0}'$ does not split any chain.
\\

\noindent
{\bf M5.2}. 
There exists a rule $(c)^{\pi}\, J_0$ below $J$.\\
Let $\calr=J_0,\ldots,J_{n-1}$ denote the rope starting with $J_{0}$. 
The rope $\calr$ need not to be a chain as contrasted with \cite{ptpi3}.
Each rule $J_p$ is a $(c)^{\sig_p}_{\sig_{p+1}}$. 
Put $\sig=\sig_{n}$. 


\[
\deduce
 {\hskip2.5cm
  \infer*{a:\Phi}
   {
    \infer[(\Sig_{N-1})^{\sig}\, J^{\prime}_{n-1}]{a_{n}:\Gamma_n}
    {
     \infer*{}
     {
      \infer[(c)^{\sig_{n-1}}_{\sig}\, J_{n-1}]{\Gamma_{n-1}'}
      {
       \infer*{a_{n-1}:\Gamma_{n-1}}
       {}
      }
     }
    }
   }
  }
{
  \deduce
   {\hskip1.9cm
    \infer[(c)^{\sig_{i}}_{\sig_{i+1}}\, J_i]{\Gamma_i'}
      {
       \infer*{a_{i}:\Gamma_i}
       {
        \infer[(c)^{\pi}_{\sig_1}\, J_0]{a_{0}\Gamma_0'}
        {
         \infer*{\Gamma_0}
         {}
         }
        }
       }
   }
   {\hskip0.1cm
    \infer[J]{r:\Gamma}
          {
          \infer*{\Gamma,A}{}
          &
          \infer*{\lnot\exists z(t<z\land A^z),\Gamma}{}
          }
   }
 }
\msten P\]
where $a_{n}:\Gam_n$ denotes the lowersequent of the trace $(\Sig_{N-1})^{\sig}\, J_{n-1}^{\prime}$ of $J_{n-1}$, and $a:\Phi$ the bar of the rule $(c)_{\sig}\, J_{n-1}$. 
Let $(\Sig_{N-1})^{\sig_{i+1}}\, J_{i}'$ denote the trace of $J_{i}$ for $0\leq i<n$. 
Put
\[
h:=\mbox{h}(a;P).
\]

By Lemma \ref{lem:5.3.18} there is no chain passing through the bar $a:\Phi$.

Let $P'$ be the following:
\[
\infer[(\Sig_{N})^{\rho}\, J_n']{a:\Phi}
{
\deduce
 {\hskip1.2cm
    \infer*{\Phi,A^{\rho}}
   {
    \infer[(c\Pi_{N})^{\sig}_{\rho}\, J_n]{\Gamma_n,A^{\rho}}
    {
     \infer{a_{n}^{l}:\Gam_n,A^{\sig}}
     {
      \infer*{}
      {
       \infer[J^{l}_{n-1}]{\Gamma_{n-1}',A^{\sig}}
       {}
      }
     }
    }
   }
  }
 {
  \deduce
   {\hskip-0.5cm
    \infer*{a_{n-1}^{l}:\Gamma_{n-1},A^{\sig_{n-1}}}
    {
     \infer[J^{l}_{i}]{\Gamma_i',A^{\sig_{i+1}}}
      {
       \infer*{a_{i}^{l}:\Gamma_i,A^{\sig_{i}}}
       {
        \infer[J_{0}^{l}]{\Gamma_0',A^{\sig_{1}}}
         {}
        }
       }
     }
   }
   {\hskip-0.4cm
    \infer*{a_{0}^{l}:\Gamma_0,A}
    {       
     \infer*{\Gamma,A}{}
    }
   }
  }
&
\deduce
 {\hskip1.3cm
  \infer*{\lnot A^{\rho},\Phi}
   {
    \infer[(w)\, J_{n}^{r}]{\lnot A^{\rho},\Gamma_n}
    {
     \infer{a_{n}^{r}:\lnot A^{\rho},\Gam_n}
     {
      \infer*{}
      {
       \infer[J_{n-1}^{r}]{\lnot A^{\rho},\Gamma_{n-1}'}
        {}
      }
     }
    }
   }
  }
 {
  \deduce
   {\hskip0.4cm
    \infer*{a_{n-1}^{r}:\lnot A^{\rho},\Gamma_{n-1}}
    {
     \infer[J^{r}_{i}]{\lnot A^{\rho},\Gamma_i'}
      {
       \infer*{a_{i}^{r}:\lnot A^{\rho},\Gamma_i}
       {
        \infer[J_{0}^{r}]{\lnot A^{\rho},\Gamma_0'}
         {}
        }
       }
      }
   }
   {\hskip0.4cm
    \infer*{a_{0}^{r}:\lnot A^{\rho},\Gamma_0}
    {
     \infer*{\lnot A^{\rho},\Gamma}{}
     }
    }
  }
}
\msfiv P'\]
For the proviso {\bf (lbranch)} in $P'$, 
any ancestor of the left cut formula of the new $(\Sig_{N})^{\rho}\, J_{n}'$ is a genuine $\Pi_{N}^{\tau}$-formula
$A^{\tau}$ for a $\tau$ with $\rho\preceq\tau$.
The formula $A^{\tau}$ is not in the branch $\calt$ from $r:\Gam$ to $a:\Phi$ in $P$
since no genuine $\Pi_{N}^{\tau}$-formula with $\tau>\Ome$
is on the rightmost branch $\calt$.
Therefore any left branch of the new $(\Sig_{N})^{\rho}\, J_{n}'$ 
is the rightmost one in the left upper part of the $J_{n}'$ in $P'$.

In $P'$, a new chain $J_{0}^{l},\ldots,J^{l}_{n-1},J_{n}$ starting with the new $J_n$ is in the chain analysis for $P'$ and $\rho=d_{\sig}^{q}\alp\in\cald_{\sig}$ is determined as follows:
\beqnarrs
&& b(\rho) = \alp=
\\
&  &
\max\{\calb_{\pi}(o(a_{n}^{l};P')),\calb_{>\sig}(\{\sig\}\cup(a_{n};P))\}+\ome^{o(a_{n}^{l};P')}+\max\{\calk_{\sig}(a_{n};P),\calk_{\sig}(h)\},
\eeqnarrs
$rg_{N-1}(\rho)=\pi$ and $st_{N-1}(\rho)=o(a_{0}^{l};P')$.

Let
\beqn
\renewcommand{\theequation}{\ref{eqn:5rope}}
0\leq n_0<n_1<\cdots<n_l=n-1\, (l\geq 0)
\eeqn
\addtocounter{equation}{-1}
be the knotting numbers of the rope $\calr$ and $K_m$ an $i_m$-knot $(\Sig_{i_m})^{\sig_{n_{m}+1}}$ of $J_{n_m}$ and $J_{n_{m}+1}$ for $m<l$. 
Let $m(i)$ denote the number
\beqn
\renewcommand{\theequation}{\ref{eqn:m(i)}}
 m(i)=\max\{m:0\leq m\leq l\spand \fal p\in [0,m)(i\leq i_p)\}.
\eeqn
\addtocounter{equation}{-1}
Then $pd_i(\rho), In(\rho), rg_i(\rho), st_i(\rho)$ are determined so as to enjoy the provisos {\bf (ch:Qpt)} and 
{\bf (st:bound)}.
 \benu
 \item 
 $pd_i(\rho)=\sig_{n_{m(i)}+1}$ for $2\leq i<N$. 
 Note that $pd_{i}(\rho)\neq\pi=\sig_{0}$ since $n_{0}\geq 0$, cf. the condition (\ref{eq:5pred}) in Section \ref{sec:WienpiN} which says that $pd_{N-1}(\rho)=\pi\Lrarw \sig=\pi$.
 \item 
 $N-1\in In(\rho)$ and $i\in In(\rho)  \Lrarw  \exi p\in [0,m(i))(i_p=i)$ for $2\leq i<N-1$.
 \item 
 Let $i\in In(\rho)\spand i\neq N-1$.
 $q$ denotes a number determined as follows.
  \bdes
  \item [Case 1] The case when there exists a $q$ such that
  \beqn
  \renewcommand{\theequation}{\ref{eqn:5rg}}
  \exi p[n_{m(i)}\geq p\geq q>n_{m(i+1)}\spand \rho\prec_{i}\sig_{p+1}\spand \sig_{q}=rg_i(\sig_{p+1})]
  \eeqn
  \addtocounter{equation}{-1}
  Then $q$ denotes the minimal $q$ satisfying (\ref{eqn:5rg}).
Note that $\rho\prec_{i}\sig_{p+1}$ is equivalent to $pd_i(\rho)=\sig_{n_{m(i)}+1}\preceq_{i}\sig_{p+1}$.
  \item [Case 2] Otherwise.
  Then set $q=n_{m(i)}+1$.
  \edes
 
 In each case set $rg_i(\rho)=\sig_q:=\kap$ for the number $q$, 
 and $st_i(\rho)=d_{\kap^{+}}\alp_{i}$ for
 \[
 \alp_{i}=B_{\kap}(a_{q}^{l};P')
 \]
where $a_{q}^{l}$ denotes the uppersequent $\Gam_q,A^{\sig_q}$ of $J^{l}_{q}$ in the left upper part of 
$(\Sig_N)^{\rho}\, J_{n}^{\prime}$ in $P'$. 

By Lemma \ref{lem:OdbndB} we have $\calb_{>\kap^{+}}(\alp_{i})\subset\calb_{>\kap}(\alp_{i})<\alp_{i}$, and hence 
$st_{i}(\rho)\in Od(\Pi_{N})$.
 \eenu

Obviously the provisos {\bf (ch:Qpt)} and {\bf (st:bound)} are enjoyed for the new chain 
$J_{0}^{l},\ldots,J^{l}_{n-1},J_{n}$.

Observe that, cf. (\ref{eq:Qrest}) in Section \ref{sec:WienpiN},
\[
\pi<\bet\in q=Q(\rho)\Rarw \bet=st_{N-1}(\rho)
.\]
\bclm\label{clm:5M.6.2}
$\rho=d_{\sig}^{q}\alp\in Od(\Pi_{N})$.
\eclm
{\bf Proof} of Claim \ref{clm:5M.6.2}.
\\
(\ref{eq:Odmu}) $\calb_{>\sig}(\{\sig,\alp\}\cup q)<\alp$:
By Lemma \ref{lem:OdbndB} we have $\calb_{>\sig}(\{\sig,\alp\})<\alp$.
It suffices to see $\calb_{>\sig}(q)<\alp$.
By the definition we have $\{pd_{i}(\rho),rg_{i}(\rho):i\in In(\rho)\}\subset \{\sig_{p},\sig_{p}^{+}:p\leq n\}$.
On the other hand we have $\calb_{>\sig}( \{\sig_{p},\sig_{p}^{+}:p\leq n\})\subset\calb_{>\sig}(\sig)$.

We have $\calb_{>\sig}(st_{N-1}(\rho))\subset\calb_{>\sig}(\alp)$.
Finally for $st_{i}(\rho)=d_{\kap^{+}}\alp_{i}$ with $i<N-1$, 
we have $\calb_{>\sig}(st_{i}(\rho))\subset\calb_{>\sig}(\{\sig,\alp_{i}\})\cup\{\alp_{i}\}$, and
$\calb_{>\sig}(\{\sig,\alp_{i}\})\subset\calb_{>\sig}(\{\sig,\alp\})$ and $\alp_{i}<\alp$.
\\

\noindent
$(\cald^{Q}.12)$:
\\
{\bf Case 2} This corresponds to $(\cald^{Q}.12.1)$, $\kap=rg_i(\rho)=pd_i(\rho)$. 
Let $\alp_{1}$ denote the diagram such that $\rho\preceq\alp_{1}\in\cald_{\kap}$.
Then
\[
\alp_1=\sig_{n_{m(i)}+2}\, ( pd_i(\rho)=\sig_{n_{m(i)}+1}\spand \sig_{n+1}=\rho)
.\]
 We have by Lemma \ref{lem:OdbndB} and {\bf (c:bound2)},
\[
\calb_{>\sig}(B_{\kap}(a_{n_{m(i)}+1};P))<B_{\kap}(a_{n_{m(i)}+1};P)\leq b(\alp_{1})
.\]
On the other hand we have for $st_{i}(\rho)=d_{\kap^{+}}\alp_{i}$
\[
\calb_{>\kap}(st_i(\rho))\subset\calb_{>\kap}(\{\kap,\alp_{i}\})\leq \calb_{>\sig}(B_{\kap}(a_{n_{m(i)}+1};P))
.\]
Thus $\calb_{>\kap}(st_i(\rho))<b(\alp_{1})$.
\\
{\bf Case 1} This corresponds to $(\cald^{Q}.12.2)$, $rg_i(\rho)=rg_i(pd_i(\rho))$ or to
$(\cald^{Q}.12.3)$, $ rg_i(pd_i(\rho))\prec_i\kap$ by Lemma \ref{lem:5.3.240}.\ref{lem:5.3.24.3}.
Let $p$ denote the maximal $p$ such that 
\[
rg_i(\sig_{p+1})=\sig_q=rg_i(\rho)\spand pd_i(\rho)\preceq_i\sig_{p+1}
.\]

 

$st_i(\rho)<st_i(pd_i(\rho))$ for the case $(\cald^{Q}.12.2)$ and 
$st_i(\rho)<st_i(\sig_{p+1})$ for the case $(\cald^{Q}.12.3)$ follow from 
Lemma \ref{lem:5.3.240}.\ref{lem:5.3.28} since for $rg_{i}(\sig_{p+1})=\sig_{q}=rg_{i}(\rho)$
\[
b(st_i(\rho))=B_{\sig_{q}}(a_{q}^{l};P')<B_{\sig_{q}}(a_{q};P)\leq b(st_i(\sig_{p+1}))
\]
and hence by Lemmata \ref{lem:rgbnd} and \ref{lem:OdbndB}
\[
st_i(\rho)<st_i(\sig_{p+1})
.\]
$(\cald^{Q}.11)$ and $(\cald^{Q}.12.3)$:
These follow from Lemma \ref{lem:5.3.240}.\ref{lem:5.3.24}.
\\
$(\cald^{Q}.2)$: 
$ \fal\tau\leq rg_i(\rho)(K_{\tau}st_i(\rho)<\rho)$.
For $\tau\leq\kap=rg_{i}(\rho)$ and $st_{i}(\rho)=d_{\kap^{+}}\alp_{i}$,
we have $K_{\tau}(st_{i}(\rho))=K_{\tau}\{\kap,\alp_{i}\}\leq\calk_{\tau}(a_{q}^{l};P')<\rho$ as in
the case {\bf M6.2} in \cite{ptpi3}.
\eprf
\\

\noindent
As in the case {\bf M6.2} in \cite{ptpi3} we see that $o(P')<o(P)$.

We have to verify that $P'$ is a proof. 
The provisos other than {\bf (uplwl)} are seen to be satisfied as in the case {\bf M5.2} of \cite{ptpi3}.
For the proviso {\bf (forerun)} see Claim \ref{clm:5M8.2.2} in the subcase {\bf M7.2} below.
It suffices to see that $P'$ enjoys the proviso {\bf (uplwl)} when the lower rule $J^{lw}$ is the new 
$(\Sig_N)^{\rho}\, J_{n}'$. For example the left rope ${}_{K_{m}}\calr$ of the 
$i_m$-knot $(\Sig_{i_m})^{\sig_{n_{m}+1}}\, K_{m}$ of $J_{n_m}$ and $J_{n_{m}+1}$ ends with 
the rule $(c)_{\sig}\, J_{n-1}$. 
We show the following claim.

\bclm\label{clm:5M6.2a}
Any left rope ${}_{J^{up}}\calr$ of a knot $J^{up}$ in the left upper part of the new 
$(\Sig_N)^{\rho}\, J_{n}'$ does not reach to $J_{n}'$. 
\eclm
{\bf Proof} of Claim \ref{clm:5M6.2a}.
Consider the original proof $P$. 
By Lemma \ref{lem:5.3.18} there is no chain passing through the bar $a:\Phi$ and hence 
it suffices to see that there is no rule $(c)^{\sig}_{\rho}$ above $a:\Phi$. 
First observe that we have $\rho<\tau$ for any rules $(c)_{\tau}$ and $(\Sig_{i})^{\tau}$ which 
are between $(\Pi_{N}\mbox{{\rm -rfl}})\, J$ and $a:\Phi$. 
Thus there is no rule $(c)^{\sig}_{\rho}$ on the branch $\calt_{0}$ from $(c)^{\pi}\, J_{0}$ to $a:\Phi$. 
Consider another branch $\calt$ above $a:\Phi$ and suppose that there is a rule $(c)^{\sig}_{\rho}\, I$ on $\calt$. 
We can assume that the merging rule $K$ of $\calt$ and $\calt_{0}$ is below $J_{0}$ and hence 
the rule $K$ is a $(\Sig_{i})^{\tau}$. 
By the proviso {\bf (h-reg)}(cf. Definition 5.4.4 in \cite{ptpi3}.)
we have $\tau\leq\sig$, i.e., 
$K$ is between $(c)^{\sig_{n-1}}_{\sig}\, J_{n-1}$ and $a:\Phi$. 
Then we have seen $\rho<\tau$ and hence the trace $(\Sig_{N-1})^{\rho}\, I_{0}$ of $(c)^{\sig}_{\rho}\, I$ is 
below $K$ by the proviso {\bf (h-reg)}. 
Therefore the chain stating with the trace $I_{0}$ passes through the left side of $K$. This is impossible by the proviso {\bf (ch:left)}.
\[
\deduce
{\hskip0.2cm
  \infer*[\calt,\calt_{0}]{a:\Phi}
  {
   \infer[(\Sig_{N-1})^{\rho}\,I_{0}]{\Psi_{0},\Phi_{0}}
   {
    \infer*{\Psi_{0},\lnot B^{\rho}}{}
   &
    \infer*{B^{\rho},\Phi_{0}}{}
   }
  }
 }
{\hskip0.0cm
 \infer[(\Sig_{i})^{\tau}\, K]{\Psi_{1},\Phi_{1}}
 {
  \infer*[\calt]{\Psi_{1},\lnot C^{\tau}}
  {
   \infer[(c)^{\sig}_{\rho}\, I]{\Psi_{2}'}{\Psi_{2}}
   }
 &
  \infer*[\calt_{0}]{C^{\tau},\Phi_{1}}
  {
   \infer[(c)^{\sig_{n-1}}_{\sig}\, J_{n-1}]{\Gamma_{n-1}'}{\Gam_{n-1}}
  }
 }
}
\msten P\]
\\

\noindent
In what follows we assume that $r*(0)=J$ is a basic rule. 
Let $v*(0)=I$ denote the vanishing cut of $r*(0)=J$. 
$v*(0)=I$ is either a $(\Sig_{i})$ or a $(cut)$.
\\

\noindent
{\bf M6}. 
$I$ is a $(\Sig_N)^{\sig}$.
\[
\infer[(\Sig_N)^{\sig}\, I]{v:\Gam,\Lam}
{
 \infer*{\Gam,A^{\sig}}{}
&
 \infer*{\lnot A^{\sig},\Lam}
  {
   \infer[(b\exi)\, J]{\exi x<\sig\lnot A_{N-1}^{\sig}(x),\Lam_0}
   {
    \infer*{\alp<\sig,\Lam_0}{}
   &
    \infer*{\lnot A_{N-1}^{\sig}(\alp),\Lam_0}{}
   }
  }
 }
 \msten P\]
where $A\equiv \fal xA_{N-1}(x)$ is a $\Pi_N$ formula. 
\\
Assuming $\alp<\sig$ let $P'$ be the following:
\[
\infer[(\Sig_{N-1})^{\sig}\, I_{N-1}]{v:\Gam,\Lam}
{
 \infer{\Gam,\Lam,\lnot A_{N-1}^{\sig}(\alp)}
 {
  \infer*{\Gam,A^{\sig}}{}
 &
  \infer*{\lnot A^{\sig},\Lam,\lnot A_{N-1}^{\sig}(\alp)}
  {
   \infer[(w)]{\lnot A^{\sig},\Lam_0,\lnot A_{N-1}^{\sig}(\alp)}
   {
   \infer*{\lnot A_{N-1}^{\sig}(\alp),\Lam_0}{}
   }
  }
 }
&
 \infer*{\Gam,A_{N-1}^{\sig}(\alp)}{}
}
\msten P'\]
where, the preproof ending with $\Gam,A_{N-1}^{\sig}(\alp)$ is obtained from the left upper part of $I$ in $P$ by inversion.

As in the case {\bf M6} of \cite{ptpi3} we see that $o(v;P')< o(v;P)$.

For the proviso {\bf (lbranch)} in $P'$, cf. the case {\bf M5.2}.
We verify that $P'$ is a proof with respect to the proviso {\bf (uplw)}. 

The proviso {\bf (uplwl)} when the lower rule $J^{lw}$ is the new $(\Sig_{N-1})^{\sig}\, I_{N-1}$: 
Consider the original proof $P$. 
By Lemma \ref{lem:(bar)} no left rope in the right upper part of $(\Sig_{N})^{\sig}\, I$ reaches to $I$. 
Also by {\bf (uplwl)} with the lower rule $J^{lw}=I$ there is no left rope of an $i$-knot $J^{up}$ reaching to $I$.

The proviso {\bf (uplwr)} when the lower rule $J^{lw}$ is the new $(\Sig_{N-1})^{\sig}\, I_{N-1}$: As above there is no left rope of an $i$-knot $J^{up}$ reaching to $I$.

The proviso {\bf (uplwr)} when the upper rule $J^{up}$ is the $(\Sig_{N-1})^{\sig}\, I_{N-1}$: $(\Sig_{N-1})^{\sig}\, I_{N-1}$ is not an $(N-1)$-knot since there is no chain passing through $(\Sig_{N})^{\sig}\, I$ by {\bf (ch:pass)}.

For the proviso {\bf (forerun)} see Claim \ref{clm:5M8.2.2} in the subcase {\bf M7.2} below.
\\

\noindent
{\bf M7}. 
$I$ is a $(\Sig_{i+1})^{\sig}$ with $1\leq i<N-1$.
\\
Then $J$ is either an $(\exi)$ or a $(b\exi)$.
Let $u_{0}:\Psi$ denote the uppermost sequent below $I$ such that $h(u_{0};P)<\sig+i$. 
Also let $u:\Phi$ denote the resolvent of $I$, cf. Definition \ref{df:5res}.
\\
{\bf M7.1} 
$u_{0}=u$.
\[
\infer*{u:\Psi}
{
 \infer[(\Sig_{i+1})^{\sig}\, I]{\Gam,\Lam}
 {
  \infer*{\Gam,\lnot A_{i+1}^{\sig}}{}
 &
  \infer*{A_{i+1}^{\sig},\Lam}
  {
   \infer[ x]{A_{i+1}^{\tau},\Lam_0}
   {
    \infer*{\alp<\tau,\Lam_0}{}
    &
    \infer*{A_{i}^{\tau}(\alp),\Lam_0}{}
   }
  }
 }
}
\msten P\]
where $A_{i+1}\equiv \exi yA_{i}(y)$ is a $\Sig_{i+1}$ formula. 
Also if $x$ is an $(\exi)$, then $\tau=\pi$ and the left upper part of the true sequent $\alp<\tau,\Lam_0$ is absent. 
Anyway $\sig\preceq\tau$.

Assuming $\alp<\tau$ and then $\alp<\sig$ by {\bf (c:bound)}, let $P'$ be the following:
\[
\infer[(\Sig_{i})^{\sig}]{u:\Psi}
{
\infer*{\Psi,A_{i}^{\sig}(\alp)}
{
 \infer{A_{i}^{\sig}(\alp),\Gam,\Lam}
 {
  \infer*{\Gam,\lnot A_{i+1}^{\sig}}{}
 &
  \infer*{A_{i+1}^{\sig}, A_{i}^{\sig}(\alp),\Lam}
  {
   \infer[(w)]{A_{i+1}^{\tau},A_{i}^{\tau}(\alp),\Lam_0}
   {
    \infer*{A_{i}^{\tau}(\alp),\Lam_0}{}
   }
  }
 }
}
&
\infer*{\lnot A_{i}^{\sig}(\alp),\Psi}
{
 \infer[(w)]{\Gam,\Lam,\lnot A_{i}^{\sig}(\alp)}
  {
   \infer*{\Gam,\lnot A_{i}^{\sig}(\alp)}{}
  }
 }
}
\msten P'\]
It is easy to see that $o(u;P')< o(u;P)$. 
For the proviso {\bf (lbranch)} in $P'$, cf. the case {\bf M5.2}.
To see that $P'$ is a proof with respect to the provisos
{\bf (forerun)},
{\bf  (uplw)}, cf. the subcase {\bf M7.2} below.
\\

\noindent
{\bf M7.2} Otherwise.\\
Let $K$ denote the lowermost rule $(\Sig_{i+1})^{\sig}$ below or equal to $I$. 
Then $u_{0}:\Psi$ is the lowersequent of $K$ by {\bf (h-reg)}.
There exists an $(i+1)$-knot $(\Sig_{i+1})^{\sig}$ which is between an uppersequent of $I$ and $u_{0}:\Psi$. 
Pick the uppermost such knot $(\Sig_{j+1})^{\sig}\, K_{-1}$ and let 
${}_{K_{-1}}\calr=J_0,\ldots,J_{n-1}$ denote the left rope of $K_{-1}$. 
Each $J_p$ is a rule $(c)^{\sig_p}_{\sig_{p+1}}$ with $\sig=\sig_0$. Let 
\beqn
\renewcommand{\theequation}{\ref{eqn:5rope}}
0\leq n_0<n_1<\cdots<n_l=n-1\, (l\geq 0)
\eeqn
\addtocounter{equation}{-1}
 be the knotting numbers of the left rope ${}_{K_{-1}}\calr$ and 
 $K_m$ an $i_m$-knot $(\Sig_{i_m})^{\sig_{n_{m}+1}}$ of $J_{n_m}$ and $J_{n_{m}+1}$ for $m<l$. 
 Put
\beqn
\renewcommand{\theequation}{\ref{eqn:m(i+1)}}
m(i+1)=\max\{m:0\leq m\leq l\spand \fal p\in [0,m)(i+1\leq i_p)\}
\eeqn
\addtocounter{equation}{-1}
Then the resolvent $u:\Phi$ is the uppermost sequent $u:\Phi$ below $J_{n_{m(i+1)}}$ such that 
\[
h(u;P)<\sig_{n_{m(i+1)}+1}+i.
\]
In the following figure of $P$ the chain $\calc_{n_{m+1}}$ starting with $J_{n_{m+1}}$ 
passes through the left side of $K_m$.
\[
  \deduce
  {\hskip1.5cm
   \infer*[\calt_1]{u:\Phi}
   {
    \infer[(c)^{\sig_{n_{m(i+1)}}}_{\sig_{n_{m(i+1)}+1}}\, J_{n_{m(i+1)}}]{\Gam_{n_{m(i+1)}}'}
    {
     \infer*{\Gam_{n_{m(i+1)}}}
     {
      \infer[(c)^{\sig_{n_{m+1}}}_{\sig_{n_{m+1}+1}}\, J_{n_{m+1}}]{\Gam_{n_{m+1}}'}
     {}
     }
    }
   }
  }
  {
  \deduce
   {\hskip0.7cm
    \infer*{\Gam_{n_{m+1}}}
    {
     \infer[(c)^{\sig_{n_{m}+1}}_{\sig_{n_{m}+2}}\, J_{n_{m}+1}]{\Gam_{n_{m}+1}'}
     {
      \infer*{\Gam_{n_{m}+1}}
      {}
     }
    }
   }
   {
   \deduce
    {\hskip0.8cm
     \infer[(\Sig_{i_{m}})^{\sig_{n_{m}+1}}\, K_{m}]{\Pi_{m},\Del_{m}}
     {
      \infer*[\calc_{n_{m+1}}]{\Pi_{m},\lnot B_{m}}{}
     &
      \infer*{B_{m},\Del_{m}}
      {
       \infer[(c)^{\sig_{n_{m}}}_{\sig_{n_{m}+1}}\, J_{n_{m}}]{\Gam_{n_{m}}'}
       {
        \infer*{\Gam_{n_{m}}}
        {}
       }
      }
     }
    }
    {
    \deduce
     {\hskip1.8cm
      \infer[(c)^{\sig}_{\sig_1}\,  J_0]{\Gam_0'}
      {
       \infer*[\calt_1]{\Gam_0}
       {}
       }
      }
      {\hskip0.6cm
        \infer[(\Sig_{i+1}^{\sig})\, I]{v:\Gam,\Lam}
        {
         \infer*[\calt_l]{\Gam,\lnot A_{i+1}^{\sig}}{}
        &
         \infer*{A_{i+1}^{\sig},\Lam}
         {
          \infer[J]{A_{i+1}^{\tau},\Lam_0}
          {
           \alp<\tau,\Lam_0
          &
           A_{i}^{\tau}(\alp),\Lam_0
          }
         }        
       }
      }
     }
    }
   }
\msfiv P\]
Assuming $\alp<\tau$ and then $\alp<\sig_{n}\leq\sig_{n_{m(i)}+1}$, let $P'$ be the following:
\[
 \infer[ I_i]{u:\Phi}
{
  \deduce
  {\hskip0.0cm
   \infer*{\Phi,A_{i}^{\sig_{n_{m(i+1)}+1}}(\alp)}
   {
    \infer[ J^{l}_{n_{m(i+1)}}]{\Gam_{n_{m(i+1)}}',A_{i}^{\sig_{n_{m(i+1)}+1}}(\alp)}
    {
     \infer*{\Gam_{n_{m(i+1)}},A_{i}^{\sig_{n_{m(i+1)}}}(\alp)}
     {
      \infer[ J^{l}_{n_{m+1}}]{\Gam_{n_{m+1}}',A_{i}^{\sig_{n_{m+1}+1}}(\alp)}
     {}
     }
    }
   }
  }
  {
  \deduce
   {\hskip-0.5cm
    \infer*{\Gam_{n_{m+1}},A_{i}^{\sig_{n_{m+1}}}(\alp)}
    {
     \infer[ J^{l}_{n_{m}+1}]{\Gam_{n_{m}+1}',A_{i}^{\sig_{n_{m}+2}}(\alp)}
     {
      \infer*{\Gam_{n_{m}+1},A_{i}^{\sig_{n_{m}+1}}(\alp)}
      {}
     }
    }
   }
   {
   \deduce
    {\hskip0.0cm
     \infer[ K^{l}_{m}]{\Pi_{m},\Del_{m},A_{i}^{\sig_{n_{m}+1}}(\alp)}
     {
      \infer*[\calc^{l}_{n_{m+1}}]{\Pi_{m},\lnot B_{m}}{}
     &
      \infer*{B_{m},\Del_{m},A_{i}^{\sig_{n_{m}+1}}(\alp)}
      {
       \infer[J^{l}_{n_{m}}]{\Gam_{n_{m}}',A_{i}^{\sig_{n_{m}+1}}(\alp)}
       {
        \infer*{\Gam_{n_{m}},A_{i}^{\sig_{n_{m}}}(\alp)}
        {}
       }
      }
     }
    }
    {
    \deduce
     {\hskip2.0cm
      \infer[ J^{l}_0]{\Gam_0',A_{i}^{\sig_{1}}(\alp)}
      {
       \infer*{\Gam_0,A_{i}^{\sig}(\alp)}
       {}
       }
      }
      {\hskip1.5cm
        \infer[I^{l}]{\Gam,\Lam,A_{i}^{\sig}(\alp)}
        {
         \infer*{\Gam,\lnot A_{i+1}^{\sig}}{}
        &
         \infer*{A_{i+1}^{\sig},\Lam,A_{i}^{\sig}(\alp)}
         {
          \infer[(w)]{A_{i+1}^{\tau},\Lam_0,A_{i}^{\tau}(\alp)}
          {
           A_{i}^{\tau}(\alp),\Lam_0
           }
          }
         }        
       }
      }
     }
    }
  &
   \deduce
  {\hskip-0.3cm
   \infer*[\calt^{r}_1\subset\calt_r]{u*(1):\lnot A_{i}^{\sig_{n_{m(i+1)}+1}}(\alp),\Phi}
   {
    \infer{\lnot A_{i}^{\sig_{n_{m(i+1)}+1}}(\alp),\Gam_{n_{m(i+1)}}'}
    {
     \infer*{\lnot A_{i}^{\sig_{n_{m(i+1)}}}(\alp),\Gam_{n_{m(i+1)}}}
     {
      \infer{\lnot A_{i}^{\sig_{n_{m+1}+1}}(\alp),\Gam_{n_{m+1}}'}
      {
       \infer*{\lnot A_{i}^{\sig_{n_{m+1}}}(\alp),\Gam_{n_{m+1}}}
      {}
      }
     }
    }
   }
  }
  {
  \deduce
   {\hskip-0.3cm   
    \infer{\lnot A_{i}^{\sig_{n_{m}+2}}(\alp),\Gam_{n_{m}+1}'}
    {
     \infer*{\lnot A_{i}^{\sig_{n_{m}+1}}(\alp),\Gam_{n_{m}+1}}
     {
      \infer{\lnot A_{i}^{\sig_{n_{m}+1}}(\alp),\Pi_{m},\Del_{m}}
      {
       \infer*[\calc^{r}_{n_{m+1}}]{\Pi_{m},\lnot B_{m}}{}
      &
       \infer*{\lnot A_{i}^{\sig_{n_{m}+1}}(\alp),B_{m},\Del_{m}}
      {}
      }
     }
    }
   }
   {
   \deduce
    {\hskip2.0cm
     \infer{\lnot A_{i}^{\sig_{n_{m}+1}}(\alp),\Gam_{n_{m}}'}
     {
      \infer*{\lnot A_{i}^{\sig_{n_{m}}}(\alp),\Gam_{n_{m}}}
      {
       \infer{\lnot A_{i}^{\sig_{1}}(\alp),\Gam_0'}
       {
        \infer*[\calt^{r}_1\subset\calt_r]{\lnot A_{i}^{\sig}(\alp),\Gam_0}
       {}
       }
      }
     }
    }
     {\hskip1.9cm
        \infer[(w)]{v^{r}:\lnot A_{i}^{\sig}(\alp),\Gam,\Lam}
        {
         \infer*[\calt_r]{\Gam,\lnot A_{i}^{\sig}(\alp)}{}
        }
      }
   }
  }
}
\]
\hspace*{\fill}$P'$
\\
Here $I_i$ denotes a $(\Sig_i)^{\sig_{n_{m(i+1)}+1}}$.

It is straightforward to see $o(u;P')<o(u;P)$. 
We show $P'$ is a proof.

First by Lemma \ref{lem:5.3.19}, in $P$ every chain passing through the resolvent $u:\Phi$ 
passes through the right side of $I$ and, by inversion, 
the right upper part of $I$ disappears in $P'$. 
Hence the new $(\Sig_i)^{\sig_{n_{m(i+1)}+1}}\, I_i$ does not split any chain.
For the proviso {\bf (lbranch)} in $P'$, cf. the case {\bf M5.2}.

\bclm\label{clm:5M8.2.2}
The proviso {\bf (forerun)} holds for the lower rule $J^{lw}=I_i$ in $P'$.
\eclm
{\bf Proof} of Claim \ref{clm:5M8.2.2}. 
Consider a right branch $\calt_{r}$ of $I_i$. 
We show that there is no rule $K$ such that $\calt_r$ passes through the left side of $K$ and 
$h(a;P')<\pi$ with the lowersequent $a$ of $K$. 
The assertion follows from this and {\bf (h-reg)}. 
The ancestors of the right cut formula $\lnot A_{i}^{\sig_{n_{m(i+1)}+1}}(\alp)$ of $I_i$ 
comes from the left cut formula $\lnot A_{i+1}^{\sig}$ of $I$ in $P$. 
Let $\calt^{r}_1$ denote the branch in $P'$
 from the lowersequent $v^{r}:\lnot A_{i}^{\sig}(\alp),\Gam,\Lam$ of the new $(w)$ 
to the right uppersequent $u*(1):\lnot A_{i}^{\sig_{n_{m(i+1)}+1}}(\alp),\Phi$ of $I_i$.
 Also let $\calt_l$ denote a (the) left branch of $I$ in $P$. 
 There exists a (possibly empty) branch $\calt_0$ such that $\calt_r=\calt_0^{\frown}\calt_l^{\frown}\calt^{r}_1$. 
 By {\bf (lbranch)} any left branch $\calt_l$ of $I$ is the rightmost one in the left upper part of $I$. 
 Therefore there is no rule $K$ such that $\calt_r$ passes through the left side of $K$ and $h(a;P')<\pi$ with the lowersequent $a$ of $K$.
\eprf

\bclm\label{clm:5M8.2.1}
The proviso {\bf (uplwr)} holds for the upper rule  $J^{up}=I_i$ in $P'$.
\eclm
{\bf Proof} of Claim \ref{clm:5M8.2.1}. 
Suppose that $I_i$ is a knot. 
Then there exists a chain $\calc_1$ starting with an $I_1$ such that $\calc_{1}$ passes through the left side of $I_i$. This chain comes from a chain in $P$ which passes through $u:\Phi$. 
Call the latter chain in $P$ $\calc_{1}$ again. 
Further assume that, in $P'$, 
 the left rope ${}_{I_{i}}\calr$ of $I_i$ reaches to a rule $(\Sig_j)^{\kap}\, J^{lw}$ with $i\leq j$. 
Let  $I_{2}$ denote the lower rule of $I_{i}$. We have to show $I_i$ foreruns $J^{lw}$. 
It suffices to show that, in $P$, any right branch $\calt$ of $J^{lw}$ passes through the right side of $I$ 
if the branch $\calt$ passes through $u:\Phi$. 
Since, by inversion, the right upper part of $I$ disappears in $P'$, 
for such a branch $\calt$ there exists a unique branch $\calt'$ corresponding to it in $P'$ 
so that $\calt'$ passes through the left side of $I_i$ and hence $\calt'$ is left to $I_i$. 
\[
\infer[(\Sig_j)^{\kap}\, J^{lw}]{\Gam_{lw},\Lam_{lw}}
{
 \infer*{\Gam_{lw},\lnot C_{lw}}{}
&
 \infer*[\calt]{C_{lw},\Lam_{lw}}
 {
  \infer*{u:\Phi}
  {
   \infer[(\Sig_{i+1})^{\sig}\, I]{\Gam,\Lam}
        {
         \infer*{\Gam,\lnot A_{i+1}^{\sig}}{}
        &
         \infer*[\calt]{A_{i+1}^{\sig},\Lam}{}
         }
    }
   }
  }
  \msten P\]
\smallskip
\[
\infer[ J^{lw}]{\Gam_{lw},\Lam_{lw}}
{
 \infer*{\Gam_{lw},\lnot C_{lw}}{}
&
 \infer*[\calt']{C_{lw},\Lam_{lw}}
 {
  \infer[ I_i]{u:\Phi}
  {
   \infer*{\Phi,A_{i}^{\sig_{n_{m(i)}+1}}(\alp)}
   {
    \infer[I^{l}]{\Gam,\Lam,A_{i}^{\sig}(\alp)}
        {
         \infer*{\Gam,\lnot A_{i+1}^{\sig}}{}
        &
         \infer*[\calt']{A_{i+1}^{\sig},\Lam,A_{i}^{\sig}(\alp)}{}
        }
    }
   &
    \infer*{\lnot A_{i}^{\sig_{n_{m(i)}+1}}(\alp),\Phi}{}
   }
  }
 }
\msten P'\]       
{\bf Case 1}. 
The case when, in $P$, there exists a member $I_3$ of the chain $\calc_1$ 
such that $I_3$ is between $u:\Phi$ and $J^{lw}$, 
and the chain $\calc_3$ starting with $I_3$ passes through the resolvent $u:\Phi$ in $P$: 
Then by Lemma \ref{lem:5.3.19} the chain $\calc_3$ passes through the right side of $I$. 
The rope $\calr_{I_3}$ starting with $I_3$ in $P$ corresponds to a part (a tail) of the left rope ${}_{I_{i}}\calr$ in $P'$. 
Thus by the assumption the rope $\calr_{I_3}$ also reaches to $J^{lw}$ in $P$.
Hence by {\bf (forerun)} there is no merging rule $K$ such that
 \benu
 \item the chain $\calc_3$ starting with $I_3$ passes through the right side of $K$, and
 \item the right branch $\calt$ of $J^{lw}$ passes through the left side of $K$.
 \eenu
Therefore the right branch $\calt$ of $J^{lw}$ passes through the right side of $I$ in $P$.
\[
\deduce
{\hskip-0.5cm
 \infer[ J^{lw}]{\Gam_{lw},\Lam_{lw}}
 {
  \infer*{\Gam_{lw},\lnot C_{lw}}{}
 &
  \infer*[\calr_{I_3}\subset{}_{I_{i}}\calr]{C_{lw},\Lam_{lw}}
  {
   \infer[I_3]{\Del_3'}
   {
    \infer*{\Del_3}
    {}
   }
  }
 }
}
{
\deduce
 {\hskip0.5cm
  \infer[I_2]{\Del_2'}
  {
   \infer*{\Del_2}
   {
    \infer[ I_i]{u:\Phi}
    {
     \infer*{\Phi,A_{i}^{\sig_{n_{m(i)}+1}}(\alp)}
     {}
    &
     \infer*{\lnot A_{i}^{\sig_{n_{m(i)}+1}}(\alp),\Phi}{}
    }
   }
  }
 }
 {\hskip-3.0cm
  \infer[I^{l}]{\Gam,\Lam,A_{i}^{\sig}(\alp)}
  {
   \infer*{\Gam,\lnot A_{i+1}^{\sig}}{}
  &
   \infer*[\calc_{3}, \, \calt']{A_{i+1}^{\sig},\Lam,A_{i}^{\sig}(\alp)}{}
   }
  }
 }
\msten P'\]     
{\bf Case 2}. Otherwise: First we show the following claim:
\bclm\label{clm:5M8.2.11}
 In $P$, we have $m(i+1)<l$ for the number of knots $l$ in (\ref{eqn:5rope}),
and $I_2$ is the lower rule of the $i_{m(i+1)}$-knot $K_{m(i+1)}$. 
Let ${}_{K_{m(i+1)}}\calr$ denote the left rope of $K_{m(i+1)}$ in $P$. 
Then ${}_{K_{m(i+1)}}\calr$ reaches to $J^{lw}$. 
\eclm
{\bf Proof} of Claim \ref{clm:5M8.2.11}. 
In $P'$, the lower rule $I_{2}$ of the knot $I_{i}$ is a member of the chain $\calc_{1}$ 
starting with $I_{1}$ and passing through the left side of $I_{i}$. 
Further $I_{2}$ is above $J^{lw}$ since the left rope ${}_{I_{i}}\calr$ of $I_{i}$ is assumed to reach to $J^{lw}$ in $P'$, 
cf. Definition \ref{df:5reach}. 
Since we are considering when {\bf Case 1} is not the case, in $P$, 
$I_{1}$ is below $J^{lw}$ and the chain $\calc_{2}$ starting with $I_{2}$ 
does not pass through $u:\Phi$, and hence 
chains $\calc_{1}$ and $\calc_{2}$ intersect as {\bf Type3 (merge)} in {\bf (ch:link)}. 
In other words there is a knot below $u:\Phi$ whose upper right rule is $(c)_{\sig_{n_{m(i+1)}+1}}\, J_{n_{m(i+1)}}$. This means that the knot is the $i_{m(i+1)}$-knot $K_{m(i+1)}$. 
Thus we have shown that $m(i+1)<l$ and $I_2$ is the lower rule of the $i_{m(i+1)}$-knot $K_{m(i+1)}$.

Next we show that the left rope ${}_{K_{m(i+1)}}\calr$ of $K_{m(i+1)}$ reaches to $J^{lw}$ in $P$. 
Suppose this is not the case. 
Let $(c)_{\kap_{4}}\, I_{4}$ denote the lowest (last) member of the left rope ${}_{K_{m(i+1)}}\calr$. 
Then $\kap<\kap_{4}$ for the rule $(\Sig_{j})^{\kap}\, J^{lw}$. 
By $\kap<\kap_{4}$, the next member $(c)^{\kap_{4}}\, I_{5}$ of the chain $\calc_{1}$ is above $J^{lw}$.
Since we are considering when {\bf Case 1} is not the case, 
the chain $\calc_{5}$ starting with $I_{5}$ does not pass through $u:\Phi$. 
By Definition \ref{df:5knot}.\ref{df:5knot.5} of left ropes and {\bf (ch:link)} there would be a knot $K'$ 
whose lower rule is $I_{5}$ and whose upper right rule is $I_{4}$. 
This is a contradiction since $I_{4}$ is assumed to be the last member of the left rope ${}_{K_{m(i+1)}}\calr$.
This shows Claim \ref{clm:5M8.2.11}.

In the following figure note that $u:\Phi$ is above $K_{m(i+1)}$ by {\bf (h-reg)} and the definition of the resolvent
 $u:\Phi$.

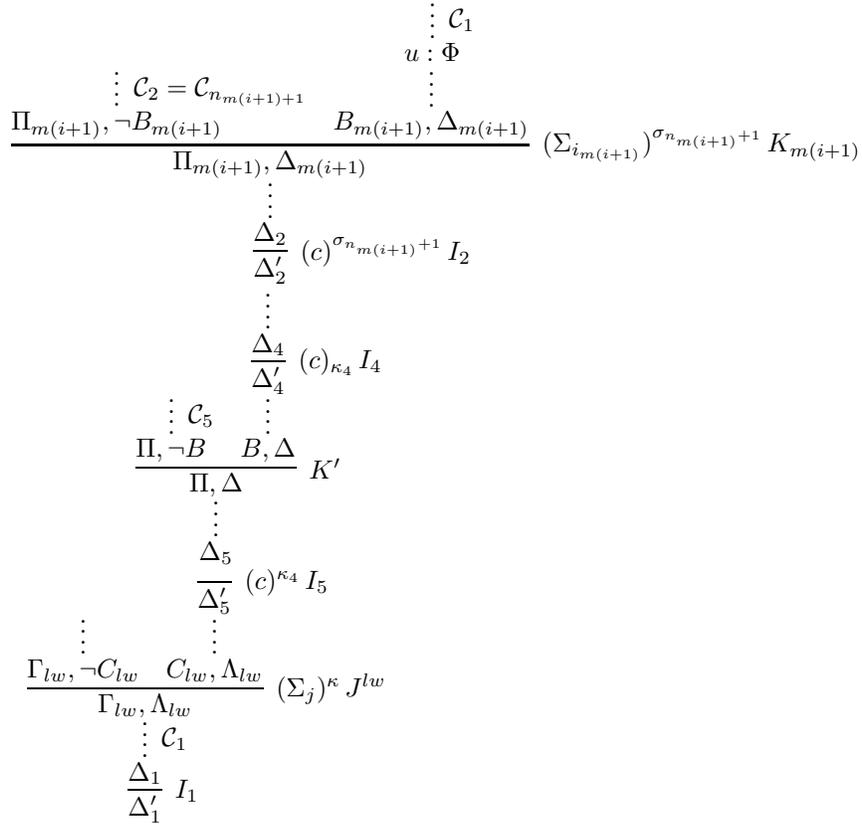
\[
\deduce
{\hskip-6.1cm
 \infer[I_{1}]{\Del_{1}'}
 {
  \infer*[\calc_{1}]{\Del_{1}}
  {
   \infer[(\Sig_j)^{\kap}\, J^{lw}]{\Gam_{lw},\Lam_{lw}}
   {
    \infer*{\Gam_{lw},\lnot C_{lw}}{}
   &
    \infer*{C_{lw},\Lam_{lw}}
    {
     \infer[(c)^{\kap_{4}}\, I_{5}]{\Del_{5}'}{}
     }
    }
   }
  }
 }
{\deduce
 {\hskip-4.7cm
  \infer*{\Del_{5}}
  {
   \infer[K']{\Pi,\Del}
   {
    \infer*[\calc_{5}]{\Pi,\lnot B}{}
   &
    \infer*{B,\Del}
    {
     \infer[(c)_{\kap_{4}}\, I_{4}]{\Del_{4}'}
     {
      \infer*{\Del_{4}}{}
     }
    }
   }
  }
 }
 {\deduce
  {\hskip0.0cm
   \infer[(c)^{\sig_{n_{m(i+1)}+1}}\, I_{2}]{\Del_{2}'}
   {
    \infer*{\Del_{2}}
    {
     \infer[(\Sig_{i_{m(i+1)}})^{\sig_{n_{m(i+1)}+1}}\, K_{m(i+1)}]{\Pi_{m(i+1)},\Del_{m(i+1)}}
     {
      \infer*[\calc_{2}=\calc_{n_{m(i+1)+1}}]{\Pi_{m(i+1)},\lnot B_{m(i+1)}}{}
     &
      \infer*{B_{m(i+1)},\Del_{m(i+1)}}{}
     }
    }
   }
  }
  {\hskip-0.1cm
   \infer*[\calc_{1}]{u:\Phi}{}
   }
 }
}
\]
\hspace*{\fill} $P$
\\
\eprf

By Claim \ref{clm:5M8.2.11}, {\bf (uplwr)} and $i_{m(i+1)}\leq i\leq j$, $K_{m(i+1)}$ foreruns $J^{lw}$ in $P$.
Therefore the right branch $\calt$ of $J^{lw}$ is left to $K_{m(i+1)}$. 
Also by {\bf (h-reg)} $K_{m(i+1)}$ is below $u:\Phi$.
 Hence $\calt$ does not pass through $u:\Phi$ in this case. 
 This shows Claim \ref{clm:5M8.2.1}. 
 In the following figure $\calc_{2}$ denotes the chain starting with $I_{2}$.

\[
 \infer[J^{lw}]{\Gam_{lw},\Lam_{lw}}
 {
  \infer*{\Gam_{lw},\lnot C_{lw}}{}
 &
  \infer*[{}_{K_{m(i+1)}}\calr]{C_{lw},\Lam_{lw}}
  {
   \infer[I_2]{\Del_2'}
   {
    \infer*{\Del_2}
    {
     \infer[(\Sig_{i_{m(i+1)}})^{\sig_{n_{m(i+1)}+1}}\, K_{m(i+1)}]{\Pi_{m(i+1)},\Del_{m(i+1)}}
     {
      \infer*[\calc_{2}]{\Pi_{m(i+1)},\lnot B_{m(i+1)}}{}
     &
      \infer*{B_{m(i+1)},\Del_{m(i+1)}}{u:\Phi}
     }
    }
   }
  }
 }
\]
\hspace*{\fill} $P$
\\
\eprf

\bclm\label{clm:5M8.2.3}
The proviso {\bf (uplw)} holds for the lower rule  $J^{lw}=I_i$ in $P'$.
\eclm
{\bf Proof} of Claim \ref{clm:5M8.2.3}. 
Let $J^{up}$ be a $j$-knot $(\Sig_j)$ above $I_i$. 
Let $H_0$ denote the lower rule of $J^{up}$. 
Assume that the left rope ${}_{J^{up}}\calr=H_0,\ldots,H_{k-1}$ of $J^{up}$ reaches to the rule $I_i$. 
We show
\[
i<j
\]
 even if $J^{up}$ is in the right upper part of $I_i$. 
 Consider the corresponding rule $J^{up}$ in $P$. 
\\
{\bf Case 1} Either $J^{up}$ is $I$ or between $I$ and $u:\Phi$: 
If either $J^{up}$ is $I$ or an $i_m$-knot $K_{m}$ with $m<m(i+1)$, then $i<i+1=j$ or $i<i_m=j$ by (\ref{eqn:m(i+1)}), resp.

 Otherwise $J^{up}$ is between $K_{m-1}$ and $J_{n_m}$ for some $m$ with $0\leq m\leq m(i+1)$. 
Then the rule $J^{up}$ is the merging rule of the chain $\calc_{n_m}$ starting with $J_{n_m}$ and 
the chain $\calc_{H_0}$ starting with $H_0$ so that $\calc_{n_m}$ passes through the right side of $J^{up}$ 
and $\calc_{H_0}$ the left side of $J^{up}$. 
Hence by {\bf (ch:link) Type3 (merge)} the rule $H_{k-1}$ is above $J_{n_m}$ and 
the left rope ${}_{H_0}\calr$ does not reach to $I_i$. 
Thus this is not the case.

\[
\deduce
{\hskip0.0cm
 \infer[J_{n_{m}}]{\Gam_{n_{m}}'}
 {
  \infer*{\Gam_{n_{m}}}
  {
   \infer[H_{k-1}]{\Lam_{k-1}'}
   {
    \infer*{\Lam_{k-1}}
    {}
   }
  }
 }
}
{
 \deduce
 {\hskip-0.3cm
  \infer[H_q]{\Lam_q'}
  {
   \infer*{\Lam_q}
   {}
  }
 }
 {
 \deduce
  {\hskip0.0cm
   \infer[J^{up}]{\Del,\Psi}
   {
    \infer*[\calc_{H_q}]{\Del,\lnot C}{}
   &
    \infer*[\calc_{n_m}]{C,\Psi}
    {}
   }
  }
  {\hskip1.3cm
    \infer[K_{m-1}]{\Pi_{m-1},\Del_{m-1}}
    {
     \infer*[\calc_{n_{m}}]{\Pi_{m-1},\lnot B_{m-1}}{}
    &
     B_{m-1},\Del_{m-1}
    }
  }
 }
}
\]
where $H_q$ denotes the lowermost member of ${}_{H_0}\calr$ such that the chain $\calc_{H_q}$ starting with $H_q$ passes through the left side of $J^{up}$. By {\bf (ch:link) Type3 (merge)} the rule $H_{q}$ is above $J_{n_m}$ and so on.
\\       

\noindent
{\bf Case 2} $J^{up}$ is in the right upper part of $I$: 
Then the left rope ${}_{H_0}\calr$ reaches to $I$. 
Hence by Lemma \ref{lem:5.3.25}, i.e., by {\bf (uplwr)} we have $i<i+1<j$.
\\

\noindent
{\bf Case 3} $J^{up}$ is in the left upper part of $I$: 
Then the left rope ${}_{H_0}\calr$ reaches to $I$. 
Hence by {\bf (uplwl)} we have $i<i+1<j$.
\\

\noindent
{\bf Case 4} 
Otherwise: 
Then there exists a rule $K$ such that $J^{up}$ is in the left upper part of $K$ and $K$ is between $I$ and $\Phi$. 
By {\bf (h-reg), (ch:pass)} $K$ is a rule $(\Sig_{p})^{\kap}$. 
The left rope ${}_{H_0}\calr=H_0,\ldots,H_{k-1}$ reaches to $K$. 
Hence by {\bf (uplwl)} we have
\beqn\label{eqn:clm:5M8.2.3}
p<j
\eeqn
\[
\infer*{\Phi}
{
 \infer[(\Sig_{p})^{\kap}\, K]{\Gam_{K},\Lam_{K}}
 {
  \infer*{\Gam_{K},\lnot D}
  {
   \infer[J^{up}]{\Del,\Psi}
   {
    \infer*{\Del,\lnot C}{}
   &
    \infer*{C,\Psi}{}
   }
  }
 &
  \infer*{D,\Lam_{K}}
  {
   \infer[(\Sig_{i+1}^{\sig})\, I]{\Gam,\Lam}
    {
     \infer*{\Gam,\lnot A_{i+1}^{\sig}}{}
    &
     \infer*{A_{i+1}^{\sig},\Lam}{}
    }
   }
  }
}
\]
{\bf Case 4.1} $H_{k-1}$ is below $K$: 
Let $K'$ denote the uppermost knot such that $K'$ is equal to or below $K$, 
and there exists a member of ${}_{H_0}\calr$ such that the chain starting with the member 
passes through the left side of $K'$. 
Let $H_q$ be the lowermost member of ${}_{H_0}\calr$ such that the chain $\calc_{H_q}$ 
starting with $H_q$ passes through the left side of $K'$. 
If there exists a member of ${}_{H_0}\calr$ such that the chain starting with the member 
passes through the left side of $K$, then $K'$ is equal to $K$.
\[
\infer[H_{q}]{\Del_{q}'}
{
 \infer*{\Del_{q}}
 {
  \infer[(\Sig_{p})^{\kap}\, K=K']{\Gam_{K},\Lam_{K}}
  {
   \infer*[\calc_{H_{q}}]{\Gam_{K},\lnot D}{J^{up}}
  &
   \infer*{D,\Lam_{K}}{}
  }
 }
}
\]
Otherwise $K'$ is below $K$ and it is a knot for the left rope ${}_{H_0}\calr$. 
Let $H_{q_{-1}}$ denote the lowermost member of ${}_{H_0}\calr$ above $K$. 
Then $H_{q_{-1}}$ is an upper right rule of the knot $K'$ and $K'$ is a rule $(\Sig_{p'})^{\kap}$ with
\beqn\label{eqn:clm:5M8.2.31}
p'\leq p
\eeqn 
by {\bf (h-reg)}.
\[
\deduce
{\hskip0.0cm
 \infer[H_{q}]{\Del_{q}'}
 {
  \infer*{\Del_{q}}
  {
   \infer[(\Sig_{p'})^{\kap}\, K']{\Gam_{K'},\Lam_{K'}}
   {
    \infer*[\calc_{H_{q}}]{\Gam_{K'},\lnot D'}{}
   &
    \infer*{D',\Lam_{K'}}{}
   }
  }
 }
}
{\hskip0.3cm
 \infer[(\Sig_{p})^{\kap}\, K]{\Gam_{K},\Lam_{K}}
 {
  \infer*{\Gam_{K},\lnot D}
  {
   \infer[H_{q_{-1}}]{\Del_{q_{-1}}'}
   {
    \infer*{\Del_{q_{-1}}}{J^{up}}
   }
  }
 &
  \infer*{D,\Lam_{K}}{I}
 }
}
\]
By Lemma \ref{lem:5.3.16}
the uppermost member of $\calc_{H_q}$ below $K'$ is the lower rule of the knot $K'$. By (\ref{eqn:clm:5M8.2.3}), (\ref{eqn:clm:5M8.2.31}) and {\bf Case 1} it suffices to show that the left rope ${}_{K'}\calr=G_0,\ldots,G_{k_0}$ of $K'$ reaches to $I_i$, i.e., to show the last member $G_{k_0}$ is equal to or below the rule $H_{k-1}$. Then we will have $i<p'\leq p<j$.

Let $G_0=H_{q_0} \, (q_0\geq q)$ denote the lower rule of $K'$ and $G_{k_1}$ the lowermost member of ${}_{K'}\calr$ such that the chain $\calc_{G_{k_1}}$ starting with $G_{k_1}$ passes through the left side of $K'$. Then by {\bf (ch:link)} $G_{k_1}$ is equal to or below $H_q$. 
\\

\noindent
{\bf Case 4.1.1} $G_{k_1}=H_q$: 
Then $G_{k_0}=H_{k-1}$, i.e., $G_{q_{1}-q_{0}}=H_{q_{1}}$ for any $q_1$ with $q_{0}\leq q_{1}<k$.
\[
\deduce
{\hskip1.8cm
 \infer[H_{q_{1}}=G_{q_{1}-q_{0}}]{\Lam_{q_1}'}
 {
  \infer*{\Lam_{q_1}}
  {
   \infer[H_{q+1}=G_{q+1-q_{0}}]{\Lam_{q+1}'}
   {
    \infer*{\Lam_{q+1}}
    {}
    }
   }
  }
 }
 {
 \deduce
  {\hskip-0.5cm
   \infer[K_{1}]{\Gam_{K_{1}},\Lam_{K_{1}}}
   {
    \infer*[\calc_{H_{q_1}}]{\Gam_{K_{1}},\lnot D_{1}}{}
   &
    \infer*{D_{1},\Lam_{K_{1}}}
    {}
    }
   }
   {
    \deduce
    {\hskip2.6cm
     \infer[H_{q}=G_{k_1}]{\Lam_{q}'}
     {
      \infer*{\Lam_{q}}
      {
       \infer[H_{q_0}=G_{0}]{\Lam_{q_0}'}
       {
        \infer*{\Lam_{q_0}}
       {}
       }
      }
     }
    }
   {\hskip0.8cm
     \infer[K']{\Gam_{K'},\Lam_{K'}}
     {
      \infer*[\calc_{H_{q}}]{\Gam_{K'},\lnot D'}{}
     &
      \infer*{D',\Lam_{K'}}
      {}
     }
    }
   }
  }
\]
where $K_{1}$ is a knot of $H_{q+1}=G_{q+1-q_{0}}$ and $H_{q}=G_{k_1}$ with $q+1-q_{0}=k_{1}+1$.
\\

\noindent
{\bf Case 4.1.2} Otherwise: Then by {\bf (ch:link)} $G_{k_1}$ is already below $H_{k-1}$.

\[
\deduce
{\hskip1.6cm
 \infer[G_{k_{1}}]{\Lam_{k_{1}+q_0}'}
 {
  \infer*{\Lam_{k_{1}+q_0}}
  {
   \infer[H_{k-1}=G_{k-1-q_0}]{\Lam_{k-1}'}
   {
    \infer*{\Lam_{k-1}}
   {}
   }
  }
 }
}
{    
\deduce
{\hskip1.7cm
 \infer[H_{q_{1}}=G_{q_{1}-q_{0}}]{\Lam_{q_1}'}
 {
  \infer*{\Lam_{q_1}}
  {
   \infer[H_{q+1}=G_{q+1-q_{0}}]{\Lam_{q+1}'}
   {
    \infer*{\Lam_{q+1}}
    {}
    }
   }
  }
 }
 {
 \deduce
  {\hskip-0.5cm
   \infer[K_{1}]{\Gam_{K_{1}},\Lam_{K_{1}}}
   {
    \infer*[\calc_{H_{q_1}}]{\Gam_{K_{1}},\lnot D_{1}}{}
   &
    \infer*{D_{1},\Lam_{K_{1}}}
    {}
    }
   }
   {
    \deduce
    {\hskip3.0cm
     \infer[H_{q}=G_{q-q_0}]{\Lam_{q}'}
     {
      \infer*{\Lam_{q}}
      {
       \infer[H_{q_0}=G_{0}]{\Lam_{q_0}'}
       {
        \infer*{\Lam_{q_0}}
       {}
       }
      }
     }
    }
   {
   \deduce
    {\hskip1.2cm
     \infer[K']{\Gam_{K'},\Lam_{K'}}
     {
      \infer*{\Gam_{K'},\lnot D'}
      {}
     &
      D',\Lam_{K'}
     }
    }
    {\hskip-0.7cm
       \infer[K_{mg}]{\Gam_{mg},\Lam_{mg}}
       {
        \infer*[\calc_{H_q}]{\Gam_{mg},\lnot D_{mg}}{}
       &
        \infer*[\calc_{G_{k_{1}}}]{D_{mg},\Lam_{mg}}{}
       }
     }
    }
   }
  }
 }
\]
where $K_{mg}$ is a merging rule of the chain $\calc_{H_q}$ starting with $H_q$ and the chain $\calc_{G_{k_{1}}}$ starting with $G_{k_{1}}$. Since the chain $\calc_{H_{q+1}}$ starting with the lower rule $H_{q+1}=G_{q+1-q_{0}}$ of $K_{1}$ passes through the left side of $K_{1}$, $G_{k_{1}}$ is not equal to $H_{q+1}$ and hence is below $H_{q+1}$ and so on.
\\

\noindent
{\bf Case 4.2} $H_{k-1}$ is above $K$: 
Then $H_{k-1}$ is a rule $(c)_{\sig_{n_{m(i+1)}+1}}$ and $K$ is a rule $(\Sig_p)^{\sig_{n_{m(i+1)}+1}}$. 
Let $d:\Gam_{K},\lnot D$ denote an uppersequent of $K$. 
By {\bf (h-reg)} and the definition of the sequent $u:\Phi$ we have 
$\sig_{n_{m(i+1)}+1}+i\leq h(d;P)\leq\sig_{n_{m(i+1)}+1}+p-1$. 
Thus by (\ref{eqn:clm:5M8.2.3}) we get $i\leq p-1<j$.

\[
\deduce
{\hskip0.6cm
 \infer[(\Sig_{i_{m(i+1)}})^{\sig_{n_{m(i+1)}+1}}\, K_{m(i+1)}]{\Pi_{m(i+1)},\Del_{m(i+1)}}
 {
  \Pi_{m(i+1)},\lnot B_{m(i+1)}
 &
  \infer*{B_{m(i+1)},\Del_{m(i+1)}}
  {
   \infer*{\Phi}
   {}
  }
 }
}
{\hskip0.0cm
 \infer[(\Sig_p)^{\sig_{n_{m(i+1)}+1}}\, K]{\Gam_{K},\Lam_{K}}
 {
  \infer*{d:\Gam_{K},\lnot D}
  {
   \infer[(c)_{\sig_{n_{m(i+1)}+1}}\, H_{k-1}]{\Lam_{k-1}'}
   {
    \infer*{\Lam_{k-1}}{J^{up}}
   }
  }
 &
  \infer*{D,\Lam_{K}}
  {
   \infer[(c)^{\sig_{n_{m(i+1)}}}_{\sig_{n_{m(i+1)}+1}}\, J_{n_{m(i+1)}}]{\Gam_{n_{m(i+1)}}'}
    {
     \infer*{\Gam_{n_{m(i+1)}}}{}
    }
   }
 }
}
\]
where the $i_{m(i+1)}$-knot $K_{m(i+1)}$ disappears when $m(i+1)=l$ in (\ref{eqn:m(i+1)}).

This shows Claim \ref{clm:5M8.2.3}.
\eprf
\smlskp
{\bf M8}. $I$ is a $(\Sig_1)^{\sig}$.\\
This is treated as in the case {\bf M8} of \cite{ptpi3}.
\smlskp
Other cases are easy.
\smlskp
This completes a proof of the Main Lemma \ref{mlem:5}.


\begin{thebibliography}{99}
\bibitem{hndWienpiN} T. Arai, A system $(O(\pi_{n}),<)$ of ordinal diagrams, handwritten note, Apr. 1994.

\bibitem{hndptpiN} T. Arai, Proof theory for reflecting ordinals IV: $\Pi_{N}$-reflecting ordinals, handwritten note, May. 1994.

\bibitem{odMahlo} T. Arai, Ordinal diagrams for recursively Mahlo universes, Archive for Mathematical Logic 39 (2000) 353-391.

\bibitem{odpi3} T. Arai, Ordinal diagrams for $\Pi_{3}$-reflection, Journal of Symbolic Logic 65 (2000) 1375-1394.

\bibitem{ptMahlo} T. Arai, Proof theory for theories of ordinals I:recursively Mahlo ordinals, 
Annals of Pure and Applied Logic vol. 122 (2003), pp. 1-85.

\bibitem{ptpi3} T. Arai, Proof theory for theories of ordinals II:$\Pi_{3}$-Reflection, 
Annals of Pure and Applied Logic vol. 129 (2004), pp. 39-92.

\bibitem{Wienpi3d} T. Arai, Wellfoundedness proofs by means of non-monotonic inductive definitions I: 
$\Pi^{0}_{2}$-operators, 
Journal of Symbolic Logic vol. 69 (2004), pp. 830-850.

\bibitem{LMPSPeking} T. Arai, Iterating the recursively Mahlo operations,
in Proceedings of the thirteenth International Congress of Logic Methodology, Philosophy of Science, 
Ed. by C. Glymour, W. Wei and D. Westerstahl,
College Publications, King's College London (2009), pp. 21-35.
arXiv:1005.1987v1[math.LO]

\bibitem{WienpiN} T. Arai, Wellfoundedness proofs by means of non-monotonic inductive definitions II: 
first order operators, to appear in Annals of Pure and Applied Logic.
arXiv:1005.2007v1[math.LO]

\bibitem{Gentzen38}G. Gentzen, Neue Fassung des Widerspruchsfreiheitsbeweises f\"ur die reine Zahlen theorie, Forschungen zur Logik und zur Grundlegung der exakter Wissenschaften, Neue Folge 4 (1938) 19-44. 

\bibitem{Rathjen94} M. Rathjen, Proof theory of reflection, Ann. Pure Appl. Logic 68 (1994), 181-224.

\bibitem{RathjenAFML1}M. Rathjen, An ordinal analysis of stability, Arch. Math. Logic 44(2005), 1-62.

\bibitem{RathjenAFML2}M. Rathjen, An ordinal analysis of parameter free $\pi^{1}_{2}$-comprehension, 
Arch. Math. Logic 44(2005), 263-362.

\bibitem{Richter-Aczel74} W.H. Richter, P. Aczel, Inductive definitions and reflecting properties of admissible ordinals, in: J.E. Fenstad and P.G. Hinman (Eds.), Generalized Recursion Theory, North-Holland, Amsterdam, 1974, pp. 301-381.

\bibitem{Takeuti87} G. Takeuti, Proof Theory, second edition. (North-Holland, Amsterdam, 1987).

\end{thebibliography}
\end{document}